\documentclass[english,11pt]{article}
\usepackage[T1]{fontenc}
\usepackage[latin9]{inputenc}
\usepackage{geometry}
\usepackage{babel}
\usepackage{amsthm}
\usepackage{amsmath}
\usepackage{amssymb}
\usepackage[unicode= true,pdfusetitle,
 bookmarks= true,bookmarksnumbered= false,bookmarksopen= false,
 breaklinks= false,pdfborder= {0 0 1},backref= false,colorlinks= false]
 {hyperref}
\usepackage{breakurl, mathrsfs, appendix}

\usepackage{fullpage}

\usepackage{graphicx}
\def \figpath {./}
\usepackage[firstinits=true, backend=bibtex]{biblatex}
\makeatletter
\theoremstyle{plain}
\newtheorem{thm}{Theorem}
\newtheorem{prop}{Proposition}
\newtheorem{lem}{Lemma}
\newtheorem{cor}{Corollary}
\theoremstyle{definition}
\newtheorem{dfn}{Definition}
\newtheorem{rem}{Remark}

\long\def\comment#1{}

\usepackage{color}

\newcommand{\silent}[1]{}
\makeatother

\newcommand{\Exs}{\ensuremath{\mathbb{E}}}
\newcommand{\defn}{\ensuremath{: \, =  }}
\newcommand{\order}{\ensuremath{\mathcal{O}}}

\newcommand{\real}{\ensuremath{\mathbb{R}}}
\newcommand{\SymMat}[1]{\ensuremath{\mathcal{S}^{#1 \times #1}}}
\newcommand{\I}{\ensuremath{ I }} 
\newcommand{\onevec}{\ensuremath{\mathbf{1}}} 
\newcommand{\OneMat}{\ensuremath{J_{\usedim}}} 

\newcommand{\Term}{\ensuremath{T}}

\DeclareMathOperator{\trace}{trace}
\DeclareMathOperator{\rank}{rank}

\DeclareMathOperator{\sgn}{sign}

\newlength{\widebarargwidth}
\newlength{\widebarargheight}
\newlength{\widebarargdepth}
\DeclareRobustCommand{\widebar}[1]{%
  \settowidth{\widebarargwidth}{\ensuremath{#1}}%
  \settoheight{\widebarargheight}{\ensuremath{#1}}%
  \settodepth{\widebarargdepth}{\ensuremath{#1}}%
  \addtolength{\widebarargwidth}{-0.3\widebarargheight}%
  \addtolength{\widebarargwidth}{-0.3\widebarargdepth}%
  \makebox[0pt][l]{\hspace{0.3\widebarargheight}%
    \hspace{0.3\widebarargdepth}%
    \addtolength{\widebarargheight}{0.3ex}%
    \rule[\widebarargheight]{0.95\widebarargwidth}{0.1ex}}%
  {#1}}

\newcommand{\inprod}[2]{\ensuremath{\langle #1 , \, #2 \rangle}}
\newcommand{\trinprod}[2]{\ensuremath{\langle \!\langle {#1}, \; {#2}
    \rangle \!\rangle}}

\newcommand{\OUTER}[1]{\ensuremath{#1 {\tiny{\otimes}} #1}}

\newcommand{\matsnorm}[2]{|\!|\!| #1 | \! | \!|_{{#2}}}
\newcommand{\vecnorm}[2]{\| #1\|_{#2}}
\newcommand{\abs}[1]{| #1 |}
\newcommand{\twonorm}[1]{\vecnorm{#1}{2}}
\newcommand{\onenorm}[1]{\vecnorm{#1}{1}}
\newcommand{\opnorm}[1]{\ensuremath{\matsnorm{#1}{\mbox{\tiny{op}}}}}
\newcommand{\nucnorm}[1]{\ensuremath{\matsnorm{#1}{\mbox{\tiny{nuc}}}}}
\newcommand{\fronorm}[1]{\ensuremath{\matsnorm{#1}{\mbox{\tiny{F}}}}}
\newcommand{\infnorm}[1]{\ensuremath{\vecnorm{#1}{\infty}}}
\newcommand{\twoinfnorm}[1]{\ensuremath{\vecnorm{#1}{2,\infty}}}
\newcommand{\twoonenorm}[1]{\ensuremath{\vecnorm{#1}{2,1}}}

\newcommand{\Ball}{\ensuremath{\mathbb{B}}}
\newcommand{\Balltwoone}{\ensuremath{ \Ball_{2,1} }}
\newcommand{\Balltwozero}{\ensuremath{ \Ball_{2,0} }}
\newcommand{\Ballfro}{\ensuremath{ \Ball_{F} }}
\newcommand{\Ballstar}[1]{\ensuremath{ \Ball_{2}(#1; \thetastar) }}

\newcommand{\sigmin}[1]{\sigma_{\min}(#1)}

\newcommand{\usedim}{\ensuremath{d}} 
\newcommand{\rdim}{\ensuremath{r}} 
\newcommand{\numobs}{\ensuremath{n}} 

\newcommand{\epsnum}{\ensuremath{\varepsilon_\numobs}}
\newcommand{\kdim}{\ensuremath{k}}

\newcommand{\smallnoisevar}{\ensuremath{ \epsilon }} 
\newcommand{\noisevar}{\ensuremath{E}} 
\newcommand{\noisestd}{\ensuremath{\sigma}} 

\newcommand{\yout}{\ensuremath{ y }}  
\newcommand{\Yout}{\ensuremath{ Y }} 

\newcommand{\Loss}{\ensuremath{\mathcal{L}}}
\newcommand{\EmpLoss}{\ensuremath{\Loss_\numobs}}
\newcommand{\PopLoss}{\ensuremath{\widebar{\Loss}}}

\newcommand{\pltheta}{\ensuremath{F}}
\newcommand{\Atwo}{\ensuremath{F}}

  \newcommand{\PlTheta}{\ensuremath{M}}
  \newcommand{\thetastar}{\ensuremath{\pltheta^*}}
\newcommand{\tdthetastar}{\PISTAR{\pltheta}}

\newcommand{\ThetaStar}{\ensuremath{\PlTheta^*}}

\newcommand{\eclass}{\ensuremath{\mathcal{E}(\ThetaStar)}}
\newcommand{\dist}[1]{\ensuremath{ \textup{d}(#1, \thetastar) }}
\newcommand{\distsqr}[1]{\ensuremath{ \textup{d}^2(#1, \thetastar) }}
\newcommand{\curv}{\alpha} \newcommand{\smoo}{\beta}
\newcommand{\LipCon}{\ensuremath{L}}
\newcommand{\smallrad}{\ensuremath{\rho}}
\newcommand{\inisig}{\ensuremath{ \tau }}
\newcommand{\starnorm}{\ensuremath{ \psi }}
\newcommand{\cond}{\ensuremath{ \kappa }}
\newcommand{\Grad}{\ensuremath{ G }}
\newcommand{\PopGrad}{\ensuremath{ \widebar{\Grad} }}

\renewcommand{\d}{\ensuremath{ \Delta }}
\newcommand{\diff}{\ensuremath{ \Lambda} }
\newcommand{\thetait}[1]{\ensuremath{\pltheta^{#1}}}
\newcommand{\thetastarit}[1]{\ensuremath{\PISTAR{\thetait{#1}}}}
\newcommand{\stepit}[1]{\ensuremath{\eta^{#1}}}

\newcommand{\preprojit}[1]{\ensuremath{\widetilde{\pltheta}^{#1}}}
\newcommand{\diffit}[1]{\ensuremath{ \Lambda^{#1} }}
\newcommand{\Thetait}[1]{\ensuremath{ \PlTheta^{#1} }}

\newcommand{\PlX}{\ensuremath{X}}

\newcommand{\PlZ}{\ensuremath{Z}}

\newcommand{\plu}{\ensuremath{U}}
\newcommand{\plv}{\ensuremath{V}}

\newcommand{\plx}{\ensuremath{x}}
\newcommand{\ply}{\ensuremath{y}}
\newcommand{\plz}{\ensuremath{z}}
\newcommand{\plphi}{\ensuremath{H}}
\newcommand{\plomega}{\ensuremath{G}}

\newcommand{\Cons}{\ensuremath{\mathcal{M}}}

\newcommand{\GenCons}{\ensuremath{\mathcal{F}}}
\newcommand{\GenConsit}[1]{\ensuremath{\mathcal{F}^{}}}

\newcommand{\Proj}{\ensuremath{\Pi}}
\newcommand{\ProjGenCons}{\ensuremath{\Pi_{\GenCons}}}

\newcommand{\Xit}[1]{X^{#1}} 
\newcommand{\Xmap}{\mathfrak{X}_n}
\newcommand{\ripparam}[1]{\delta_{#1}}

\newcommand{\inco}{\ensuremath{\mu}} 
\newcommand{\pobs}{\ensuremath{p}} 
\newcommand{\Obs}{\ensuremath{\Omega}} 
\newcommand{\ProjObs}{\ensuremath{\Proj_{\Obs}}}
\newcommand{\Tspace}{\ensuremath{ \mathcal{T} }}
\newcommand{\ProjTspace}{\ensuremath{ \Proj_{\Tspace} }}

\newcommand{\CovMat}{\ensuremath{\Sigma}}
\newcommand{\SigHat}{\ensuremath{\widehat{\CovMat}_\numobs}}
\newcommand{\Tset}{\ensuremath{R}}
\newcommand{\TsetComp}{\ensuremath{\Tset^c}}
\newcommand{\snr}{\ensuremath{\gamma}}
\newcommand{\Ckset}[1]{\ensuremath{ \mathbb{C}(#1) }}
\newcommand{\noisemat}{\ensuremath{ W }}  

\newcommand{\Shift}{\ensuremath{S}}
\newcommand{\Adj}{\ensuremath{A}}
\newcommand{\csize}{\ensuremath{k}}  
\newcommand{\ShiftBar}{\ensuremath{\widebar{\Shift}}}
\newcommand{\pin}{\ensuremath{p}}
\newcommand{\qout}{\ensuremath{q}}

\newcommand{\obsnr}{\ensuremath{\nu}}
\newcommand{\obupper}[1]{\ensuremath{L_{#1}}}
\newcommand{\oblower}[1]{\ensuremath{\ell_{#1}}}
\newcommand{\Obsind}{\ensuremath{ Q }}

\newcommand{\SStar}{\ensuremath{ S^{*} }}
\newcommand{\ConsS}{\ensuremath{ \mathcal{S} }}
\newcommand{\ProjConsS}{\ensuremath{ \Proj_{\ConsS} }}
\newcommand{\PlS}{\ensuremath{ S }}
\newcommand{\dS}{\ensuremath{ \Delta_{S} }}
\newcommand{\StarSupp}{\ensuremath{ \Phi^* }}

\newcommand{\PISTAR}[1]{\ensuremath{#1_{\pi^*}}}

\newcommand{\LossTil}{\ensuremath{\widetilde{\Loss}_\numobs}}
\newcommand{\GradTil}{\ensuremath{\widetilde{G}}}
\newcommand{\GradTilit}[1]{\ensuremath{\widetilde{G}^{#1}}}
\newcommand{\LossBar}{\ensuremath{\widebar{\Loss}}}

\newcommand{\Event}{\ensuremath{\mathcal{E}}}

\newcommand{\mprob}{\ensuremath{\mathbb{P}}}

\newcommand{\Sstar}{\SStar}

\newcommand{\Rtil}{\ensuremath{\widetilde{R}}}

\newcommand{\frobnorm}[1]{\ensuremath{\fronorm{#1}}}

\newcommand{\MYPSI}{\ensuremath{\psi}}

\newcommand{\vtiny}{\vspace*{.05in}}
\newcommand{\vsmall}{\vspace*{.1in}}

\makeatletter
\long\def\@makecaption#1#2{
        \vskip 0.8ex
        \setbox\@tempboxa\hbox{\small {\bf #1:} #2}
        \parindent 1.5em  
        \dimen0=\hsize
        \advance\dimen0 by -3em
        \ifdim \wd\@tempboxa >\dimen0
                \hbox to \hsize{
                        \parindent 0em
                        \hfil 
                        \parbox{\dimen0}{\def\baselinestretch{0.96}\small
                                {\bf #1.} #2
                                } 
                        \hfil}
        \else \hbox to \hsize{\hfil \box\@tempboxa \hfil}
        \fi
        }
\makeatother

\begin{document}

\begin{center} {\LARGE{\bf{
Fast low-rank estimation by projected gradient descent: \\ General
statistical and algorithmic guarantees}}} \\

\vspace*{.3in}

{\large{
\begin{tabular}{ccc}
Yudong Chen$^\ast$ && Martin J. Wainwright$^{\dagger,\ast}$ 
\end{tabular}

\vspace*{.1in}

\begin{tabular}{ccc}
Department of Statistics$^{\dagger}$
& &  Department of Electrical Engineering and Computer Sciences$^{\ast}$
\end{tabular}
\begin{tabular}{c}
University of California, Berkeley \\ Berkeley, CA 94720
\end{tabular}

\vspace*{.2in}

\begin{tabular}{c}
{\texttt{yudong.chen@cornell.edu $ \quad $ wainwrig@berkeley.edu}}
\end{tabular}
}}

\vspace*{.2in}

\today

\vspace*{.2in}
\begin{abstract}
Optimization problems with rank constraints arise in many
applications, including matrix regression, structured PCA, matrix
completion and matrix decomposition problems.  An attractive heuristic
for solving such problems is to factorize the low-rank matrix, and to
run projected gradient descent on the nonconvex factorized
optimization problem.  The goal of this problem is to provide a
general theoretical framework for understanding when such methods work
well, and to characterize the nature of the resulting fixed point.  We
provide a simple set of conditions under which projected gradient
descent, when given a suitable initialization, converges geometrically
to a statistically useful solution.  Our results are applicable even
when the initial solution is outside any region of local convexity,
and even when the problem is globally concave.  Working in a
non-asymptotic framework, we show that our conditions are satisfied
for a wide range of concrete models, including matrix regression,
structured PCA, matrix completion with real and quantized
observations, matrix decomposition, and graph clustering problems.
Simulation results show excellent agreement with the theoretical
predictions.
\end{abstract}
\end{center}


\section{Introduction}

\label{SecIntro}

There are a variety of problems in statistics and machine learning
that require estimating a matrix that is assumed---or desired---to be
low-rank.  For high-dimensional problems, the low-rank property is is
useful as a form of regularization, and also can lead to more
interpretable results in scientific settings.  Low-rank matrix
estimation can be formulated as a nonconvex optimization problem
involving a cost function, measuring the fit to the data, along with a
rank constraint.  Even when the cost function is convex---such as in
the ubiquitous case of least-squares fitting---solving a
rank-constrained problem can be computationally difficult, with many
interesting special cases known to have NP-hard complexity in the
worst-case setting.  However, statistical settings lead naturally to
random ensembles, in which context such complexity concerns have been
assuaged to some extent by the use of semidefinite programming (SDP)
relaxations.  These SDP relaxations are based on replacing the
nonconvex rank constraint with a convex constraint based on the
trace/nuclear norm.  For many statistical ensembles of problems, among
them multivariate regression, matrix completion and matrix
decomposition, such types of SDP relaxations have been shown to have
near-optimal performace (e.g., see the papers~\cite{candes2009exact,
  recht2010guaranteed, negahban2010noisymc, negahban2009estimation,
  chen2011LSarxiv, koltchinskii2011matrixCompletion} and references
therein).  Although in theory, any SDP can be solved to $\epsilon$
accuracy in polynomial-time~\cite{NesNem87}, the associated
computational cost is often too high in practice.  Letting $\usedim$
denote the dimension of the matrix, it can be as high as $\usedim^{6}$
using standard interior point methods~\cite{Boyd02,NesNem87}; such a
scaling is not practical for many real-world applications involving
high-dimensional matrices.  More recent work has developed algorithms
that are specifically tailored to certain classes of SDPs; however,
even such specialized algorithms require at least $\usedim^{2}$ time,
since solving the SDP involves optimizing over the space of $\usedim
\times \usedim$ matrices.

In practice, researchers often resort instead to heuristic methods
that directly optimize over the space of low-rank matrices, using
iterative algorithms such as alternating minimization, power
iteration, expectation maximization~(EM) and projected gradient
descent. Letting $\rdim$ denote the rank, these factorized
optimization problems live in an $\order(\rdim \usedim)$ dimensional
space, as opposed to the $\order(\usedim^2)$ space of the original
problem.  Such heuristic methods are quite effective in practice for
some problems, but sometimes can also suffer from local optima.  These
intriguing phenomena motivate a recent and evolving line of work on
understanding such iterative methods in the low-rank space.  As we
discuss in detail below, recent work has studied some of these
algorithms in a number of specific settings.  A natural question then
arises: is there a general theory for understanding when low-rank
iterative methods will succeed?

In this paper, we make progress on this general question by focusing
on projected gradient descent in the low-rank space.  We characterize
a general set of conditions that govern the computational and
statistical properties of the solutions, and then specialize this
general theory to obtain corollaries for a broad range of problems.
In more detail, suppose that we write a rank-$\rdim$ matrix $\PlTheta
\in \real^{\usedim \times\usedim}$ in its factorized form
$\OUTER{\pltheta} = \pltheta \pltheta^\top$, where $\pltheta \in
\real^{\usedim \times \rdim}$, and consider projected gradient descent
methods in the variable $\pltheta$.  The matrix quadratic form
$\OUTER{\pltheta}$ makes the problem inherently nonconvex, and in many
cases, the problem is not even locally convex. Nevertheless, our
theory shows that given a suitable initialization, projected gradient
descent converges geometrically to a statistically useful solution,
under conditions that are much more general than convexity. Our
results are applicable even when the initial solution is outside any
region of local convexity, or when the problem is globally
concave. Each iteration of projected gradient descent typically takes
time that is linear in $\usedim \rdim$, the degrees of freedom of a
low-rank matrix, as well as in the input size.  Therefore, by directly
enforcing low-rankness, our method simultaneously achieves two goals:
we not only attain \emph{statistical} consistency in the
high-dimensional regime, but also gain \emph{computational} advantages
over convex relaxation methods that lift the problem to the space of
$\usedim \times\usedim$ matrices.

For this approach to be relevant, an equally important question is
when the above conditions for convergence are satisfied. We verify
these conditions for a broad range of statistical and machine learning
problems, including matrix sensing, matrix completion in both its
standard and one-bit forms, sparse principal component analysis
(SPCA), graph clustering, and matrix decomposition or robust PCA. For
each of these problems, we show that a suitable initialization can be
obtained efficiently using simple methods, and the projected gradient
descent approach has sample complexity and statistical error bounds
that are comparable (and sometimes better) to the best existing
results (which are often achieved by convex relaxation
methods). Notably, our approach does not require using fresh samples
in each iteration---a heuristic known as sample splitting that is
often used to simplify analysis---nor does it involve the computation
of multiple singular value decompositions (SVDs).

Let us now put our contributions in a broader historical context.  The
seminal work in~\cite{burer2005LRSDP} studies the problem of obtaining
low-rank solutions to SDPs using gradient descent on the factor
space. Several subsequent papers aim to obtain rigorous guarantees for
nonconvex gradient descent focused on specific classes of matrix
estimation problems. For instance, the recent
papers~\cite{zheng2015convergent,tu2015procrustes} study exact
recovery in the setting of noiseless matrix sensing (i.e., solving
linear matrix equalities with  random designs).  Focusing on the
rank-one setting, De et al.~\cite{de2014global,de2015wild} study the
noiseless matrix completion problem, and a stochastic version of
nonconvex gradient descent; they prove global convergence with a constant success probability, assuming
independence between the samples used by each iteration. The recent
manuscript~\cite{sun2014mc} studies several variants of nonconvex
gradient descent algorithms, again for noiseless matrix
completion. Another line of
work~\cite{candes2014wirtinger,chen2015quadratic} considers the phase
retrieval problem, which can be reformulated as recovering a rank one
($\rdim = 1$) matrix from random quadratic measurements. The
regularity conditions imposed in this work bear some similarity with
our conditions, but their validation requires a very different
analysis.  An attractive feature of phase retrieval is that it is known to be locally convex around the global optimum under certain
settings~\cite{SoltanolkotabiThesis,white2015local_convex}.

The work in this paper develops a unified framework for analyzing the
behavior of projected gradient descent in application to low-rank
estimation problems, covering many of the models described above as
well as various others.  Our theory applies to matrices of arbitrary
rank $r$, and is framed in the statistical setting of noisy
observations, allowing for noiseless observations as a special case,
When specialized to particular models, our framework yields a variety
of corollaries providing guarantees for concrete statistical models
that have not been studied in the work above.  Notably, our general
conditions \emph{do not} depend on local convexity, and thus can be
applied to models such as sparse PCA and clustering in which no form
of local convexity holds.  (In fact, our results apply even when the
loss function is globally concave).  In addition, we impose only a
natural gradient smoothness condition that is much less restrictive
than the vanishing gradient condition imposed in other work.  Thus,
one of the main contributions of this paper is to illuminate to
weakest known conditions under which nonconvex gradient descent can
succeed, and also allows for applications to several problems that
lack local convexity and vanishing gradients.

It is also worth noting that other types of algorithms for nonconvex
problems have also been analyzed, including alternating
minimization~\cite{jain2013low,hardt2014understanding,hardt2014fast},
EM algorithms~\cite{balakrishnan2014EM,wang2014HighDimEM} and power
methods~\cite{hardt2014power}, various hard-thresholding and singular
value
projection~\cite{jain2010svp,netrapalli2014nonconvexRPCA,jain2014hard,bhatia2015rob_regression},
gradient descent for nonconvex regression and spectrally sparse recovery problems~\cite{loh2013regularized,wang2013sparseNonConvex,cai2015hankel}, as well as
gradient descent on Grassmannian
manifolds~\cite{keshavan2009matrixafew,zhang2015grassmannian}.
Finally, there is a large body of work on convex-optimization based
approach to the concrete examples considered in this paper. We compare
our statistical guarantees with results of these types after the
statements of each of our corollaries.

\paragraph{Notation:} The $i$-th row and $j$-th column of a matrix
$\PlZ$ are denoted by $\PlZ_{i\cdot}$ and $\PlZ_{\cdot j}$,
respectively. The spectral norm $\opnorm{\PlZ}$ is the largest
singular value of $ \PlZ$. The nuclear norm $\nucnorm{\PlZ}$ is the
sum of the singular values of $\PlZ$.  For parameters $1 \le a,b\le
\infty$ and a matrix $\PlZ$, the $\ell_{a}/\ell_{b}$ norm of $ \PlZ $
is $\matsnorm{\PlZ}{b,a} = \big( \sum_{i} \vecnorm{\PlZ_{i\cdot}}b^{a}
\big)^{\frac{1}{a}}$---that is, the $\ell_{a}$ norm of the vector of
the $\ell_{b}$ norms of the rows. Special cases include the Frobenius
norm $ \fronorm{\PlZ} = \matsnorm{\PlZ}{2,2}$, the elementwise
$\ell_{1}$ norm $\onenorm{\PlZ} = \matsnorm{\PlZ}{1,1}$ and the
elementwise $ \ell_{\infty} $ norm $ \infnorm{\PlZ} =
\matsnorm{\PlZ}{\infty, \infty}$.  For a convex set $T$, we use
$\Proj_{T}$ to denote the Euclidean projection onto $T$.


\section{Background}

We begin by setting up the class of matrix estimators to be studied in
this paper, and then providing various concrete examples of specific
models to which our general theory applies.

\subsection{Matrix estimators in the factorized formulation}

Letting $\SymMat{\usedim}$ denote the space of all symmetric
$\usedim$-dimensional matrices, this paper focuses on a class of
matrix estimators that take the following general form.  For a given
sample size $\numobs \geq 1$, let $\EmpLoss: \SymMat{\usedim}
\rightarrow \real$ be a cost function. It is a random function, since
it depends (implicitly in our notation) on the observed data, and the
function value $\EmpLoss(\PlTheta)$ provides some measure of fit of
the matrix $\PlTheta$ to the given data.  For a given convex set
$\Cons \subseteq \SymMat{\usedim}$, we then consider a minimization
problem of the form
\begin{align}
\label{EqnEmpiricalSDP}
\min_{\PlTheta \in \SymMat{\usedim}} \EmpLoss(\PlTheta) \qquad
\mbox{such that \ensuremath{\PlTheta\succeq0} and \ensuremath{\PlTheta
    \in \Cons}.}
\end{align}
The goal of solving this optimization problem is to estimate some
unknown target matrix $\ThetaStar$.  Typically, the target matrix is a
(near)-minimizer of the population version of the program---that is, a
solution to the same constrained minimization problem with $\EmpLoss$
replaced by its expectation $\PopLoss(\PlTheta) =
\Exs[\EmpLoss(\PlTheta)]$.  However, our theory does not require that
$\ThetaStar$ minimizes this quantity, nor that the gradient $\nabla
\PopLoss(\ThetaStar)$ vanish.

In many cases, the matrix $\ThetaStar$ either has low rank, or can be
well-approximated by a matrix of low rank.  Concretely, if the target
matrix $\ThetaStar$ has rank $\rdim < \usedim$, then it can be written
in the outer product form $\ThetaStar= \OUTER{\thetastar}$ for some
other matrix $\thetastar \in \real^{\usedim \times \rdim}$ with orthogonal columns. This
factorized representation motivates us to consider the function
$\LossTil(\pltheta) \defn \EmpLoss(\OUTER{\pltheta})$, and the
factorized formulation
\begin{align}
\label{EqnEmpiricalFactorized}
\min_{\pltheta \in \real^{\usedim \times \rdim}} \quad &
\LossTil(\pltheta) \qquad \mbox{such that \ensuremath{\pltheta \in
    \GenCons},}
\end{align}
where $\GenCons$ is some convex set that contains $\thetastar$, and
for which the set $\big \{ \OUTER{\pltheta} \, \mid \, \pltheta \in
\GenCons \}$ acts as a surrogate for $\Cons$.  Note that due to the
factorized representation of the low-rank matrix, this factorized
program is (in general) nonconvex, and is typically so even if the
original program~\eqref{EqnEmpiricalSDP} is convex.

Nonetheless, we can apply a projected gradient descent method in order
to compute an approximate minimizer. For this particular problem, the projected gradient
descent updates take the form
\begin{align}
\label{EqnGenProjGradOpt}
\thetait{t+1} & = \ProjGenCons\Big(\thetait{t}-\stepit{t}\nabla
\LossTil( \thetait{t} ) \Big)
\end{align}
where $\stepit{t} >0$ is a step size parameter, $\ProjGenCons$ denotes
the Euclidean projection onto the set $\GenCons$, and the
gradient\footnote{This gradient takes the simpler form $ \nabla
  \LossTil(\pltheta) = 2 \nabla_{\PlTheta} \EmpLoss(\OUTER{\pltheta})
  \pltheta $ whenever $ \nabla \EmpLoss (\OUTER{\pltheta}) $ is
  symmetric, which is the case in the concrete examples that we
  treat.}  is given by $\nabla \LossTil(\pltheta) = \big[
  \nabla_{\PlTheta} \EmpLoss(\OUTER{\pltheta}) + (\nabla_{\PlTheta}
  \EmpLoss(\OUTER{\pltheta}) )^\top \big]\pltheta$.  The main goal of
this paper is to provide a general set of sufficient conditions under
which---up to a statistical tolerance term $\epsnum$---the sequence
$\{\thetait{t}\}_{t=0}^\infty$ converges to some $\thetastar$ such
that $\OUTER{\thetastar} = \ThetaStar$.

A significant challenge in the analysis is the fact that there are
\emph{many} possible factorizations of the form $\ThetaStar=
\OUTER{\thetastar}$.  In order to address this issue, it is convenient
to define an equivalent class of valid solutions as follows
\begin{align}
\label{EqnDefnEclass}
\eclass & \defn \{\thetastar \in \real^{\usedim \times \rdim} \,\mid\,
\OUTER{\thetastar} = \ThetaStar, \, {\thetastar_{\cdot i}}^\top
\thetastar_{\cdot j} =0, \forall i\neq j \, \}.
\end{align}
For the applications of interest here, the underlying goal is to
obtain a good estimate of \emph{any} matrix in the set ~$\eclass$. In
particular, such an estimate implies a good estimate of~$\ThetaStar$
itself as well as the column space and singular values of all the
members of the class $\eclass$.  Accordingly, we define the
pseudometric
\begin{align}
\label{EqnDefnDist}
\dist{\pltheta} \defn \min_{\thetastar \in
  \eclass}\fronorm{\pltheta-\thetastar}.
\end{align}
Note that all matrices $\thetastar \in \eclass$ have the same singular
values, so that we may write the singular values
$\sigma_{1}(\thetastar) \geq \cdots \geq \sigma_{\rdim} (\thetastar) >
0 $ as well as $\opnorm{\thetastar}$ and $\fronorm{\thetastar}$
without any ambiguity.  In fact, this invariant property holds more
generally for any function of the sorted singular values and column
space of $ \thetastar $ (e.g., any unitarially invariant norm).


\subsection{Illustrative examples}
\label{SecIllustrative}
Let us now consider a few specific models to illustrate the general
set-up from the previous section.  We return to demonstrate
consequences of our general theory for these (and other) models in
Section~\ref{SecConcrete}.

\subsubsection{Matrix regression}
\label{SecMatrixSensingExample}

We begin with a simple example, namely one in which we make noisy
observations of linear projections of an unknown low-rank matrix
$\ThetaStar \in \SymMat{\usedim}$.  In particular, suppose that we are
given $\numobs$ i.i.d. observations $\{(y_i, X_i)\}_{i=1}^\numobs$ of
the form
\begin{align}
\label{EqnMatrixSensing}
y_i & = \trace(X^i \ThetaStar) + \smallnoisevar_i \quad \mbox{for $i =
  1, \ldots, \numobs$,}
\end{align}
and $\{\smallnoisevar_i\}_{i=1}^\numobs$ is some i.i.d. sequence of
zero-mean noise variables.  The paper~\cite{negahban2009estimation}
provides various examples of such matrix regression problems,
depending on the particular choice of the regression matrices
$\{\Xit{i}\}_{i=1}^\numobs$.

\paragraph{Original estimator:}
Without considering computational complexity, a reasonable estimate of
$\ThetaStar$ would be based on minimizing the least-squares cost
\begin{align}
\label{EqnLeastSquares}
\EmpLoss(\PlTheta) & \defn \frac{1}{2 \numobs} \sum_{i=1}^\numobs
\big(y_i - \trace(X^i \PlTheta) \big)^2
\end{align}
subject to a rank constraint.  However, this problem is
computationally intractable in general due to the nonconvexity of the rank function.  A
standard convex relaxation is based on the nuclear norm
\mbox{$\nucnorm{\PlTheta} \defn \sum_{j=1}^\usedim
  \sigma_j(\PlTheta)$,} corresponding to the sum of the singular
values of the matrix.  In the symmetric PSD case, it is equivalent to
the trace of the matrix.  Using the nuclear norm as regularizer leads
to the estimator
\begin{align*}
\min_{\PlTheta \in \SymMat{\usedim}} \Big \{ \frac{1}{2 \numobs}
\sum_{i=1}^\numobs \big(y_i - \trace(X^i \PlTheta) \big)^2 \Big \}
\qquad \mbox{such that $\PlTheta \succeq 0$ and $\nucnorm{\PlTheta}
  \leq R$},
\end{align*}
where $R > 0$ is a radius to be chosen.  This is a special case of our
general estimator~\eqref{EqnEmpiricalSDP} with $\EmpLoss$ being the
least-squares cost~\eqref{EqnLeastSquares}, and the constraint set
$\Cons = \{ \PlTheta \in \SymMat{\usedim} \, \mid \,
\nucnorm{\PlTheta} \leq R \}$.

\paragraph{Population version:}   Suppose that the 
noise variables $\smallnoisevar_i$ are i.i.d. zero-mean with variance
$\sigma^2$, and the regression matrices $\{X_i\}_{i=1}^\numobs$ are
also i.i.d., zero-mean and such that $\Exs[\trace(X_i \PlTheta)^2] =
\fronorm{\PlTheta}^2$ for any matrix $\PlTheta$. Under these
conditions, an easy calculation yields that the population cost
function is given by $\LossBar(\PlTheta) = \frac{1}{2}
\fronorm{\PlTheta - \ThetaStar}^2 + \frac{1}{2} \sigma^2$.  For this
particular case, note that $\ThetaStar$ is the unique minimizer of the
population cost.

\paragraph{Projected gradient descent:}  The factorized cost function
is given by
\begin{align}
\label{EqnMatrixSensingLossTil}
\LossTil(\pltheta) & = \frac{1}{2 \numobs} \sum_{i=1}^\numobs \Big \{
y_i - \trace \big( \Xit{i} (\OUTER{\pltheta}) \Big \}^2,
\end{align}
and has gradient $\nabla \LossTil(\pltheta) = \frac{2}{\numobs}
\sum_{i=1}^\numobs ( y_i - \trace \big( \Xit{i} (\OUTER{\pltheta} )
\big) (\Xit{i})^T \pltheta$ assuming  each $ \Xit{i} $ is symmetric.  Setting $\GenCons = \real^{\usedim \times \rdim}$, the projected gradient descent
updates~\eqref{EqnGenProjGradOpt} reduce to usual gradient
descent---that is,
\begin{align*}
\thetait{t+1} & = \thetait{t}- \stepit{t} \nabla
\LossTil(\thetait{t}), \qquad \mbox{for $t = 0, 1, \ldots$.}
\end{align*}
We return to analyze these updates in Section~\ref{SecMatrixSensing}.


\subsubsection{Rank-$\protect\rdim$ PCA with row sparsity}
\label{SecSparsePCAExample}

Principal component analysis is a widely used method for
dimensionality reduction.  For high-dimensional problems in which
$\usedim \gg \numobs$, it is well-known that classical PCA is
inconsistent~\cite{JohLu09}. Moreover, minimax lower bounds show that
consistent eigen-estimation is impossible in the absence of structure
in the eigenvectors.  Accordingly, a recent line of work
(e.g.,~\cite{JohLu09,AmiWai09,cai2013adaptive,berthet2013lowerSparsePCA,vu2013minimax,BirJohNadPau12})
has studied different forms of PCA with structured eigenvectors.

Here we consider one such form of structured PCA, namely a rank
$\rdim$ model with row-wise sparsity.  For a given signal-to-noise
ratio $\snr >0$ and an orthonormal matrix $\thetastar \in
\real^{\usedim\times \rdim}$, consider a covariance matrix of the form
\begin{align}
 \label{EqnSpikedCovariance}
\CovMat & = \snr \underbrace{\big(\OUTER{\thetastar}
  \big)}_{\ThetaStar} + \I_{\usedim}.
\end{align}
By construction, the columns of $\thetastar$ span the top rank-$\rdim$
eigenspace of $\CovMat$ with the corresponding maximal eigenvalues
$\snr + 1$. In the row-sparse version of this
model~\cite{vu2013minimax}, this leading eigenspace is assumed to be
supported on $\kdim$ coordinates---that is, the matrix $\thetastar$
has at most $\kdim$ non-zero rows. Given $\numobs$ i.i.d. samples
$\{x_{i}\}_{i= 1}^{\numobs}$ from the Gaussian distribution
$N(0,\CovMat)$, the goal of sparse PCA is to estimate the sparse
eigenspace spanned by~$\thetastar$.

\paragraph{Original estimator:}
A natural estimator is based on a semidefinite program, referred to as
the Fantope relaxation in the paper~\cite{vu2013fantope}, given by
\begin{align}
\label{EqnFantopeSDP}
\min_{ \substack{\PlTheta \in \SymMat{\usedim} \\ 0 \preceq \PlTheta
    \preceq \I_{\usedim} } } \Big\{ - \trace(\SigHat\PlTheta) \Big\} &
\qquad \mbox{such that \ensuremath{\trace(\PlTheta)\leq\rdim} and
  \ensuremath{\onenorm{\PlTheta} \leq R}},
\end{align}
where $\SigHat$ is the empirical covariance matrix, and $R > 0$ is a
radius to be chosen. This is a special case of our general set-up with
$\EmpLoss(\PlTheta) = -\trace(\SigHat \PlTheta)$ and
\begin{align*}
\Cons & \defn \Big \{ \PlTheta \in \SymMat{\usedim} \, \mid \, 0
\preceq \PlTheta \preceq \I_\usedim, \mbox{
  \ensuremath{\trace(\PlTheta)\leq\rdim} and
  \ensuremath{\onenorm{\PlTheta} \leq R}} \; \Big \}.
\end{align*}

\paragraph{Population version:}
Since $\Exs[\SigHat] = \CovMat$, the population cost function is given
by
\begin{align*}\LossBar(\PlTheta) = \Exs[ \EmpLoss(\PlTheta)] = - \trace(\CovMat
\PlTheta).
\end{align*}
Thus, by construction, for any radius $R \geq
\onenorm{\OUTER{\thetastar}}$, the matrix $\ThetaStar=
\OUTER{\thetastar}$ is the unique minimizer of the population version
of the problem~\eqref{EqnFantopeSDP}, subject to the constraint
$\PlTheta \in \Cons$.

\paragraph{Projected gradient descent:}
For a radius $\Rtil$ to be chosen, we consider a factorized version of
the SDP
\begin{align}
\label{EqnLossTilPCA}
\LossTil(\pltheta) \defn - \trinprod{\SigHat}{\OUTER{\pltheta}},
\qquad \GenCons \defn \big\{\pltheta \in \real^{\usedim\times \rdim}
\, \mid \, \opnorm{\pltheta} \leq 1, \, \twoonenorm{\pltheta} \leq
\Rtil \big\},
\end{align}
where we recall that $\twoonenorm{\pltheta} = \sum_{i= 1}^{\usedim}
\twonorm{\pltheta_{i\cdot}}$.  This norm is the appropriate choice for
selecting matrices with sparse rows, as assumed in our initial set-up.
We return in Section~\ref{SecSparsePCA} to analyze the
projected gradient updates~\eqref{EqnGenProjGradOpt} applied to pair
$(\LossTil, \GenCons)$ in equation~\eqref{EqnLossTilPCA}.

As a side-comment, this example illustrates that our theory does not
depend on local convexity of the function $\LossTil$.  In this case,
even though the original function $\EmpLoss$ is convex (in fact,
linear) in the matrix $\PlTheta \in \SymMat{\usedim}$, observe that
that the function $\LossTil$ from equation~\eqref{EqnLossTilPCA} is
never locally convex in the low-rank matrix $\pltheta \in
\real^{\usedim \times \rdim}$; in fact, since $\SigHat$ is positive
semidefinite, it is a globally concave function.



\subsubsection{Low-rank and sparse matrix decomposition}
\label{SecMatrixDecompositionExample}
There are various applications in which it is natural to model an
unknown matrix as the sum of two matrices, one of which is low-rank
and the other of which is sparse.  Concretely, suppose that we make
observations of the form $\Yout = \ThetaStar + \Sstar + \noisevar$
where $\ThetaStar$ is low-rank, the matrix $\Sstar$ is symmetric and 
elementwise-sparse, and $\noisevar$ is a symmetric matrix of noise variables.
Many problems can be cast in this form, including robust forms of PCA,
factor analysis, and Gaussian graphical model estimation; see the
papers~\cite{candes2009robustPCA,agarwal2012decomposition,chen2011LSarxiv,chandrasekaran2011siam,Hsu2010RobustDecomposition} and references
therein for further details on these and other applications.

\paragraph{Original estimator:}
Letting $\PlS_j \in \real^\usedim$ denote the $j^{th}$ column of a
matrix $\PlS \in \real^{\usedim \times \usedim}$, define the set of
matrices $\ConsS \defn \{ \PlS \in \real^{\usedim \times \usedim} \,
\mid \, \onenorm{\PlS_{j}} \le R_j \quad \mbox{for $j =
  1,2,\ldots,\usedim$} \}$, where $(R_1, \ldots, R_\usedim)$ are
user-defined radii. Using the nuclear norm and $\ell_{1}$ norm as
surrogates for rank and sparsity respectively, a popular convex
relaxation approach is based on the SDP
\begin{align*}
\min_{\PlTheta \in \SymMat{\usedim}} \Big \{ \frac{1}{2} \min_{\PlS
  \in \ConsS} \fronorm{\Yout - (\PlTheta + \PlS)}^{2} \Big \} \qquad
\mbox{subject to $\PlTheta \succeq 0$ and $\nucnorm{\PlTheta} \leq
  R$,}
\end{align*}
This is a special case of our general estimator with
$\EmpLoss(\PlTheta) \defn \frac{1}{2} \min \limits_{\PlS \in \ConsS}
\fronorm{\Yout - (\PlTheta + \PlS)}^{2}$, and the constraint set
$\Cons \defn \{ \PlTheta \in \SymMat{\usedim} \, \mid \,
\nucnorm{\PlTheta} \leq R \}$.

\paragraph{Population version:}
In this case, the population function is given by
\begin{align*}
\LossBar(\PlTheta) & \defn \Exs \Big[ \frac{1}{2} \min \limits_{\PlS
    \in \ConsS} \fronorm{\Yout - (\PlTheta + \PlS)}^{2} \Big \} \Big],
\end{align*}
where the expectation is over the random noise matrix $\noisevar$.  In
general, we are not guaranteed that $\ThetaStar$ is the unique
minimizer of this objective, but our analysis shows that (under
suitable conditions) it is a near-minimizer, and this is adequate for
our theory.

\paragraph{Projected gradient descent:}
In this paper, we analyze a version of gradient descent that operates
on the pair $(\LossTil, \GenCons)$ given by
\begin{align}
\label{EqnMatrixDecompositionPD}
\LossTil(\pltheta) = \frac{1}{2} \min_{\PlS \in \ConsS} \fronorm{\Yout
  - ( (\OUTER{\pltheta}) + \PlS)}^{2}, \quad \mbox{and} \quad \GenCons
\defn \Big\{ \pltheta \in \real^{\usedim \times \rdim} \,\mid\,
\twoinfnorm{\pltheta} \le \sqrt{\frac{2\inco}{\usedim}}
\fronorm{\thetait{0}} \Big\}.
\end{align}
Here $\thetait{0}$ is the initialization of the algorithm, and the
parameter $\inco > 0$ controls the matrix incoherence.	See
Sections~\ref{SecMC} and~\ref{SecMatrixDecomposition} for discussion
of matrix incoherence parameters, and their necessity in such
problems.  The gradient of $\LossTil$ takes the form
\begin{align*}
\nabla\LossTil(\pltheta) & = 2 \Big \{ \ProjConsS \big(\Yout -
(\OUTER{\pltheta}) \big) - \big(\Yout - (\OUTER{\pltheta}) \big) \Big
\} \pltheta,
\end{align*}
where $\ProjConsS$ denotes projection onto the constraint set
$\ConsS$.  This projection is easy to carry it, as it simply involves
a soft-thresholding of the columns of the matrix.  Likewise, the
projection onto the set $\GenCons$ from
equation~\eqref{EqnMatrixDecompositionPD} is easy to carry out.  We
return to analyze these projected gradient updates in
Section~\ref{SecMatrixDecomposition}. \\

\vspace*{.2in}

In addition to the three examples introduced so far, our theory also
applies to various other low-rank estimation problems, including that
of matrix completion with real-valued observations
(Section~\ref{SecMC}) and binary observations (Section~\ref{SecOB}),
as well as planted clustering problems (Section~\ref{SecClustering}).

\section{Main results}
\label{SecMain}

In this section, we turn to the set-up and statement of our main
results on the convergence properties of projected gradient descent
for low-rank factorizations.  We begin in Section~\ref{SecConditions}
by stating the conditions on the function $\LossTil$ and $\GenCons$
that underlie our analysis.  In Section~\ref{SecSublinear}, we state a
result (Theorem~\ref{ThmGeneralLip}) that guarantees sublinear
convergence, whereas Section~\ref{SecLinear} is devoted to a result (Theorem~\ref{ThmGeneralSmooth}) 
that guarantees faster linear convergence under slightly stronger
assumptions.  In Section~\ref{SecConcrete} to follow, we derive
various corollaries of these theorems for different concrete versions
of low-rank estimation.

Given a radius $\smallrad >0$, we define the ball
$\Ballstar{\smallrad} \defn \big\{\pltheta \in \real^{\usedim \times
  \rdim} \, \mid \,\dist{\pltheta} \leq \smallrad \big\}$.  At a high
level, our goal is to provide conditions under which the projected
gradient sequence $\{\thetait{t}\}_{t=0}^\infty$ converges some
multiple of the ball $\Ballstar{ \epsnum}$, where $\epsnum > 0$ is a
\emph{statistical tolerance}.


\subsection{Conditions on the pair $(\LossTil, \GenCons)$}
\label{SecConditions}

Recall the definition of the set $\eclass $ of equivalent orthogonal
factorizations of a given matrix $\ThetaStar$.	We begin with a
condition on $\GenCons$ that guarantees that it respects the structure
of this set.

\paragraph{$\ThetaStar$-faithfulness of $\GenCons$:}
For a radius $\smallrad$, the constraint set $\GenCons$ is said to be
\emph{$\ThetaStar$-faithful} if for each matrix $\pltheta \in \GenCons
\cap \Ball_{2}(\smallrad;\thetastar)$, we guaranteed that
\begin{align}
\label{EqnFaithful}
\arg \min_{A \in \eclass}\fronorm{A - \pltheta} & \subseteq \GenCons.
\end{align}
Of course, this condition is implied by the inclusion $\eclass
\subseteq \GenCons$. The $\ThetaStar$-faithfulness condition is
natural for our setting, as our goal is to estimate the eigen
structure of $ \ThetaStar $, and the set $ \GenCons $ should therefore
represent prior knowledge of this structure and be independent of a
specific factorization of $ \ThetaStar$.

\paragraph{Local descent condition:}  Our next condition provides
a guarantee on the cost improvement that can be obtained by taking a
gradient step when starting from any matrix $\Atwo$ that is
``sufficiently'' far away from the set $\eclass$.

\begin{dfn}[Local descent condition]
\label{DefLtilCurveCon}
For a given radius $\smallrad > 0$, curvature parameter $\curv > 0$
and statistical tolerance $\epsnum \geq 0$, a cost function $\LossTil$
satisfies a \emph{local descent condition} with parameters
$(\curv,\smoo, \epsnum,\smallrad)$ over $\GenCons$ if for each $\Atwo \in
\GenCons \cap \Ballstar{\smallrad}$, there is some $\tdthetastar \in
\arg \min \limits_{A \in \eclass} \fronorm{A - \Atwo}$ such that
\begin{align}
\label{EqnLtilCurveCond}
\trinprod{\nabla \LossTil(\pltheta)}{ \pltheta - \thetastar } &
\geq \curv \fronorm{\Atwo - \tdthetastar}^{2} -
 \frac{\smoo^2 }{\curv} \fronorm{\tdthetastar-\thetastar}^2
- \curv\epsnum^2, \quad \forall  \thetastar \in \eclass.
\end{align}
\end{dfn}

In order to gain intuition for this condition, note that by a
first-order Taylor series expansion, we have $\LossTil(\pltheta) -
\LossTil(\tdthetastar) \approx \trinprod{\nabla \LossTil(\pltheta)}{
  \pltheta - \PISTAR{\pltheta}}$, so that this inner product measures
the potential gains afforded by taking a gradient step.  Now consider
some matrix $\Atwo$ such that $\fronorm{\Atwo - \tdthetastar} >
\sqrt{2} \epsnum$, so that its distance from $\eclass$ is larger than
the statistical precision. The lower bound~\eqref{EqnLtilCurveCond}
with $\thetastar = \tdthetastar $then implies that
\begin{align*}
\LossTil(\pltheta) - \LossTil(\tdthetastar) \approx \trinprod{\nabla
  \LossTil(\pltheta)}{ \pltheta - \PISTAR{\pltheta}} & \geq
\frac{\curv}{2} \fronorm{\pltheta - \PISTAR{\pltheta}}^2,
\end{align*}
which guarantees a quadratic descent condition. Note that the
condition~\eqref{EqnLtilCurveCond} actually allows for additional
freedom in the choice of $\thetastar$ so as to accommodate the
non-uniqueness of the factorization.

One way in which to establish a bound of the
form~\eqref{EqnLtilCurveCond} is by requiring that $\LossTil$ be
locally strongly convex, and that the gradient $\nabla
\LossTil(\tdthetastar)$ approximately vanishes.  In particular,
suppose $\LossTil$ is $2\curv$-strongly convex over the set $\GenCons
\cap \Ballstar{\smallrad}$, in the sense that
\begin{align*}
\trinprod{\nabla \LossTil(\pltheta) - \nabla
  \LossTil(\PISTAR{\pltheta})}{ \pltheta - \PISTAR{\pltheta}} & \geq
2\curv \fronorm{\pltheta - \PISTAR{\pltheta}}^2 \qquad \mbox{for all
  $\pltheta \in \GenCons \cap \Ballstar{\smallrad}$.}
\end{align*}
If we assume that $\frobnorm{\nabla \LossTil(\tdthetastar)} \leq \curv
\epsnum$, then some simple algebra yields that the lower
bound~\eqref{EqnLtilCurveCond} holds.

However, it is essential to note that our theory covers several
examples in which a lower bound~\eqref{EqnLtilCurveCond} of the form
holds, even though $\LossTil$ fails to be locally convex, and/or the
gradient $\nabla \LossTil(\tdthetastar)$ does not approximately
vanish.\footnote{We note that the vanishing gradient condition is needed in all existing work on nonconvex gradient descent~\cite{candes2014wirtinger,sun2014mc,zheng2015convergent}.}  Examples include the problem of sparse PCA, previously
introduced in Section~\ref{SecSparsePCAExample}; in this case, the
function $\LossTil$ is actually globally concave, but nonetheless our
analysis in Section~\ref{SecSparsePCA} shows that a local descent
condition of the form~\eqref{EqnLtilCurveCond} holds.  Similarly, for the
planted clustering model studied in Section~\ref{SecClustering}, the
same form of global concavity holds.  In addition, for the matrix
regression problem previously introduced in
Section~\ref{SecMatrixSensingExample}, we prove in
Section~\ref{SecMatrixSensing} that the condition~\eqref{EqnLtilCurveCond}
holds over a set over which $\LossTil$ is nonconvex.  The generality
of our condition~\eqref{EqnLtilCurveCond} is essential to accommodate
these and other examples. \\

\vtiny

\paragraph{Local Lipschitz condition:}
Our next requirement is a straightforward local Lipschitz property:
\begin{dfn}[Local Lipschitz]  
\label{DefLtilLipCon}
The loss function $\LossTil$ is locally Lipschitz in the sense that
\begin{align}
\label{EqnLtilLipCon}
\fronorm{\nabla \LossTil(\pltheta)} & \leq 
\LipCon\opnorm{\thetastar}.
\end{align}
for all $\pltheta \in \GenCons \cap \Ball_{2}(\smallrad;\thetastar)$.
\end{dfn}

\vtiny

\paragraph{Local smoothness:}
Our last condition is \emph{not} required to establish convergence of
projected gradient descent, but rather to guarantee a faster
geometric rate of convergence.  It is a condition---complementary to
the local descent condition---that upper bounds the behavior of the
gradient map $\nabla \LossTil$.

\begin{dfn}[Local smoothness] 
\label{DefLtilSmoothCon} For some curvature
and smoothness parameters $\curv$ and $\smoo$, statistical tolerance
$\epsnum$ and radius $\smallrad$, we say that the loss function
$\EmpLoss$ satisfies a\emph{ local smoothness condition} with
parameters $(\curv,\smoo,\epsnum,\smallrad)$ over $\GenCons$ if for
each $\pltheta,\pltheta' \in \GenCons \cap
\Ballstar{\smallrad}$ and $\thetastar \in \eclass$,
\begin{align}
\label{EqnLtilSmoothCond}
\abs{\trinprod{\nabla_{\PlTheta} \LossTil( \pltheta ) -
    \nabla_{\PlTheta} \LossTil( \pltheta')}{\pltheta - \thetastar}} 
& \leq
\big(\smoo\fronorm{\pltheta - \pltheta'} + \curv\epsnum \big) \fronorm{\pltheta
  - \thetastar}.
\end{align}
\end{dfn}

The above conditions are stated in terms of the loss function
$\LossTil$ for the factor matrix $\pltheta$. Alternatively, one may
restate these conditions in terms of the loss function $\EmpLoss$ on
the original space, and we make use of this type of reformulation in
parts of our proofs.  For instance, see Section~\ref{SecProofCor} for
details.


\subsection{Sublinear convergence under Lipschitz condition}
\label{SecSublinear}

With our basic conditions in place, we are now ready to state our
first main result.  It guarantees a sublinear rate of convergence
under the $\ThetaStar$-faithfulness, local  descent, and local
Lipschitz conditions.

More precisely, for some descent and Lipschitz parameters $\curv \le
\LipCon $, a statistical tolerance $\epsnum \geq 0$, and a constant
$\inisig \in (0, \frac{1}{2})$, suppose that $\epsnum \leq \frac{1 -
  \inisig}{2} \sigma_{\rdim}(\thetastar)$, the cost functions
$\LossTil$ satisfies the local descent and Lipschitz conditions
(Definitions~\ref{DefLtilCurveCon} and~\ref{DefLtilLipCon}) with
parameters $\curv,\LipCon,\epsnum$ and \mbox{$\smallrad= (1 -
  \inisig)\sigma_{\rdim}(\thetastar)$,} and the constraint set $
\GenCons $ is $ \ThetaStar $-faithful and convex.  Let $ \cond =
\cond(\thetastar) \defn
\frac{\sigma_1(\thetastar)}{\sigma_\rdim(\thetastar)} $ be the
condition number of $ \thetastar $. We then have the following
guarantee:
\begin{thm}
\label{ThmGeneralLip}
Under the previously stated conditions, given any initial point
$\thetait{0}$ belonging to the set $\GenCons \cap \Ball_{2}((1 -
  \inisig)\sigma_{\rdim}(\thetastar);\thetastar)$, the projected gradient
iterates $\{\thetait{t}\}_{t= 1}^{\infty}$ with step size $\stepit{t}
= \frac{1}{\curv(t + 20\cond^2\LipCon^{2}/\curv^{2})}$ satisfy the bound
\begin{align}
\label{EqnLipGuarantee}
\distsqr{\thetait{t}} &
\leq \frac{20 \LipCon^{2} \opnorm{\thetastar}^{2}}{t\curv^{2}} + 4\epsnum^{2} 
\qquad \mbox{for all iterations $t = 1, 2, \ldots$.}
\end{align}
\end{thm}
\noindent See Section~\ref{SecProofThmGeneralLip} for the proof of
this claim.

As a minor remark, we note that the assumption $\epsnum\le \frac{1 -
  \inisig}{2} \sigma_{\rdim}(\thetastar)$ entails no loss of
generality---if it fails to hold, then the initial solution
$\thetait{0}$ already satisfies an error bound better than what is
guaranteed for subsequent iterates.

Conceptually, Theorem~\ref{ThmGeneralLip} provides a minimal set of
conditions for the convergence of projected gradient descent using the
nonconvex factorization $\PlTheta = \OUTER{\pltheta}$.  The first
term on the right hand side of equation~\eqref{EqnLipGuarantee}
corresponds to the \emph{optimization error}, whereas the second
$\epsnum^2$ term is the \emph{statistical error.} The
bound~\eqref{EqnLipGuarantee} shows that the distance between
$\thetait{t}$ and $\thetastar$ drops at the rate $\order(\frac{1}{t})$
up to the statistical limit $\epsnum^2$ that is determined by the
sample size and the signal-to-noise ratio (SNR) of the problem.  We
see concrete instances of this statistical error in the examples to
follow.


\subsection{Linear convergence under smoothness condition}
\label{SecLinear}

Although Theorem~\ref{ThmGeneralLip} does guarantee convergence, the
resulting rate is sublinear ($\order(1/t)$), and hence rather slow.
In this section, we show that if in addition to the local Lipschitz
and descent conditions, the function $\LossTil$ satisfies the local
smoothness conditions in Definition~\ref{DefLtilSmoothCon}, then much
faster convergence can be guaranteed.

More precisely, suppose that for some numbers $\curv$, $\smoo$,
$\LipCon$, $\epsnum$ and $\inisig$ with $0 < \curv \le \smoo =
\LipCon$, $ 0< \inisig <1 $ and $\epsnum \le \frac{1 - \inisig}{4} \sigma_{\rdim}(\thetastar)$,
the loss function $\LossTil$ satisfies the local descent, Lipschitz and smoothness conditions in
Definitions~\ref{DefLtilCurveCon}--\ref{DefLtilSmoothCon} over $\GenCons$ with
parameters $\curv$, $\smoo$, $\LipCon$, $\epsnum$ and $\smallrad= (1 -
\inisig^{2})\sigma_{\rdim}(\thetastar)$, and that the set
$\GenCons$ is $\ThetaStar$-faithful and convex. 

\begin{thm}
\label{ThmGeneralSmooth}
Under the previously stated conditions, there is a
constant $0<c_{\inisig} < 1$ depending only on $\inisig$ such that
given an initial matrix $\thetait{0}$ in the set $\GenCons \cap
\Ball_{2}((1 - \inisig) \sigma_{\rdim}(\thetastar);\thetastar)$, the
projected gradient iterates $\{\thetait{t}\}_{t= 1}^{\infty}$ with
step size $\stepit{t} = c_{\inisig} \frac{\curv}{\cond^6\smoo^{2}}$ satisfy
the bound
\begin{align}
\label{EqnFinalLinearGuarantee}
\distsqr{\thetait{t}} & \leq \Big( 1-c_{\inisig} \frac{\curv^{2}}{\cond^6\smoo^{2}} \Big)^{t}\distsqr{\thetait{0}} + 16 \epsnum^{2}
\qquad \mbox{for all iterations $t = 1, 2, \ldots$.}
\end{align}
\end{thm}
\noindent See Section~\ref{SecProofThmGeneralSmooth} for the proof of
this claim.

The right hand side of the bound~\eqref{EqnFinalLinearGuarantee} again
consists of an optimization error term and a statistical error
term. The theorem guarantees that the optimization error converges
linearly at the geometric rate $ \order((1-c)^t) $ up to a statistical
limit.  Note that the theorem requires the initial solution
$\thetait{0} $ to lie within a ball around $ \thetastar $ with radius
$(1 - \inisig) \sigma_{\rdim}(\thetastar) $, which is slightly smaller
than the radius $ \smallrad = (1-\inisig^2) \sigma_{\rdim}(\thetastar)
$ for which the local descent, Lipschitz and smoothness conditions
hold. Moreover, the step size and the convergence rate depend on the
condition number of $ \thetastar $ as well as the quality of the
initialization through $ \inisig $. We did not make an attempt to
optimize this dependence, but improvement in this direction, including
adaptive choices of the step size, is certainly an interesting problem
for future work.


\section{Concrete results for specific models}
\label{SecConcrete}

In this section, we turn to the consequences of our general theory for
specific models that arise in applications.  Throughout this section,
we focus on geometric convergence guaranteed by
Theorem~\ref{ThmGeneralSmooth} using a constant step size.  The main
technical challenges are to verify the local descent, local Lipschitz
and local smoothness assumptions that are needed to apply this result.
Since Theorem~\ref{ThmGeneralLip} depends on weaker assumptions---it
does not need the local smoothness property---it should be understood
also that our analysis can be used to derive corollaries based on
Theorem~\ref{ThmGeneralLip} as well.

\paragraph{Note:}  In all of the analysis to follow,
we adopt the shorthand $\sigma_j = \sigma_j(\thetastar)$ for the
singular values of $\thetastar$, and $\cond =
\frac{\sigma_1}{\sigma_r} $ for its condition number.


\subsection{Noisy matrix completion}
\label{SecMC}

We begin by deriving a corollary for the problem of noisy matrix
completion.  Since we did not discuss this model in
Section~\ref{SecIllustrative}, let us provide some background here.
There are a wide variety of matrix completion problems
(e.g.,~\cite{Lau01}), and the variant of interest here arises when the
unknown matrix has low rank.  More precisely, for an unknown PSD and
low-rank matrix $\ThetaStar \in \SymMat{\usedim}$, suppose that we are
given noisy observations of a subset of its entries.  In the so-called
Bernoulli model, the random subset of observed entries is chosen
uniformly at random---that is, each entry is observed with some
probability~$\pobs$, independently of all other entries.  We can
represent these observations by a random symmetric matrix $Y \in
\SymMat{\usedim}$ with entries of the form
\begin{align}
 \label{EqnMC}
Y_{ij} & = \begin{cases} \PlTheta_{ij}^{*} + \noisevar_{ij}, &
  \mbox{with probability \ensuremath{\pobs}, and}\\ * &
  \text{otherwise.}
\end{cases}, \quad \text{for each }i\ge j.
\end{align}
Here the variables $\{\noisevar_{ij},i\ge j\}$ represent a form of
measurement noise.

A standard method for matrix completion is based on solving the
semidefinite program
\begin{align}
\label{EqnMCLoss}
\min_{\PlTheta \in \SymMat{\usedim}} \Big \{
\underbrace{\frac{1}{2\pobs}\sum_{(i,j) \in \Omega}
  \left(\PlTheta_{ij} - \Yout_{ij}\right)^{2}}_{\EmpLoss(\PlTheta)}
\Big \} \qquad \mbox{such that \ensuremath{\PlTheta \succeq 0 \quad
    \mbox{and} \quad \nucnorm{\PlTheta} \leq R},}
\end{align}
where $R > 0$ is a radius to be chosen.  As noted above, the PSD
constraint and nuclear norm bound are equivalent to the trace
constraint $\trace(\PlTheta) \leq R$. In either case, this is a
special case of our general estimator~\eqref{EqnEmpiricalSDP}.

The SDP-based estimator~\eqref{EqnMCLoss} is known to have good
performance when the underlying matrix $\ThetaStar$ satisfies certain
matrix incoherence conditions.  These conditions involve its leverage
scores, defined in the following way. Here we consider a simplified
setting where the eigenvalues of $ \ThetaStar $ are equal.  By
performing an eigendecomposition, we can write $\ThetaStar = U D U^T$
where $D \in \real^{\rdim \times \rdim}$ is a diagonal matrix of
eigenvalues (a constant multiple of the identity when they are
constant), and take $\thetastar = U D^{1/2}$.  With this notation, the
\emph{incoherence} parameter of $\ThetaStar= \OUTER{\thetastar}$ is
given by
\begin{equation}
 \label{EqnIncoherence}
\inco \defn \frac{\usedim \max \limits_{i= 1,\ldots,\usedim}
  \twonorm{\thetastar_{i\cdot}}^{2}}{ \rdim \opnorm{\thetastar}^{2}} =
\frac{\usedim\twoinfnorm{\thetastar}^{2}}{\rdim\opnorm{\thetastar}^{2}}.
\end{equation}

Since we already enforce low-rankness in the factorized formulation,
we can drop nuclear norm constraint. The generalized projected
gradient descent~\eqref{EqnGenProjGradOpt} is specified by letting $ \LossTil $
 and $\GenCons$ set
\begin{align*}
\LossTil(\pltheta) \defn \frac{1}{2\pobs}\sum_{(i,j) \in \Omega}
  \big( (\OUTER{\pltheta})_{ij} - \Yout_{ij} \big)^{2}
\quad\mbox{and}\quad
\GenCons  \defn \Big\{\pltheta \in \real^{\usedim\times \rdim} \,
\mid\, \twoinfnorm{\pltheta} \le
\sqrt{\frac{2\inco}{\usedim}}\fronorm{\thetait{0}} \Big\}.
\end{align*}
Note that $ \GenCons $ is convex, and depends on the initial solution
$ \Thetait{0} $.  The gradient of $\EmpLoss$ is $\nabla_{\PlTheta}
\EmpLoss(\PlTheta)= \frac{1}{\pobs}\ProjObs(\PlTheta - \Yout)$, and
the projection $\ProjGenCons$ is given by the row-wise ``clipping''
operation
\begin{align*}
[\ProjGenCons(\theta)]_{i\cdot} = 
\begin{cases} 
\pltheta_{i\cdot}\;, & \twonorm{\pltheta_{i\cdot}} \le
\sqrt{\frac{2\inco\rdim}{\usedim}}
\opnorm{\thetait{0}},\\ \pltheta_{i\cdot}
\sqrt{\frac{2\inco\rdim}{\usedim}}
\frac{\opnorm{\thetait{0}}}{\twonorm{\pltheta_{i\cdot}}} \;, &
\twonorm{\pltheta_{i\cdot}} > \sqrt{\frac{2\inco\rdim}{\usedim}}
\opnorm{\thetait{0}},
\end{cases}
\quad \text{for \ensuremath{i = 1,2,\ldots,\usedim}.}
\end{align*}
This projection ensures that the iterates of gradient
descent~\eqref{EqnGenProjGradOpt} remains incoherent.

\begin{rem}
\label{RemThetaStarRowNorm}
Note that $\twonorm{\thetastar_{i\cdot}}^{2} =
\twonorm{\Pi_{\text{col}(\thetastar)}(e_{i})}^{2}$ ($e_{i}$ is the
$i$-th standard basis vector and $\text{col}(\thetastar)$ is the
column space of~$\thetastar$), so the values of
$\twoinfnorm{\thetastar}$ and $\inco$ depend only on
$\text{col}(\thetastar)$ and are the same for any $\thetastar$ in
$\eclass$.
\end{rem}

With this notation in place, we are now ready to apply
Theorem~\ref{ThmGeneralSmooth} to the noisy matrix completion problem.
As we show below, if the initial matrix $\thetait{0}$ satisfies the
bound $\dist{\thetait{0}} \le \frac{1}{5} \opnorm{\thetastar}$, then
the set $ \GenCons $ is $\ThetaStar $-faithful. Moreover, if the
expected sample size satisfies $\numobs= \pobs\usedim^{2} \succsim
\max\{\inco\rdim\usedim\log\usedim, \inco^{2}\rdim^{2}\usedim\}$, then
with probability at least $1-4 \usedim^{-3}$ the loss
function~$\LossTil$ satisfies the local descent, Lipschitz and
smoothness conditions with parameters
\begin{align*}
\smallrad = \frac{3}{5} \opnorm{\thetastar}, \quad \curv =
\frac{2}{25} \opnorm{\thetastar}^{2}, \quad \LipCon = \smoo =
c_{2}\inco\rdim\opnorm{\thetastar}^{2} \quad \text{and} \quad \epsnum=
100\frac{\sqrt{\rdim}
  \opnorm{\ProjObs(\noisevar)}}{\pobs\opnorm{\thetastar}}.
\end{align*} 

Using this fact, we have the following consequence of
Theorem~\ref{ThmGeneralSmooth}, which holds when the sample size size
$ \numobs $ satisfies the bound above and is large enough to ensure
that $\epsnum \le \frac{1}{10} \opnorm{\thetastar}$.
\begin{cor}
\label{CorMC}
Under the previously stated conditions, if we are given an initial
matrix $\thetait{0}$ satisfying the bound $\dist{\thetait{0}} \le
\frac{1}{5} \opnorm{\thetastar}$, then with probability at least $1-4
\usedim^{-3}$, the gradient iterates $\{\thetait{t}\}_{t= 1}^{\infty}$
with step size $\stepit{t} = c_{3} \frac{1}{\inco^{2}\rdim^{2}
  \opnorm{\thetastar}^{2}}$ satisfy the bound
\begin{align}
\label{EqnMCBound}
\distsqr{\thetait{t}} & \le \Big( 1-c_{4} \frac{1}{\inco^{2}\rdim^{2}}
\Big)^{t}\distsqr{\thetait{0}} + c_{5}
\frac{\rdim\opnorm{\ProjObs(\noisevar)}^{2}}{\pobs^{2}
  \opnorm{\thetastar}^{2}}.
\end{align}
\end{cor}
\noindent See Section~\ref{SecProofCorMC} for the proof of this
claim.\\

Even though Corollary~\ref{CorMC} is a consequence of our general
theory, it leads to results for exact/approximate recovery in the
noiseless/noisy setting that are as good as or better than known
results.  In the noiseless setting ($\noisevar = 0$), our sample size
requirement and contraction factor are sharper than those in the
paper~\cite{sun2014mc} by a polynomial factor in the rank
$\rdim$. Turning to the noisy setting, suppose the noise matrix
$\noisevar$ has independent sub-Gaussian entries with parameter
$\noisestd^{2}$. A technical result to be proved later (see
Lemma~\ref{LemCensorGaussMat}) guarantees that given a sample size
$\numobs= \pobs\usedim^{2} \succsim c_{1}\usedim\log^{2}\usedim$, we
have the operator norm bound $\opnorm{\ProjObs(\noisevar)} \precsim
\noisestd\sqrt{\pobs\usedim}$ with probability at least $1 -
\usedim^{-12}$.  Together with the bound~\eqref{EqnMCBound}, we
conclude that
\begin{align}
\label{EqnMC_consequence}
\frac{1}{d^{2}} \fronorm{\OUTER{\thetait{\infty}} -
  \OUTER{\thetastar}}^{2} \le \frac{3}{\usedim^{2}}
\opnorm{\thetastar}^{2} \fronorm{\thetait{\infty} -
  \thetastarit{\infty}}^{2} \precsim
\frac{\noisestd^{2}\rdim}{\pobs\usedim} =
\frac{\noisestd^{2}\rdim\usedim}{\numobs}.
\end{align}
The scaling $\frac{\sigma^{2}\rdim\usedim}{\numobs}$ is better than
the results in past work
(e.g.,~\cite{negahban2010noisymc,koltchinskii2011matrixCompletion,keshavan2010noisy})
on noisy matrix completion by a $\log\usedim$ factor; in fact, it
matches the minimax lower bounds established in the
papers~\cite{negahban2010noisymc,koltchinskii2011matrixCompletion}.
Thus, Corollary~\ref{CorMC} in fact establishes that the projected
gradient descent method yields minimax-optimal estimates.


\begin{figure}
\centering
\begin{tabular}{ccc}		
\includegraphics[scale=0.5, clip, trim = 0 0 0 0]{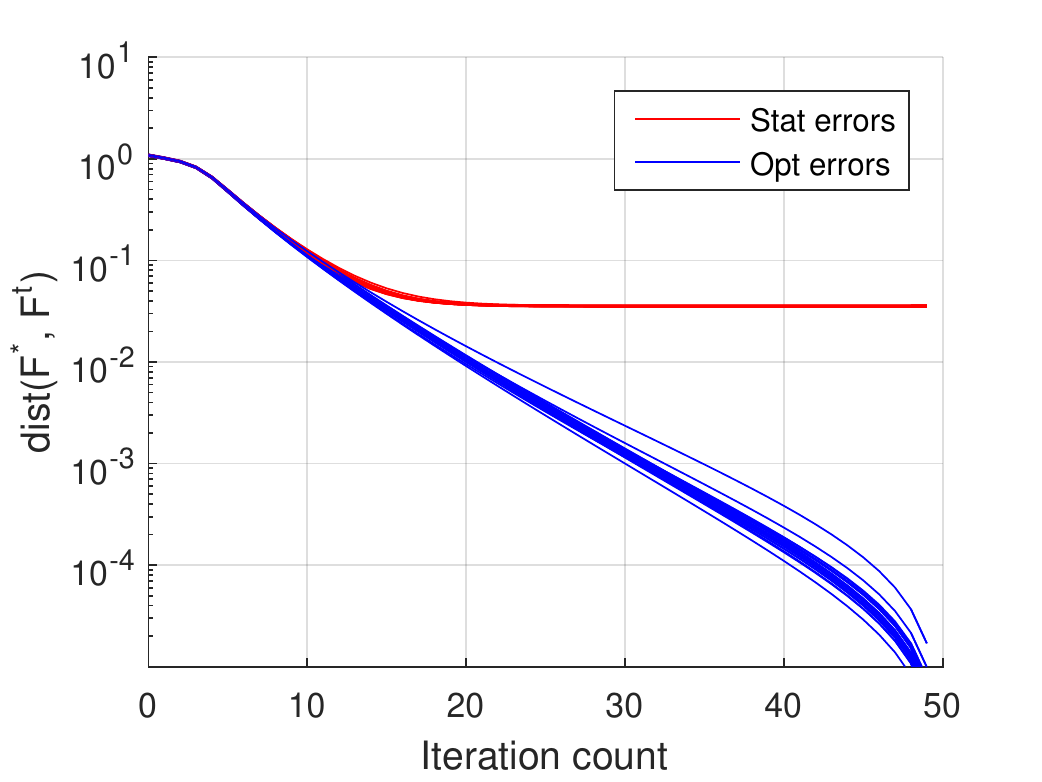} &
\includegraphics[scale=0.5, clip, trim = 0 0 0 0]{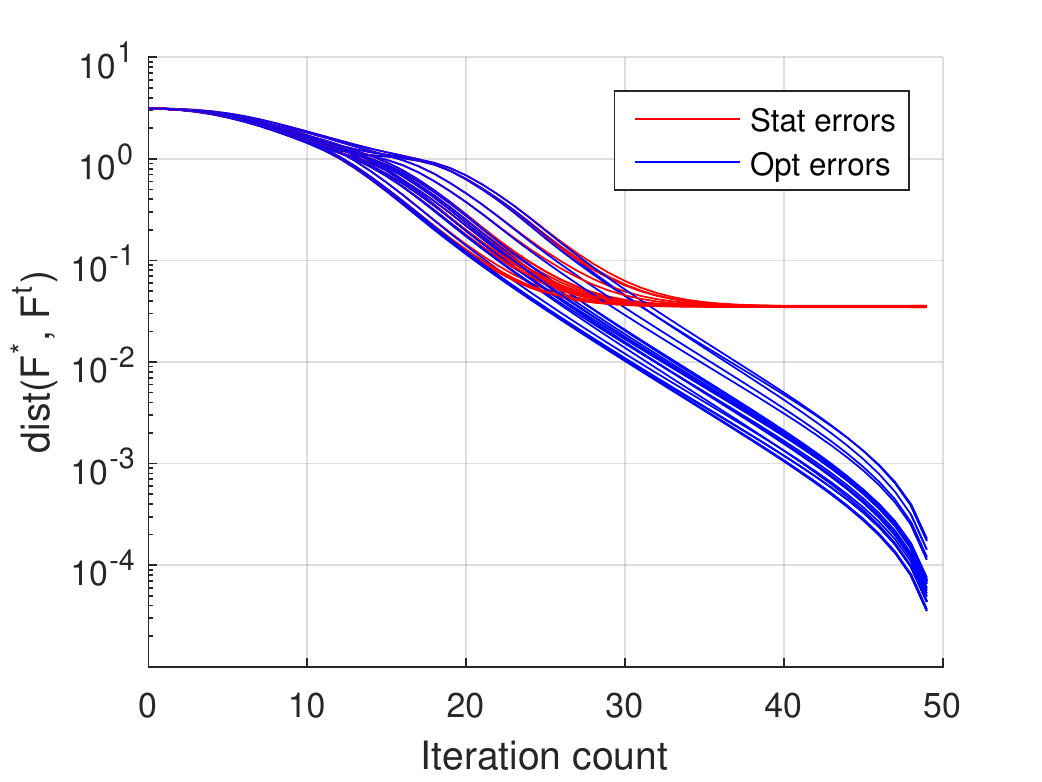} &
\includegraphics[scale=0.5, clip, trim = 0 0 0 0]{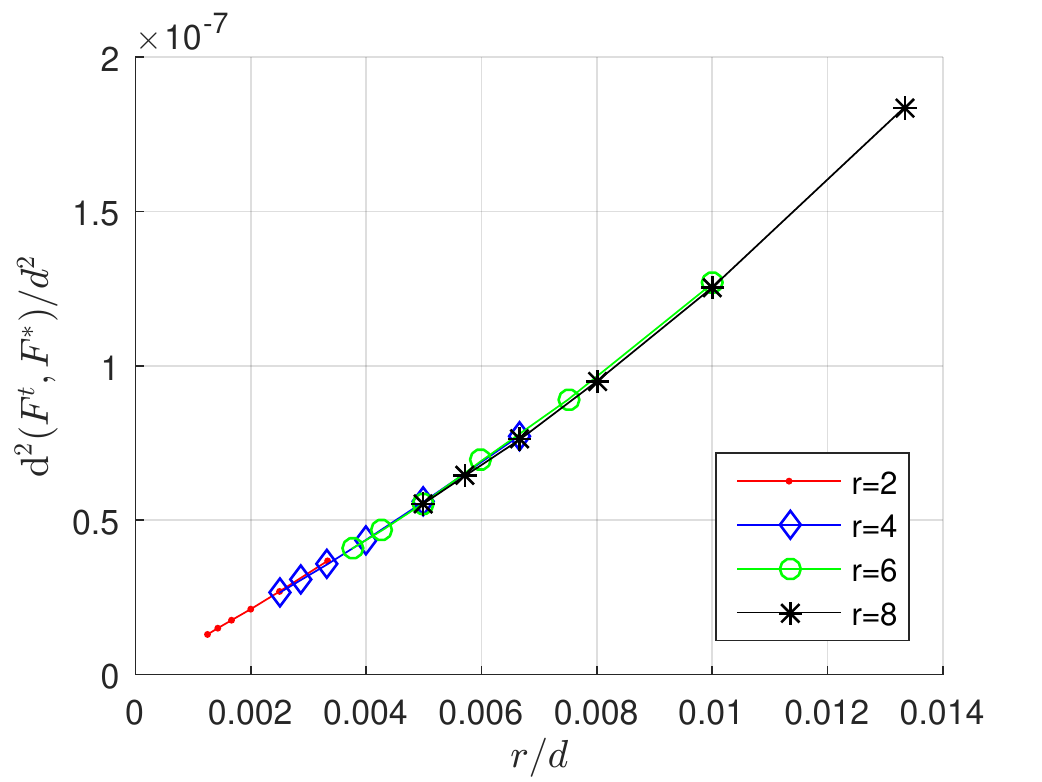} \\
(a) & (b) & (c) 
\end{tabular} 
\caption{Simulation results for matrix completion.  (a) Plots of
  optimization error $ \textup{d} (\thetait{t}, \thetait{T} ) $ and
  statistical error $ \dist{\thetait{t}} $ versus the iteration number
  $ t$ using SVD initialization. Panel (b): same plots using a random
  initialization. The simulation is performed using $ \usedim = 1000$,
  $\rdim=10$, $ \pobs=0.1 $ and $ \noisestd = 0.01 \cdot
  \frac{\rdim}{\usedim}$.  Panel (c): plots of per-entry estimation
  error $ \frac{1}{\usedim^2} \dist{\widehat{\pltheta}} $ versus $
  \frac{\rdim}{\usedim} $, for different values of $ (\usedim, \rdim)
  $ using SVD-based initialization. Each point represents the average
  over $ 20 $ random instances. The simulation is performed using $
  \pobs=0.1 $ and $ \noisestd = 0.001 $.}
\label{FigMC}
\end{figure}

\paragraph{Initialization:}

Suppose the rank-$\rdim$ SVD of the matrix $\ProjObs(M)$ is given by $
USV^\top $. We can take $ \thetait{0} = \ProjGenCons(US^{\frac{1}{2}})
$. Under the previously stated condition on the sample size $ \numobs
$, the matrix $ \thetait{0} $ satisfies the requirement in
Corollary~\ref{CorMC} as shown in, e.g.,~\cite{keshavan2009matrixafew}
(combined with the above bound on $ \opnorm{\ProjObs(\noisevar)} $).

\paragraph{Computation: }

Computing the gradient $\nabla_{\pltheta} \LossTil(\pltheta) =
\frac{2}{\pobs}\ProjObs(\OUTER{\pltheta} - \Yout)\pltheta$ takes time
$\order(\rdim^{2}\abs{\Obs})$. The projection $\ProjGenCons(\pltheta)$
can be computed in time $\order(\rdim\usedim)$.

\paragraph{Simulations:}
In order to illustrate the predictions of Corollary~\ref{CorMC}, we
performed a number of simulations.  Since the distance measures
$\dist{\pltheta}$ and $ \fronorm{\OUTER{\pltheta} -
  \OUTER{\thetastar}} $ are difficult to compute, so we instead use
the subspace distance
\begin{align}
\label{EqnSineApprox}
\dist{\pltheta} \approx \fronorm{\sin \angle (\pltheta,
  \thetastar)}^2,
\end{align}
as an approximation.\footnote{This approximation is valid up to a
  constant of $ 2 $ if both $ \pltheta $ and $ \thetastar $ are
  orthonormal~(cf.\ Proposition~2.2 in \cite{vu2013minimax}).}  Here $
\sin \angle (\pltheta, \thetastar) $ is the vector of principal angles
between the column spaces of matrices $ \pltheta $ and $ \thetastar $.
For each example and given values of model parameters $ \usedim,
\rdim, \numobs, \noisestd $ etc., we generate a random instance by
sampling the true matrix $ \thetastar $ and the problem data randomly
from the relevant model, and then run our projected gradient descent
algorithm with $ T=50 $ iterations.

In the matrix completion case, we sampled the true matrix $\thetastar$
uniformly at random from all $ \usedim\times\rdim $ orthonormal matrix
uniformly at random, generated a noise matrix $ \noisevar $ with
i.i.d. $N(0, \noisestd^2)$ entries, and chose the observed entries
randomly according to the Bernoulli model with probability~$ \pobs
$. We considered two approaches for obtaining the initial matrix $
\thetait{0} $: (a) the SVD-based procedure described in
Section~\ref{SecMC}, and (b) random initialization, where $
\thetait{0} $ is a random $ \usedim \times \rdim $ orthonormal matrix
projected onto the associated constraint set $ \GenCons $. The step
size for projected gradient descent is fixed at $ \stepit{t} \equiv
\frac{0.5}{\pobs}$.  Panels (a) and (b) in Figure~\ref{FigMC} show the
resulting convergence behavior of the algorithm, which confirm the
geometric convergence (and threshold effect for the statistical error)
that is predicted by our theory.

For these random ensembles, our theory predicts thatwith high
probability the per-entry error of the output $ \widehat{\pltheta} $
satisfies a bound of the form
	\begin{align*}
	\frac{1}{\usedim^2}\distsqr{\widehat{\pltheta}} \precsim
        \frac{\noisestd^2 \rdim}{\pobs \usedim};
	\end{align*}
	cf.\ equation~\eqref{EqnMC_consequence}. Therefore, with $
        \pobs $ and $ \noisestd $ fixed, the ratio
        $\frac{1}{\usedim^2}\distsqr{\widehat{\pltheta}} $ should be
        proportional to $ \frac{\rdim}{\usedim} $.

\subsection{Matrix regression}
\label{SecMatrixSensing}

Recall the matrix regression model previously introduced in
Section~\ref{SecMatrixSensingExample}.	In order to simplify notation,
it is convenient to introduce define a linear mapping \mbox{$\Xmap:
  \real^{\usedim\times \usedim} \mapsto \real^{\numobs}$} via
$[\Xmap(\PlTheta)]_{i} \defn \trinprod{\Xit{i}}{\PlTheta}$ for $i = 1,
2, \ldots, \numobs$.  Note that the adjoint operator
$\Xmap^{*}:\real^{\usedim} \mapsto \real^{\usedim\times \usedim}$ is
given by $\Xmap^{*}(u)= \sum_{i= 1}^{\numobs} \Xit{i} u_{i}$.  With
this notation, we have the compact representation
\begin{align*}
\nabla \LossTil(\pltheta) & = \frac{2}{\numobs} \Big( \Xmap^{*}(\Xmap(
\OUTER{\pltheta}) - \yout) \Big) \pltheta.
\end{align*}
Since $ \trinprod{\Xit{i}}{ \OUTER{\pltheta} } = \trinprod{(\Xit{i} + {\Xit{i}}^\top)/2}{
 \OUTER{\pltheta} }$, we may assume without loss of generality that the matrices $\{ \Xit{i} \}$ are symmetric.
 
In this case, projected gradient descent can be performed with
$\GenCons = \real^{\usedim \times \rdim}$, so that the
$\ThetaStar$-faithfulness condition holds trivially.  It remains to
verify that the cost function $\LossTil$ from
equation~\eqref{EqnMatrixSensingLossTil} satisfies the local descent,
local Lipschitz and local smoothness properties, and these properties
depend on the structure of the operator $\Xmap$.  For instance, one
way in which to certify the conditions of
Theorem~\ref{ThmGeneralSmooth} is via a version of \emph{restricted
  isometry property} (RIP) applied to the operator $\Xmap$.

\begin{dfn}[Restricted isometry property]
\label{DefRIP}
The operator $\Xmap: \real^{\usedim \times \usedim} \rightarrow
\real^\numobs$ is said to satisfied the restricted isometry property with
parameter $\ripparam{k}$ if
\begin{align*}
(1 - \ripparam{k}) \fronorm{\PlTheta}^{2} \le \frac{1}{\numobs}
  \twonorm{\Xmap(\PlTheta)}^{2} \le (1 + \ripparam{k})
  \fronorm{\PlTheta}^{2}, \quad \mbox{for all $\usedim$-dimensional
    matrices with $\rank(\PlTheta) \leq k$}.
\end{align*}
\end{dfn}
It is well known~\cite{recht2010guaranteed,negahban2009estimation}
that RIP holds for various random ensembles.  For instance, suppose
that the entries of $\Xit{i}_{j\ell}$ are i.i.d. zero-mean unit variance
random variables, satisfying a sub-Gaussian tail bound.  Examples of
such ensembles include the standard Gaussian case ($\Xit{i}_{j\ell} \sim
N(0,1)$) as well as Rademacher variables ($\Xit{i}_{j\ell} \in \{-1, 1\}$
equiprobably).	For such ensembles, it is known that with high
probability, a RIP condition of order $\rdim$ holds with a sample size
$\numobs \succsim \rdim \usedim$. \\

The RIP condition provides a straightforward way of verifying the
conditions of Theorem~\ref{ThmGeneralSmooth}.  More precisely, as we
show in the proof of Corollary~\ref{CorMatrixSensing}, if the operator
$\Xmap$ satisfies RIP with parameter $ \ripparam{4\rdim} \in [0,
\frac{1}{12} )$, then the loss function $\LossTil$ satisfies the local
descent, descent and smoothness conditions with parameters
\begin{align*}
\smallrad = (1-12\ripparam{4\rdim}) \sigma_\rdim, 
\quad \curv = 6\ripparam{4\rdim} \sigma_\rdim^{2}, 
\quad \LipCon = \smoo= 64 \cond^2 \sigma_\rdim^{2} 
\quad \text{and} \quad \epsnum =
\frac{2\sqrt{\rdim} \cond
  \opnorm{\numobs^{-1}\Xmap^{*}(\smallnoisevar)}}{\ripparam{4\rdim}
  \sigma_\rdim}.
\end{align*}

Using this fact, we have the following corollary of
Theorem~\ref{ThmGeneralSmooth}.  We state it assuming that the operator $ \Xmap $ satisfies RIP with parameter $ \ripparam{4\rdim} \in [0,
\frac{1}{12} )$, and  the sample size $\numobs$ is large enough to ensure that $\epsnum \leq \frac{1-\sqrt{12\ripparam{4\rdim}}}{4} \sigma_\rdim$.
\begin{cor}
\label{CorMatrixSensing}
Under the previously stated conditions, there is a universal function
$\MYPSI:[0, 1/12] \rightarrow (0,1)$ such that given any initial
matrix $\thetait{0}$ satisfying the bound \mbox{$\dist{\thetait{0}} \le
  (1-\sqrt{12\ripparam{4\rdim}}) \sigma_\rdim$,} the projected
gradient iterates $\{\thetait{t}\}_{t= 1}^{\infty}$ with step size
$\stepit{t} =
\frac{\MYPSI(\ripparam{4\rdim})} {\ripparam{4\rdim} \cond^{10} \sigma_\rdim^{2}}$
satisfy the bound
\begin{align*}
\distsqr{\thetait{t}} & \le \Big( 1-\frac{\MYPSI(\ripparam{4\rdim})}{\cond^{10}}
\Big)^{t} \distsqr{\thetait{0}} + c_{0} 
\frac{\rdim \cond^2  \opnorm{\Xmap^{*}(\smallnoisevar)}^{2}}{\numobs^{2} \ripparam{4\rdim}^2 \sigma_\rdim^{2}},
\end{align*}
\end{cor}
\noindent See Section~\ref{SecProofCorMatSensing} for the proof of
this claim.\\


Note that the radius of the region of convergence
$\Ballstar{(1-\sqrt{12\ripparam{4\rdim}}) \sigma_\rdim}$ can be arbitrarily close to $ \sigma_\rdim $ with $ \ripparam{4\rdim} $ sufficiently small. Moreover, the function $\LossTil$ need not be convex in this region.  As a simple example, consider
the scalar case with $\usedim= \rdim = 1$, with noiseless observations
($\smallnoisevar= 0$) of the target parameter $\thetastar= 1$.	A
simple calculation then yields that $\LossTil(\pltheta) = c \,
(\pltheta^{2}-1)^{2}$ for some  constant $c > 0$, which is
nonconvex outside of the ball $\Ballstar{\frac{1}{\sqrt{3}}}$.

Specified to the noiseless setting with $\smallnoisevar= 0$, Corollary~\ref{CorMatrixSensing} is similar to the results for nonconvex gradient descent in~\cite{tu2015procrustes,zheng2015convergent}. In the more general noisy setting, our statistical error rate $\epsnum$ is consistent with the results in~\cite{negahban2010noisymc}.  For a more concrete example, suppose $ \cond = \order(1) $, and  each $\Xit{i}$ and $\smallnoisevar$ have i.i.d.\ Gaussian entries with $\Xit{i}_{j\ell} \sim \mathcal{N}(0,1)$ and $\smallnoisevar_{i} \sim \mathcal{N}(0,\noisestd^{2})$. It can be
shown that as long as $\numobs \succsim \rdim\usedim\log\usedim$, RIP holds with $\ripparam{4\rdim} < \frac{1}{192}$ and $\opnorm{\Xmap^{*}(\smallnoisevar)} \precsim
\noisestd\sqrt{\numobs\usedim}$.  The bound in Corollary~\ref{CorMatrixSensing} therefore implies a constant contraction factor and that 
\begin{align*}
 \fronorm{\OUTER{\thetait{\infty}} - \OUTER{\thetastar}}^{2} &
 \le 3\distsqr{\thetait{\infty}} \opnorm{\thetastar}^{2} \precsim
 \noisestd^{2} \frac{\rdim\usedim}{\numobs}.
\end{align*}

\paragraph{Initialization:}

Suppose the rank-$ \rdim $ SVD of the matrix $\frac{1}{\numobs}\Xmap^{*}(y)$ is given by $ USV^\top $. We can take $\thetait{0} = US^{\frac{1}{2}}$. Under the
above Gaussian example, it can be shown the condition on
the initial solution is satisfied if $\numobs \succsim
\usedim\rdim^{2} \cond^4 \log\usedim$ and $ \noisestd $ is small enough~\cite{jain2010svp,jain2013low}.

The sample size required for this initialization scales quadratically
in the rank $\rdim$, as compared to the linear scaling that is the
best possible~\cite{recht2010guaranteed,negahban2009estimation}.  This
looseness is a consequence of requiring the initialization error to
satisfy a Frobenius norm bound instead of an operator norm one. It can
be avoided by using a more sophisticated initialization
procedures---for instance, one based on a few iterations of the
singular value projection (SVP) algorithm~\cite{jain2010svp}.  In the
current setting, since our primary focus is on understanding
low-complexity algorithms via gradient descent, we do not pursue this
direction further.

\paragraph{Computation:}

Let $T_{\text{mul}}$ be the maximum time to multiply $\Xit{i}$ with a
vector in $\real^{\usedim}$. Finding the initial solution as above
requires computing the rank-$\rdim$ SVD of the $\usedim \times
\usedim$ matrix $\frac{1}{\numobs}\Xmap^{*}(\yout)$, which can be done
in time $\order(\numobs\rdim T_{\text{mul}} + \usedim\rdim^{2})$;
cf.~\cite{halko2011finding}.  The gradient
$\frac{1}{\numobs}\Xmap^{*}(\Xmap(\PlTheta) - \yout)$ and can be
computed in time $\order(\numobs\rdim T_{\text{mul}} + \usedim\rdim)$.
Therefore, the overall time complexity is $\order(\numobs\rdim
T_{\text{mul}} + \usedim\rdim^{2})$ times the number of iterations.


\subsection{Rank-$\protect\rdim$ PCA with row sparsity}
\label{SecSparsePCA}

Recall the problem of sparse PCA previously introduced in
Section~\ref{SecSparsePCAExample}.  In this section, we analyze the
projected gradient updates applied to this problem, in particular with
the loss function $\LossTil$ from equation~\eqref{EqnLossTilPCA}, and
the constraint set
\begin{align*}
 \GenCons \defn \big\{\pltheta \in \real^{\usedim\times \rdim} \, \mid
 \, \opnorm{\pltheta} \leq 1, \, \twoonenorm{\pltheta} \leq
 \twoonenorm{\thetastar} \big\}.
\end{align*}
To be clear, this choice of constraint set is somewhat unrealistic,
since it assumes knowledge of the norm $\twoonenorm{\thetastar}$.
This condition could be removed by analyzing instead a penalized form
of the estimator, but as our main goal is to illustrate the general
theory, we remain with the constrained version here.  

We now apply Theorem~\ref{ThmGeneralSmooth} to this problem. As we show in the proof of Corollary~\ref{CorSparsePCA}, the set $ \GenCons $ is $ \ThetaStar $-faithful. Moreover, for each $0 < \inisig < 1$, suppose that the SNR satisfies $\snr >\frac{2}{\inisig^{2}}$ in the row-sparse spiked covariance model,  then with probability at least $1-2\usedim^{-3}$, the loss
function $\LossTil$ satisfies the local descent, smoothness
conditions and the relaxed Lipschitz condition~(\ref{EqnWeakLipCon})
with parameters
\begin{align}
\smallrad = 1 - \inisig^{2}, \quad
\curv = \frac{\snr\inisig^{2}}{4}, \quad
\LipCon = \smoo= 4(\gamma + 1)\sqrt{\rdim} \quad \mbox{and} \quad \epsnum = c_{1} \frac{\snr + 1}{\snr\inisig^{2}}\sqrt{\rdim} \max \Big\{\sqrt{\frac{\kdim\log\usedim}{\numobs}}, \frac{\kdim\log\usedim}{\numobs} \Big\}.
\label{EqnSparsePCACurve}
\end{align} 
Using these fact, we have the following corollary of Theorem~\ref{ThmGeneralSmooth}. We state it assuming that the SNR obeys the bound  $\snr >\frac{2}{\inisig^{2}}$ and the sample size $ \numobs $ is large enough to ensure  $\epsnum\le \frac{1 - \inisig}{20}$
\begin{cor}
\label{CorSparsePCA}
Under the previously stated conditions, there is a function $ \MYPSI:
(0, 1) \rightarrow (0,1) $ such that given any initial matrix
$\thetait{0} \in \GenCons \cap \Ballstar{1 - \inisig}$, with
probability at least $1-2\usedim^{-3}$, the projected gradient
iterates $\{\thetait{t}\}_{t= 1}^{\infty}$ with step size $\stepit{t}
= \MYPSI(\inisig) \frac{\snr}{(\snr + 1)^{2}\rdim}$ satisfy the bound
\begin{align*}
\distsqr{\thetait{t}}^{2} & \leq \Big( 1-\MYPSI(\inisig)
\frac{\snr^{2}}{(\snr + 1)^{2}\rdim} \Big)^{t}\distsqr{\thetait{0}} +
c_{2}\cdot\frac{(\snr + 1)^{2}\rdim}{\snr^{2}\inisig^{4}}\max
\Big\{\frac{\kdim\log\usedim}{\numobs},
\frac{\kdim^{2}\log^{2}\usedim}{\numobs^{2}} \Big\}.
\end{align*}
\end{cor}
\noindent
See Section~\ref{SecProofCorSparsePCA} for the proof of this
corollary.

\begin{rem}
It is noteworthy that $ \LossTil(\pltheta)$ is in fact \emph{globally
  concave} in $\pltheta$.  In order to see this fact, consider the
scalar case with $\usedim= \rdim= 1$, where $\LossTil(\pltheta) =
-C\pltheta^{2}$ for some $C >0$.
\end{rem}
The error rate $\epsnum$ is in fact minimax optimal (up to a
logarithmic factor) with respect to $\numobs$, $\usedim$, $\kdim$ as
well as the rank $\rdim$; for instance, see the
paper~\cite{vu2013minimax,cai2013rank}, where the upper bound is
achieved using computationally intractable estimators. Similar error
rates are obtained in~\cite{cai2013adaptive,ma2013iterative,wang2014pca} using
more sophisticated algorithms, but under a scaling of the sample
size---in particular one that is quadratic in sparsity (see
below)---that allows for a good initialization.


\paragraph{Initialization:}

The above results require an initial solution $\thetait{0}$ with
$\dist{\thetait{0}} < 1$.  Under the spiked covariance model and given
a sample size $\numobs \succsim \kdim^2 \log \usedim$, such solution
can be found by the diagonal thresholding
method~\cite{JohLu09,ma2013iterative}.  Here the quadratic dependence
on $\kdim$ is related to a computational
barrier~\cite{berthet2013lowerSparsePCA}, and thus may not improvable
using a polynomial-time algorithm.

\begin{figure}
\centering
\begin{tabular}{ccc}
\includegraphics[scale=0.5, clip, trim = 0 0 0 0]{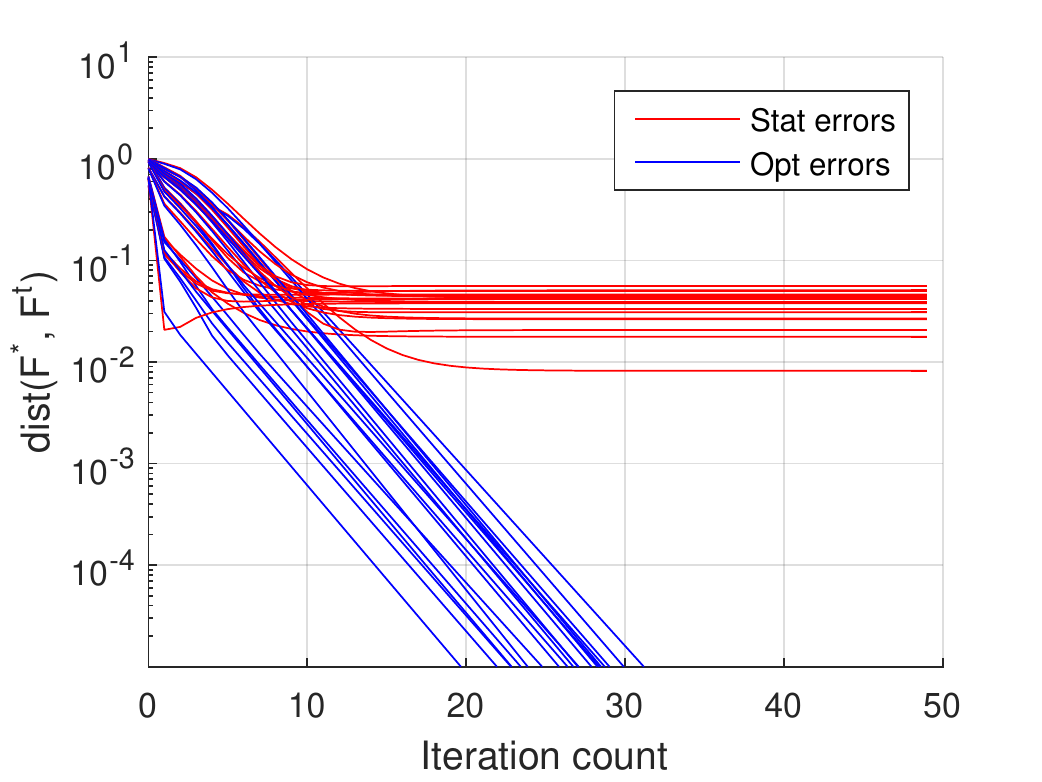} &
\includegraphics[scale=0.5, clip, trim = 0 0 0 0]{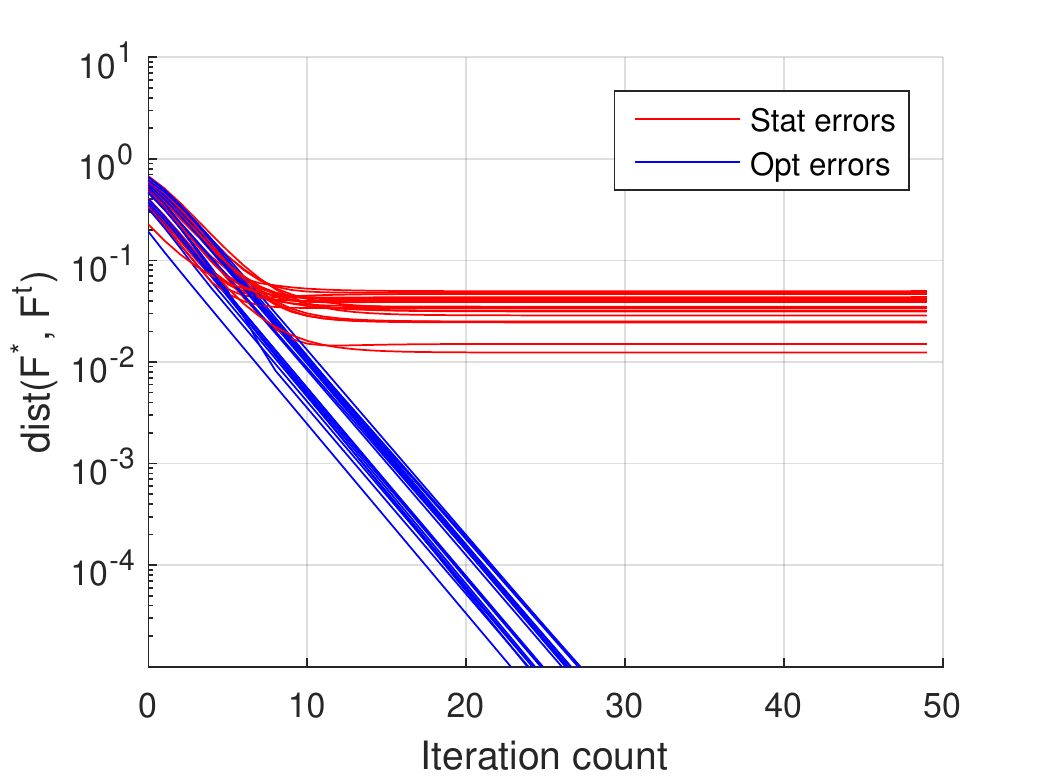} &
\includegraphics[scale=0.5, clip, trim = 0 0 0 0]{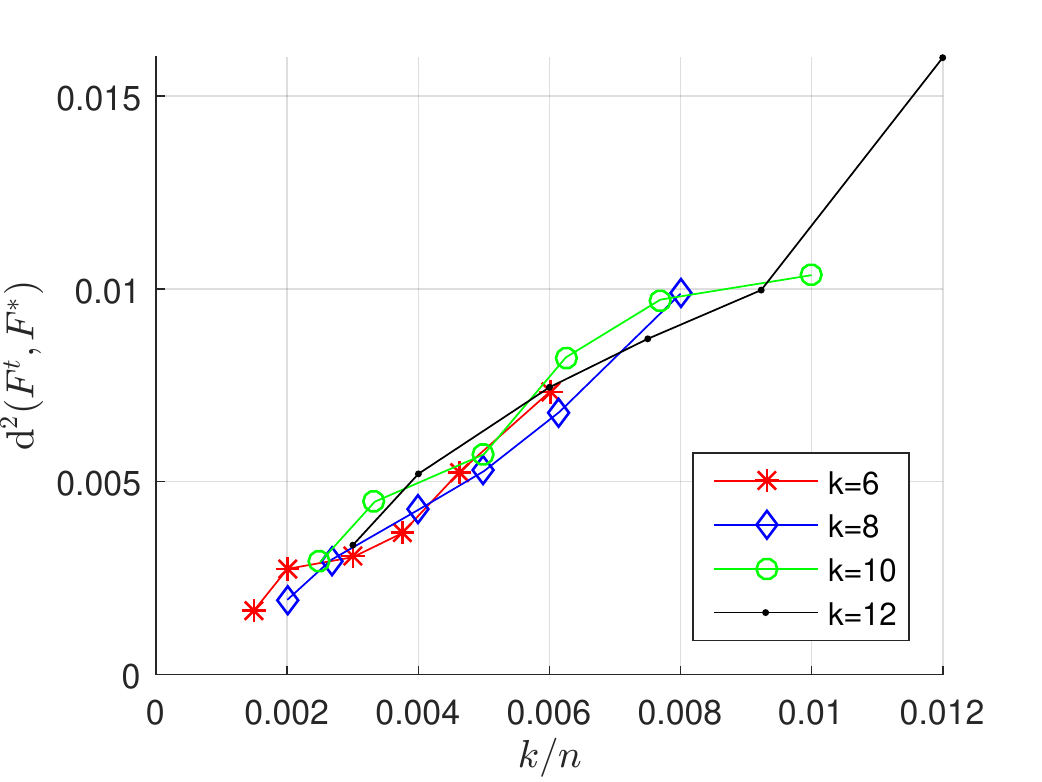}\\
(a) & (b) & (c)
	\end{tabular} 
	\caption{Simulation results for sparse PCA. Panel (a); plots
          of optimization error $ \textup{d} (\thetait{t}, \thetait{T}
          ) $ and statistical error $ \dist{\thetait{t}} $ versus the
          iteration number $ t$, using diagonal thresholding
          initialization.  Panel (b): same plots using perturbation
          initialization. For both panels (a) and (b), simulations are
          performed using $ \usedim = 5000$, $ \rdim = 1 $, $\kdim=5$,
          $ \snr = 4 $ and $ \numobs = 4000 $.  Panel (c): plot of
          estimation error $ \dist{\widehat{\pltheta}} $ versus $
          \frac{\kdim}{\numobs} $, for different values of $
          (\kdim,\numobs) $ using diagonal thresholding
          initialization. Each point represents the average over $20$
          random instances. The simulation is performed using $
          \usedim = 5000 $, $ \rdim = 1 $ and $ \snr = 4 $.}
	\label{FigSPCA}
\end{figure}

\paragraph{Computation:}

The algorithm requires projection onto the intersection of the
spectral norm and $\ell_{1}/\ell_{2}$ norm balls. In the rank one case
($\rdim= 1$), it reduces to projecting to the intersection of the
vector $\ell_{2}$ and $\ell_{1}$ balls, which can be done
efficiently~\cite{su2012projection}.  In the general case with $\rdim
> 1$, it can be done by alternating projection.  The speed of
convergence depends on the eigengap $\snr$, exhibiting similarity to
the standard power method for finding eigenvectors.  

\paragraph{Simulations:}  We performed experiments under the same
general set-up as the matrix completion (see the discussion
surrounding equation~\eqref{EqnSineApprox}). For sparse PCA, we
generated random ensembles of problems by fixing the rank $\rdim = 1$,
and choosing a random unit-norm $\thetastar \in \real^{\usedim}$
supported on $ \kdim $ randomly chosen coordinates.  Using this random
vector, we formed the spiked covariance matrix $\CovMat$ with top
eigenvector $ \thetastar $ and SNR $ \snr$. We considered two
approaches for initialization: (a) diagonal thresholding as described
in the papers~\cite{JohLu09,ma2013iterative}, and (b) choosing
$\thetait{0}$ to be the perturbed version $\thetait{0}_{\Tset} =
\thetastar_{\Tset} + \frac{1}{\sqrt{2}} E_1 $ and $
\thetait{0}_{\TsetComp} = \thetastar_{\TsetComp} + \frac{1}{\sqrt{2}}
E_2 $, where $ \Tset =\text{support}(\thetastar)$ and $ E_1 $ and $
E_2 $ are random unit norm vectors with the appropriate
dimensions. The step size is fixed at $ \stepit{t} \equiv
\frac{0.5\snr}{(\snr+1)^2} $.

Panels (a) and (b) of Figure~\ref{FigSPCA} show the convergence rates
of the optimization and statistical error using these two different
types of initializations.  Consistent with our theory, we witness an
initially geometric convergence in terms of statistical error followed
by an error floor at the statistical precision.  In panel (c), we
study the scaling of the estimation error.  Our theory predicts that
given a suitable initialization and sample size $ \numobs \succsim
\kdim \log\usedim $, then with high probability the output $
\widehat{\pltheta} $ satisfies
\begin{align*}
\distsqr{\widehat{\pltheta}} \precsim \frac{(\snr+1)^2 \rdim}{\snr^2}
\cdot \frac{\kdim \log \usedim}{\numobs}.
\end{align*}
Therefore, with the triplet of parameters $(\usedim, \rdim, \snr)$
fixed, the error $ \distsqr{\widehat{\pltheta}}$ should grow
proportionally with the ratio $\frac{\kdim}{\numobs}$, a prediction
that is confirmed in Figure~\ref{FigSPCA}(c).


\subsection{Planted densest subgraph}
\label{SecClustering}

The planted densest subgraph problem is a generalization of the
planted clique problem; it can be viewed as a single cluster (or rank
one) version of the more general planted partition problem. For a
collection of $\usedim$ vertices, there is an unknown subset of size
$\csize$ which forms a cluster. Based on this cluster and two
probabilities $\pin > \qout$, a random symmetric matrix $A \in
\{0,1\}^{\usedim\times \usedim}$, which we think of as the adjacency
matrix of the observed graph, is generated in the following way:
\begin{itemize}
\item for each pair of vertices $i,j$ in the cluster, $A_{ij} = 1$ with
probability $p$, and zero otherwise.
\item for all other pairs of vertices, $A_{ij} = 1$ with probability $q$,
and zero otherwise.
\end{itemize}
Let $\thetastar \in  \big\{0,1 \big\}^{\usedim}$ be the cluster membership
vector: i.e., $\pltheta_{j}^{*} = 1$ if and only if vertex $j$ belongs
to the cluster.

A previous approach is to recover the cluster matrix $\ThetaStar= \OUTER{\thetastar}$
by solving a particular SDP, derived as a relaxation of the MLE. Let
$\Shift \defn \Adj - \frac{\pin + \qout}{2}\OneMat$ be a shifted version
of the adjacency matrix, where $\OneMat$ is the $\usedim\times \usedim$
all one matrix. Consider the semidefinite program
\begin{align}
\min_{\PlTheta \in \SymMat{\usedim}} \Big\{ - \trinprod{\Shift}{\PlTheta} \Big\} \qquad \mbox{such that \ensuremath{\PlTheta\succeq0}, \ensuremath{\sum_{i,j}\PlTheta_{ij} = k^{2}} and  \ensuremath{\PlTheta\in[0,1]^{\usedim\times \usedim}}.} \label{EqnClusteringSDP}
\end{align}
It is known~\cite{chen2012sparseclustering} that with probability
at least $1 - \usedim^{2}$, the true cluster matrix $\ThetaStar$ is
the unique optimal solution to this program when
\begin{equation}
\frac{(\pin - \qout)^{2}}{\pin} \ge c_{1} \Big( \frac{\log\usedim}{\csize} + \frac{\usedim}{\csize^{2}} \big), \label{EqnClusteringCond}
\end{equation}
for some universal constant $c_{1} > 0$. When $\pin= 1= 2\qout$, this
condition reduces to the well-known $\csize \succsim \sqrt{\usedim}$
tractability region for the planted clique problem~\cite{alon1998hiddenClique}.

Alternatively, we may solve the factorized formulation by projected
gradient decent~\eqref{EqnGenProjGradOpt}, as applied to the problem
\begin{align*}
\LossTil( \pltheta ) = \trinprod{ - \Shift}{\OUTER{\pltheta}},
\qquad \GenCons = \{ \pltheta \,\mid\, \pltheta \in [0,1]^{\usedim},\,
       {\textstyle \sum_{i=1}^{\usedim} \pltheta_{i} = \csize \}.}
\end{align*}
This setting is a $\rdim= 1$ special case of our general framework.
In this case $\eclass= \{\pm\thetastar\}$ contains only two elements, and can be verified to be $ \ThetaStar $-faithful. 

We now ready to apply our general theory to this problem. As we show in proof of Corollary~\ref{CorClustering}, if the model parameters satisfy the condition~(\ref{EqnClusteringCond}), then with probability at least $1 -
\usedim^{-3}$,  the loss function $\LossTil $ satisfies the local 
descent and smoothness conditions and the relaxed Lipschitz
condition~(\ref{EqnLtilWeakLipCon}) with parameters
\begin{align*}
\smallrad = \frac{2}{5}\sqrt{\csize}, \quad \curv = \frac{1}{20}(\pin
- \qout)\csize, \quad \smoo = 12(\pin - \qout)\csize \quad \text{and}
\quad \epsnum= 0.
\end{align*}
Using this fact, we have the following corollary of Theorem~\ref{ThmGeneralSmooth}. We state it assuming that the condition condition~(\ref{EqnClusteringCond}) holds
\begin{cor}
\label{CorClustering} 
Under the previously stated conditions, given 
an initial vector $\thetait{0} \in \GenCons \cap
\Ballstar{\frac{1}{5}\sqrt{\csize}}$, the projected gradient iterates
$\{\thetait{t}\}_{t= 1}^{\infty}$ with step size $\stepit{t} = c_{2}
\frac{1}{(\pin - \qout)\csize}$ satisfy the bound
\begin{align*}
\distsqr{\thetait{t}} & \leq \Big( 1-c_{3}
\Big)^{t}\distsqr{\thetait{0}}.
\end{align*}
\end{cor}
\noindent See Section~\ref{SecProofCorClustering} for the proof of
this claim.

The corollary guarantees exact recovery of $ \ThetaStar $ when $ t \to
\infty $. The condition~\eqref{EqnClusteringCond} matches the best
existing results; see e.g.,~\cite{chen2012sparseclustering} and the
references therein.


\begin{figure}
\centering
\begin{tabular}{cc}
\includegraphics[scale=0.5, clip, trim = 0 0 0 0]{\figpath
  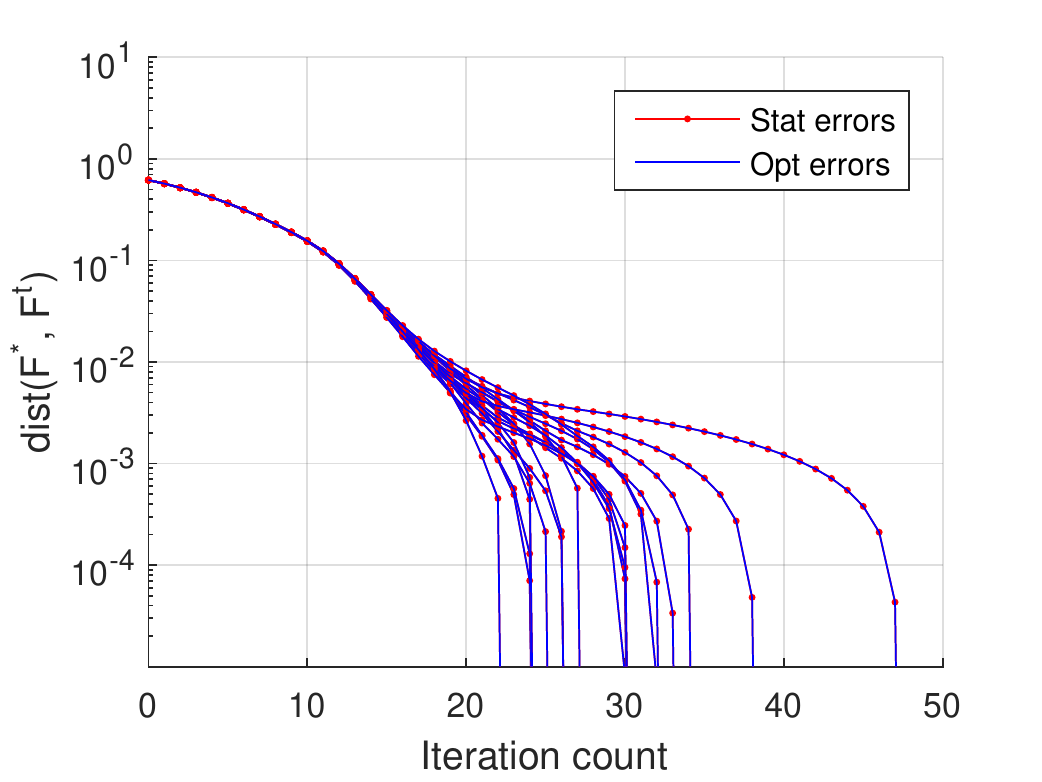} &
\includegraphics[scale=0.5, clip, trim = 0 0 0
  0]{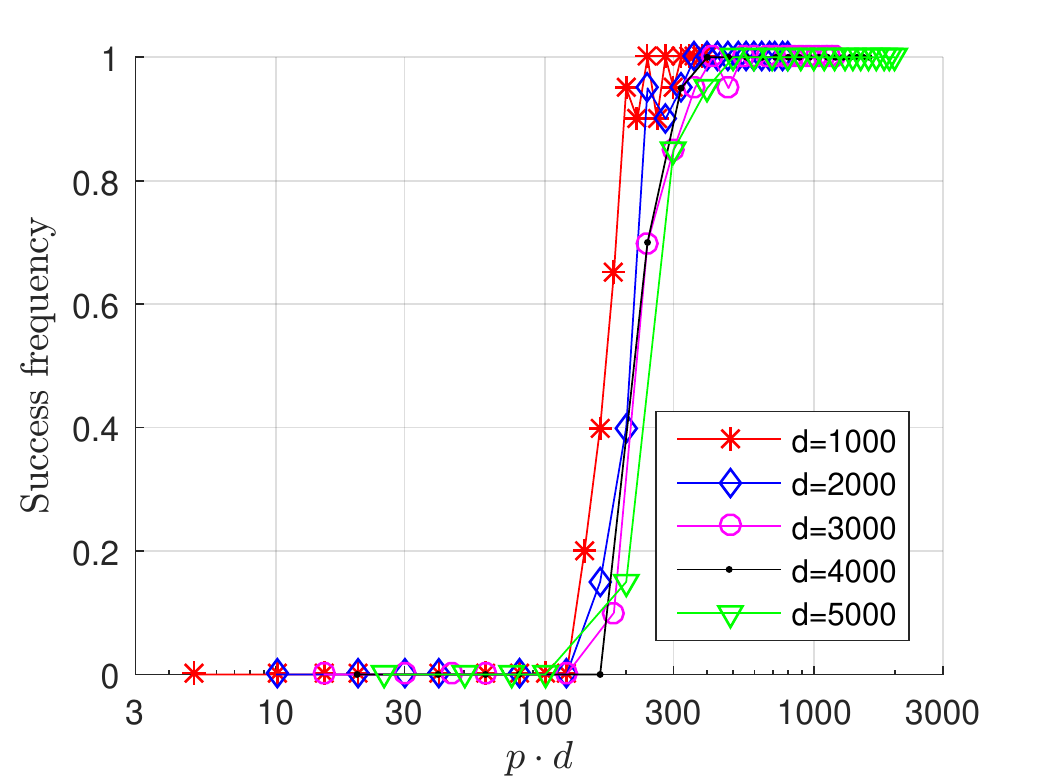} \\
(a) & (b)
	\end{tabular} 
	\caption{Simulations for planted densest subgraph. Panel (a):
          plots of optimization error $ \textup{d} (\thetait{t},
          \thetait{T} ) $ and statistical error $ \dist{\thetait{t}} $
          versus the iteration number $ t$, using SVD-based
          initialization. The simulation is performed using $ \usedim
          = 8000 $, $ \kdim =2000 $, $ \pin = 0.13 $ and $\qout = 0.05
          $.  Panel (b): plot of the probability of successful exact
          recovery of $ \thetastar $ versus $ \pin \usedim $, for
          different values of $ (\usedim, \pin) $ using SVD-based
          initialization.  We declare exact recovery if $
          \dist{\widehat{\pltheta}} \le 2\times 10^{-3} $, and each
          point represents frequency of exact recovery over $ 20 $
          random instances.  The simulation is performed with $ \qout
          = \frac{\pin}{4} $ and $ \kdim = \frac{\usedim}{2} $.  }
	\label{FigPlanted}
\end{figure}

\paragraph{Initialization}

Set $\thetait{0}$ to be the top left singular vector of $\Adj -
\qout\OneMat$ projected onto the set $ \GenCons $. Note that
$\thetastar$ is a left singular vector of the matrix $\Exs[\Adj -
  \qout\OneMat]$ corresponding to the only non-zero singular
value~$(\pin - \qout)\kdim$.  Under the
condition~(\ref{EqnClusteringCond}), Proposition~\ref{PropBandeira}
ensures that $\opnorm{(\Adj - \qout\OneMat) - \Exs[\Adj -
    \qout\OneMat]} \le \frac{1}{4}(\pin - \qout)\kdim$ with
probability at least $1 - \usedim^{-3}$. On this event, applying
Wedin's $\sin\Theta$ theorem~\cite{golub_matrixcomp} guarantees that
$\thetait{0}$ satisfies the requirement in
Corollary~\ref{CorClustering}.

\paragraph{Computation:}

The set $\GenCons$ is the intersection of a hyperplane and a box in $
\real^{\usedim} $, so the associated projection $\ProjGenCons$ can be
computed in time $ \order(\usedim)
$~\cite{maculan2003project}. Computing the gradient $\nabla \LossTil
(\pltheta) = - 2\Shift\pltheta$ only requires matrix-vector
multiplication with the matrix $\Shift$, which is the sum of a
rank-$1$ matrix and the (usually sparse) graph adjacency matrix
$A$. In contrast, solving the SDP in equation~\eqref{EqnClusteringSDP}
using ADMM requires multiple full SVD of dense matrices even when the
graph is sparse.

\paragraph{Simulations:}
We performed experiments under the same general set-up as the matrix
completion (see the discussion surrounding
equation~\eqref{EqnSineApprox}).  The $ 0 $-$ 1 $ cluster indicator
matrix $ \thetastar \in \real^{\usedim\times 1} $ is supported on $
\csize $ coordinates, and we sampled the graph adjacency matrix $ \Adj
$ from the planted densest subgraph model with edge probabilities $
(\pin, \qout) $ and cluster size $ \kdim $. The initial matrix $
\thetait{0} $ is obtained using the SVD-based procedure described in
Section~\ref{SecClustering}.  The step size is fixed at $ \stepit{t}
\equiv \frac{0.1}{(\pin - \qout)\kdim }$.  Panel (a) of
Figure~\ref{FigPlanted} shows plots of the optimization and
statistical errors versus the iteration number; consistent with
Corollary~\ref{CorClustering}, these iterates converge at least
geometrically.

In terms of the scaling of the sample size required for exact
recovery, we know that if $ \kdim \precsim \frac{\usedim}{\log\usedim}
$, then the fixed point of the algorithm $ \widebar{\pltheta}$ will be
equal to $\thetastar$ with high probability provided that $\frac{(\pin
  - \qout)^2}{\pin} \succsim \frac{\usedim}{\kdim^2}$.  In particular,
see equation~\eqref{EqnClusteringCond}.  Therefore, with $ \qout =
\frac{\pin}{3}$ and $ \kdim = \frac{\usedim}{2} $, exact recovery of $
\thetastar $ can be achieved with probability close to one as soon as
$ \pin \usedim $ is above a constant threshold.  This theoretical
prediction is confirmed in panel (b) of Figure~\ref{FigPlanted}.


\subsection{One-bit matrix completion}
\label{SecOB}

Let us now turn to an extension of the standard (linear) matrix
completion model studied in Section~\ref{SecMC}.  It provides a more
challenging problem to analyze, and our general theory provides
(to the best of our knowledge) the first known polynomial-time
algorithm for achieving the minimax rate in the case of rank
$\rdim$ matrices.

In order to set up the problem, suppose that $\thetastar \in
\real^{\usedim\times \rdim}$ is an orthonormal matrix and has
incoherence parameter $\inco$ as previsouly defined in
equation~\eqref{EqnIncoherence}.  Given a set $\Obs \subseteq
[\usedim] \times [\usedim]$ of observed elements, a noise parameter
$\noisestd > 0$ and a differentiable function $f : \real \mapsto
[0,1]$ with Lipschitz derivative, we observe a binary symmetric matrix
$\Yout \in \{-1,1\}^{\usedim\times \usedim}$ such that for each $(i,j)
\in \Obs$ with $ i\ge j $,
\begin{align*}
\Yout_{ij} &=
\begin{cases}
1, & \text{with probability \ensuremath{f(\ThetaStar_{ij}/\noisestd)}}, \\
-1, & \text{with probability \ensuremath{1-f(\ThetaStar_{ij}/\noisestd)}}.
\end{cases}
\end{align*}
We further assume that the observation set $\Obs$ is symmetric and
generated by the Bernoulli model with parameter $\pobs$, that is,
$\mprob\big( (i,j), (j,i) \in \Omega \big)= \pobs$ independently for
each $(i,j)$ with $ i\ge j $ The goal is to estimate $\ThetaStar$
given the binary observations $\Yout$.  Examples of the function $f$
include the \emph{logistic model} with $f(x)= \frac{\exp(x)}{1 +
  \exp(x)}$; the \emph{probit model} with $f(x)= \Phi(x)$, where
$\Phi(x)$ is the cumulative distribution function of a standard
Gaussian; and the \emph{Laplacian model} in which $f$ is the
cumulative distribution function of Laplacian($0,1$): variable.  See
the papers~\cite{cai2013maxnormOneBitMC,davenport2014_onebitMC} for
more details on these choices.

For a given $f$, consider thee negative log-likelihood of~$\PlTheta$,
given by
\begin{align}
\EmpLoss(\PlTheta) & = -2\sum_{(i,j) \in \Obs} \Big[ \frac{1 +
    \Yout_{ij}}{2} \log f(\PlTheta_{ij}/\noisestd) + \frac{1 -
    \Yout_{ij}}{2} \log \big( 1-f(\PlTheta_{ij}/\noisestd) \big) \Big]
\nonumber \\
& = - \trinprod{\ProjObs(\OneMat + \Yout)}{\log f(\PlTheta/\noisestd)}
- \trinprod{\ProjObs(\OneMat - \Yout)}{\log \big(
  1-f(\PlTheta/\noisestd) \big) }, \label{EqnLossOB}
\end{align}
where $\OneMat$ is the $\usedim \times \usedim$ all one matrix,
$\circ$ denotes the Hadamard product, and functions are applied to a
matrix element-wise. As in matrix completion, we use the set $\GenCons
= \big \{\pltheta \,\mid\, \twoinfnorm{\pltheta} \le
\sqrt{\frac{2\inco}{\usedim}} \fronorm{\thetait{0}} \big \}$ Note that
the gradient of the loss function is given by
\begin{align*}
\nabla_{\PlTheta} \EmpLoss(\PlTheta) & = - \frac{1}{\noisestd}
\ProjObs \Big[ \frac{ f'(\PlTheta/\noisestd) \circ (\Yout -
    2f(\PlTheta/\noisestd) + \OneMat)}{f(\PlTheta/\noisestd) \circ (1
    - f(\PlTheta/\noisestd))} \Big],
\end{align*}
where the fraction are also element-wise.

Since the function $f$ is differentiable with a Lipschitz derivative
$f'$, Rademacher's theorem guarantees that the second derivative $f''$
is defined almost everywhere.  Our corollary depends function $f$
through the following two quantities, defined for each $a > 0$:
\begin{equation}
\label{EqnFlatness}
\begin{aligned}\obupper{a}  \defn  & \max \Big\{ \sup_{\abs{x} < a} 
\frac{\abs{f'(x)}}{f(x)(1-f(x))}, \; \sup_{\abs{x} < a}
\frac{f'(x)^{2}}{f(x)^{2}(1-f(x))^{2}}, \; \sup_{\abs{x} < a}
\frac{\abs{f''(x)}}{f(x)(1-f(x))} \Big\}, \quad \mbox{and} \\
\oblower{a} \defn & \sup_{\abs{x} < a}
\frac{f(x)(1-f(x))}{f'(x){}^{2}}.
\end{aligned}
\end{equation}
These quantities are similar to those in the
paper~\cite{davenport2014_onebitMC}, along with the additional control
over the second derivative~$f''$ required for proving (fast) geometric
convergence.  We introduce the shorthand $\obsnr \defn
\frac{\inco\rdim}{\usedim\noisestd} =
\frac{\infnorm{\ThetaStar}}{\noisestd}$, which we think of as a
measure of SNR. In the constant SNR setting $\obsnr= \Theta(1)$, the
quantities $\obupper{4\obsnr}$ and $\oblower{4\obsnr}$ are positive
universal constants independent of the other model parameters
$\usedim,\rdim,\pobs$ etc.

We now apply Theorem~\ref{ThmGeneralSmooth} to the one-bit matrix completion problem. Set 
\[
\smallrad= c_{1}\max
\Big\{1,\frac{1}{\oblower{4\obsnr}\obupper{4\obsnr}} \Big\}, \quad
\alpha= c_{2} \frac{\pobs}{\oblower{4\obsnr}\noisestd^{2}}, \quad
\LipCon= \beta= c_{3}
\frac{\obupper{4\obsnr}\pobs\inco\rdim}{\noisestd^{2}} \quad
\text{and} \quad \epsnum=
c_{4}\noisestd\obupper{4\obsnr}\oblower{4\obsnr}(1 +
\obsnr){\textstyle \sqrt{\frac{\usedim\rdim}{\pobs}}}.
\] As we shown in the proof of Corollary~\ref{CorOB}, if the initial matrix satisfies the condition  $\dist{\thetait{0}} \le1
- \sqrt{1 - \smallrad}$, then the set $\GenCons$ is
$\ThetaStar$-faithful. Moreover, if the expected sample size satisfies the bound $\numobs= \pobs\usedim^{2} \ge c_{5} \max\{\inco\rdim\usedim\log\usedim, \usedim\log^{2}\usedim,\inco^{2}\rdim^{2}\usedim\}$ and is large enough to ensure $\epsnum < \frac{1}{20}(1 - \sqrt{1 - \smallrad})$, then with probability at least $ 1-c_6 \usedim^{-3} $, the loss function $ \LossTil $ associated with~\eqref{EqnLossOB} satisfies the local descent, Lipschitz and smoothness conditions with parameters $ \smallrad $, $ \curv$, $\LipCon$, $\smoo $ and $ \epsnum $ given above.
Using these facts, we obtain the following guarantee, which we state assuming that the sample size $ \numobs $ satisfies the above conditions.
\begin{cor}
\label{CorOB} 
Under the previously stated conditions,  if we are
given an initial matrix $\thetait{0}$ with $\dist{\thetait{0}} \le1
- \sqrt{1 - \smallrad}$, then with probability at least $ 1-c_6 \usedim^{-3} $,  the gradient descent iterates
$\{\thetait{t}\}_{t= 1}^{\infty}$ with step size $\stepit{t} = c_{7}
\frac{\noisestd^{2}}{\oblower{4\obsnr}\obupper{4\obsnr}^{2}\pobs\inco\rdim}$ satisfy the bound
\begin{align*}
\distsqr{\thetait{t}} & \le \Big( 1-c_{8}
\frac{1}{\oblower{4\obsnr}^{2}\obupper{4\obsnr}^{2}\inco\rdim}
\Big)^{t}\distsqr{\thetait{0}} +
c_{9}\noisestd^{2}\obupper{4\obsnr}^{2}\oblower{4\obsnr}^{2}(1 +
\nu)^{2} \frac{\usedim\rdim}{\pobs}.
\end{align*}
\end{cor}
\noindent See Section~\ref{SecProofCorOB} for the proof of this claim.

In order to interpret the above result, let us consider the setting
with a constant SNR $\obsnr= \Theta(1)$, in which case $\noisestd=
\frac{\inco\rdim}{\usedim\obsnr}\asymp\infnorm{\OUTER{\thetastar}}$.
Corollary~\ref{CorOB} guarantees that given an initial matrix
$\thetait{0}$ within a constant radius of $\thetastar$, the projected
gradient descent converges geometrically and has per-entry error
\begin{align}
\label{EqnOB_consequence}
\frac{1}{\usedim^{2}} \fronorm{\OUTER{\thetait{\infty}} -
  \OUTER{\thetastar}}^{2} & \le \frac{3}{\usedim^{2}}
\opnorm{\thetastar}^{2}\distsqr{\thetait{\infty}} \precsim
\frac{\usedim\rdim}{\numobs}\noisestd^{2}\asymp\frac{\usedim\rdim}{\numobs}\infnorm{\OUTER{\thetastar}}^{2}.
\end{align}
This bound has the same form as that in Section~\ref{SecMC} for
standard matrix completion, with the important difference that it is
an essentially \emph{multiplicative} bound where the pre-factor
depends on the SNR~$\obsnr$.

\begin{figure}
\centering
\begin{tabular}{cc}
\includegraphics[scale=0.5, clip, trim = 0 0 0 0]{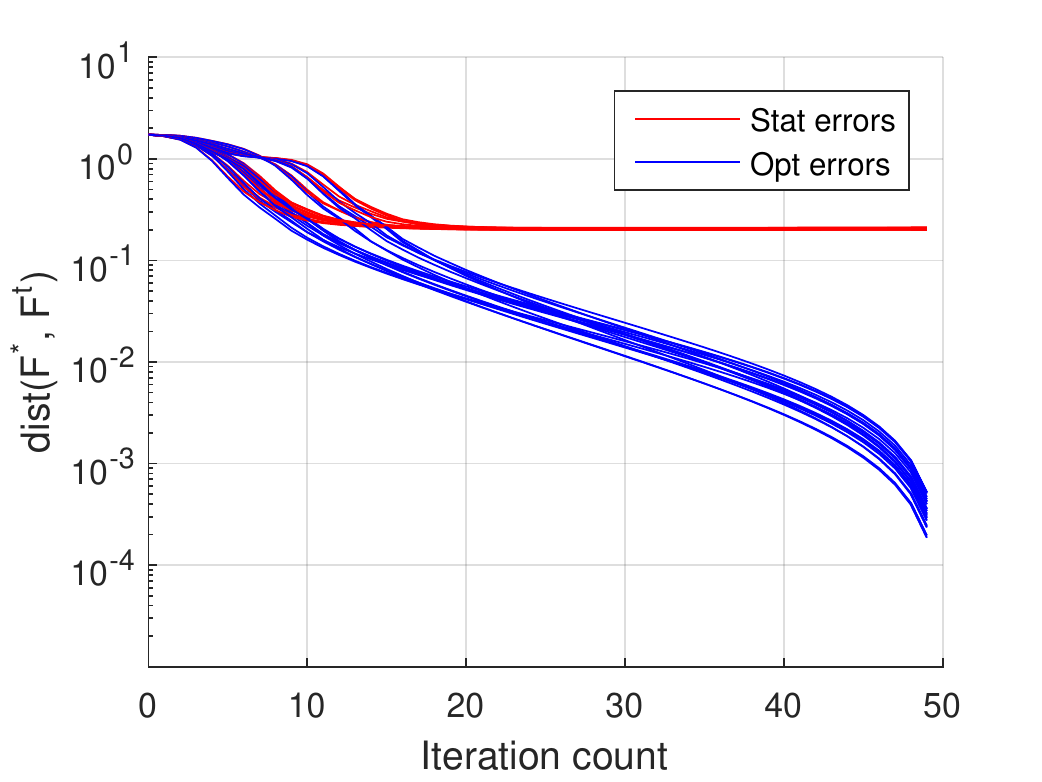} & \includegraphics[scale=0.5, clip, trim = 0 0 0
  0]{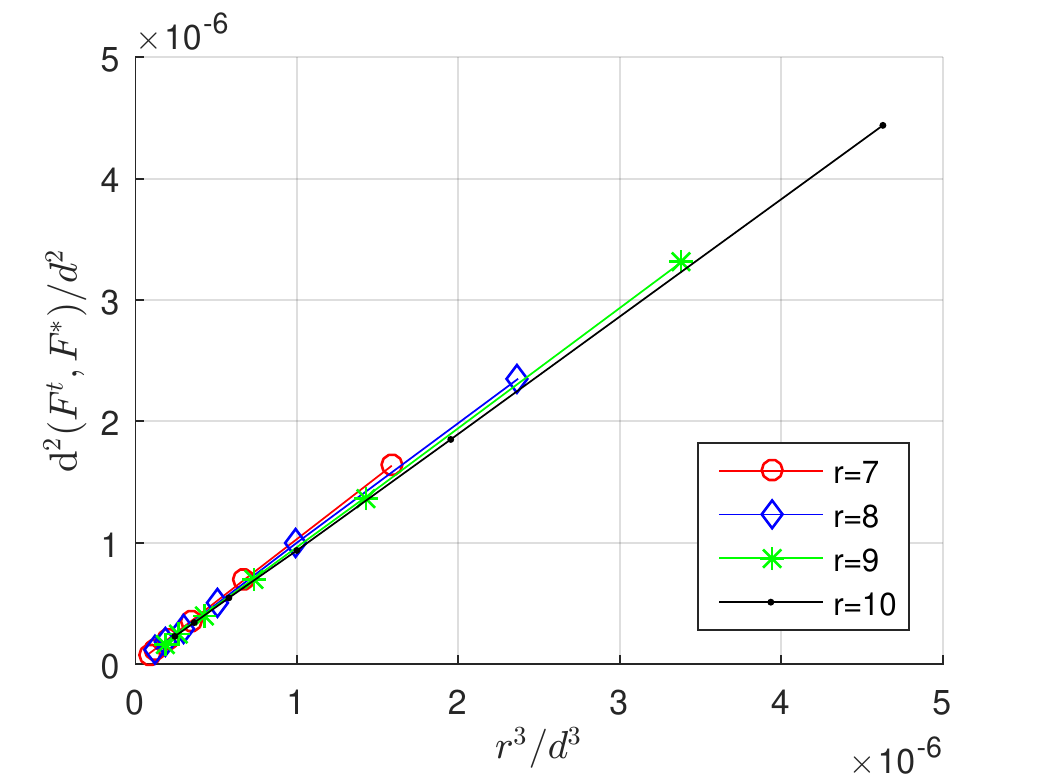} \\
(a) & (b)
\end{tabular} 
\caption{Simulation results for one-bit matrix completion.  Panel (a);
  plots of OB optimization error $ \textup{d} (\thetait{t},
  \thetait{T} ) $ and statistical error $ \dist{\thetait{t}} $ versus
  the iteration number $ t$, using random initialization.  The
  simulation is performed using $ \usedim = 1000$, $\rdim = 3$ and $
  \pobs =0.5$.  Panel (b) plot of per-entry estimation error $
  \frac{1}{\usedim^2} \dist{\widehat{\pltheta}} $ versus $
  \frac{\rdim^3}{\usedim^3} $, for different values of $ (\usedim,
  \rdim) $ using random initialization. Each point represents the
  average over $ 20 $ random instances. The simulation is performed
  using $ \pobs=0.5 $ and $ \noisestd = \frac{0.5\rdim}{\usedim} $.}
	\label{FigOneBit}
\end{figure}

It is worth comparing our error bounds with previous results under the
setting $\nu= \Theta(1)$. One body of past
work~\cite{davenport2014_onebitMC,cai2013maxnormOneBitMC} has studied
the recovery of \emph{approximately low-rank} matrices with bounded
nuclear norm---that is, matrices whose vectors of singular values are
in the $\ell_{q}$ ball with $q= 1$.  This is a milder sparsity
assumption, and so leads to the slower error rate $O \Big(
\sqrt{\frac{\usedim\rdim}{\numobs}} \Big)$.  The result here applies
to exactly low-rank matrices ($q= 0$), and so leads to the faster rate
$\frac{\usedim\rdim}{\numobs}$.  Both of these scalings are be
minimax-optimal in the simpler linear
setting~\cite{negahban2010noisymc}.  On the other hand, Bhaskar et
al.~\cite{bhaskar2015oneBit} also analyze the case of exactly low-rank
matrices, but their algorithm relies on rank-constrained optimization
and does not have convergence guarantees in polynomial time.
Moreover, their error rate scales as
$\frac{\usedim^{7}\rdim^{3}}{\numobs^{4}}$, and thus has a worse
dependence on $\rdim$, $\usedim$ and $\numobs$ as compared to ours.

\paragraph{Initialization and time complexity:}

In theory, we can obtain a good initial solution $\thetait{0}$ by
solving one of the convex programs in the
papers~\cite{davenport2014_onebitMC,cai2013maxnormOneBitMC} followed
by a projection onto the set $\GenCons$. Since we only need the
initial error to be a constant, it suffices to have $\numobs \succsim
\usedim\rdim + \usedim\log\usedim$ observations.  In fact, in our
simulations, we find that a randomly chosen initial matrix
$\thetait{0}$ is often good enough (see
Figure~\ref{FigOneBit}(a)). Given such an initial solution, the
projected gradient iterates converges geometrically with a contraction
factor $1 - \frac{c_{8}}{\inco\rdim}$, so we need
$\order(\inco\rdim\log(1/\delta))$ iterations to compute a
$\delta$-accurate solution. Therefore, we can achieve the
$\order(\frac{\usedim\rdim}{\numobs})$ error rate in polynomial time;
to the best of our knowledge, this polynomial-time guarantee for
achieving the minimax-rate in the exact low-rank case is the first
such result in the literature.

\paragraph{Simulations:}
We performed experiments under the same general set-up as the matrix
completion (see the discussion surrounding
equation~\eqref{EqnSineApprox}).  The matrix $ \thetastar $ is random
orthonormal, and the observations are generated using the Bernoulli
model with observation probability $ \pobs $ and the standard Gaussian
CDF as the link function $ f $ with noise magnitude $ \noisestd =
\frac{2\rdim}{\usedim} $. The initial matrix $ \thetait{0} $ is
obtained by random initialization.  The step size is fixed at $
\stepit{t} \equiv \frac{0.5\noisestd^2}{\pobs} $.  Panel (a) of
Figure~\ref{FigOneBit} illustrates the geometric convergence
of the algorithm.

In terms of the scaling of the estimation error, with $ \noisestd =
\frac{2\rdim}{\usedim} $, $ \numobs = \pobs \usedim^2 $ and $ \pobs $
fixed, the per-entry error of the output $ \widehat{\pltheta} $
satisfies
\begin{align*}
\frac{1}{\usedim^2}\distsqr{\widehat{\pltheta}} \precsim
\frac{\usedim\rdim}{\numobs} \propto \frac{\rdim^3}{\usedim^3}
\end{align*}
with high probability;
cf.\ equation~\eqref{EqnOB_consequence}. Therefore, we should expect
that the squared error
$\frac{1}{\usedim^2}\distsqr{\widehat{\pltheta}} $ scales
proportionally with the ratio $\frac{\rdim^3}{\usedim^3}$, a
prediction that is confirmed in Figure~\ref{FigOneBit}(b).

\subsection{Low-rank and sparse matrix decomposition}
\label{SecMatrixDecomposition}

Recall from Section~\ref{SecMatrixDecompositionExample} the problem of
noisy matrix decomposition, in which we observe a noisy sum of the
form $\Yout= \OUTER{\thetastar} + \SStar + \noisevar$, where
$\noisevar$ is a symmetric noise matrix.  Our goal is to estimate
$\thetastar$, and in this section, we analyze a version of this model
in which the factor matrix $\thetastar \in \real^{\usedim\times
  \rdim}$ has equal eigenvalues and incoherence parameter $\inco$ as
defined in equation~\eqref{EqnIncoherence}, and the perturbing matrix
$\SStar \in \SymMat{\usedim}$ is element-wise sparse.

One line of work concerns the setting where the non-zero entries of
$\SStar$ are randomly
located~\cite{candes2009robustPCA,chen2011LSarxiv}, whereas another
line of work focuses on deterministic
models~\cite{chandrasekaran2011siam,Hsu2010RobustDecomposition,
  agarwal2012decomposition}.  We focus on one version of the
deterministic setting, in which each row/column of the matrix $\SStar$
has at most $\kdim$ non-zero entries, whose locations and values are
otherwise arbitrary.  In light of keeping the presentation as simple
as possible, we assume here the values of
$\vecnorm{\SStar_{i\cdot}}1$, the $\ell_{1}$ norm of each row
of~$\SStar$ are known.\footnote{This is unrealistic and could be
  relaxed, albeit at the price of more involved analysis of the
  Lagrangian version instead of the constrained version.}

Using the nuclear norm and $\ell_{1}$ norms as surrogates for rank and
sparsity (respectively), the constrained version of a popular convex
relaxation approach is based on the SDP
\begin{align*}
\min_{\PlTheta \in \SymMat{\usedim}} \Big \{ \frac{1}{2}\big( \min_{\PlS \in
  \ConsS} \fronorm{\Yout - (\PlTheta + \PlS)}^{2} \big) +
\lambda \nucnorm{\PlTheta} \Big \},
\end{align*}
where $\ConsS \defn \{ \PlS \in \real^{\usedim\times \usedim} \, \mid
\, \onenorm{\PlS_{i\cdot}} \le \onenorm{\SStar_{i\cdot}},i=
1,2,\ldots,\usedim\}$.	Alternatively, we may drop the nuclear norm
regularizer and solve the factorized formulation by projected gradient
descent, as applied to the problem
\begin{align*}
\EmpLoss(\PlTheta) = \frac{1}{2}\min_{\PlS \in \ConsS} \fronorm{\PlTheta + \PlS -
  \Yout}^{2}, \qquad \GenCons = \Big\{ \pltheta \,\mid\,
\twoinfnorm{\pltheta} \le \sqrt{\frac{2\inco}{\usedim}}
\fronorm{\thetait{0}} \Big\}.
\end{align*}
Note that $\EmpLoss(\PlTheta)$ is the squared Euclidean distance
between the point $\Yout - \PlTheta$ and the closed convex
set~$\ConsS$.  Therefore, the function $\EmpLoss$ is convex and has
gradient
\begin{align*}
\nabla_{\PlTheta} \EmpLoss(\PlTheta) & = \PlTheta + \ProjConsS(\Yout -
\PlTheta) - \Yout.
\end{align*}

We now derive a guarantee for this problem using our
Theorem~\ref{ThmGeneralSmooth}. As we show in the proof of
Corollary~\ref{CorLS}, if the initial matrix $\thetait{0}$ satisfies
$\dist{\thetait{0}} \le \frac{1}{5}$, then the constraint set $
\GenCons$ is $ \ThetaStar $-faithful. Moreover, if the sparsity of the
matrix $\SStar$ satisfies $\frac{\inco\rdim\kdim}{\usedim} \le c_{1}$
and the noise matrix satisfies $\fronorm{\noisevar} \le c_{2}
\opnorm{\thetastar}$, then the loss function and the feasible set
satisfy the local descent, Lipschitz and smoothness conditions with
parameters
\begin{equation*}
\smallrad = \frac{3}{5} \opnorm{\thetastar}, \quad \curv = \frac{1}{10}
\opnorm{\thetastar}^{2}, \quad \LipCon = \smoo =
48\opnorm{\thetastar}^{2} \quad \text{and} \quad \epsnum=
128\frac{\fronorm{\noisevar}}{\opnorm{\thetastar}}.
\end{equation*}
Using these facts, we have the following guarantee, which is stated assuming that the matrices $ \SStar $ and $ \noisestd $ satisfy the assumptions above.
\begin{cor}
\label{CorLS} 
Under the previously stated conditions, given any initial matrix
$\thetait{0}$ satisfying the bound \mbox{$\dist{\thetait{0}} \le
  \frac{1}{5}$,} the gradient iterates $ \big\{ \thetait{t} \big\}
_{t= 1}^{\infty}$ with step size $\stepit{t} = c_{3}
\frac{1}{\opnorm{\thetastar}^{2}}$ satisfy the bound
\begin{align}
\distsqr{\thetait{t}} & \le \big( 1-c_{4} \big) ^{t}
\distsqr{\thetait{0}} + c_{5}
\frac{\fronorm{\noisevar}^{2}}{\opnorm{\thetastar}^{2}}.
\end{align}
\end{cor}

\noindent
See Section~\ref{SecProofCorLS} for the proof.

The condition $\frac{\inco\rdim\kdim}{\usedim} \le c_{1}$ matches the
best existing results
in~\cite{Hsu2010RobustDecomposition,chen2011LSarxiv} for the
deterministic setting of matrix decomposition. As a passing
observation, the above results can be applied to matrix completion
with \emph{adversarial} missing entries---by arbitrarily filling in
the missing entries and treating them as sparse corruption,
Corollary~\ref{CorLS} guarantees recovery when each row/column has at
most $ \kdim \le c_1 \frac{\usedim}{\inco \rdim} $ missing entries,
whose locations can be arbitrary.

\begin{figure}
\centering
\begin{tabular}{ccc}
	\includegraphics[scale=0.5, clip, trim = 0 0 0 0]{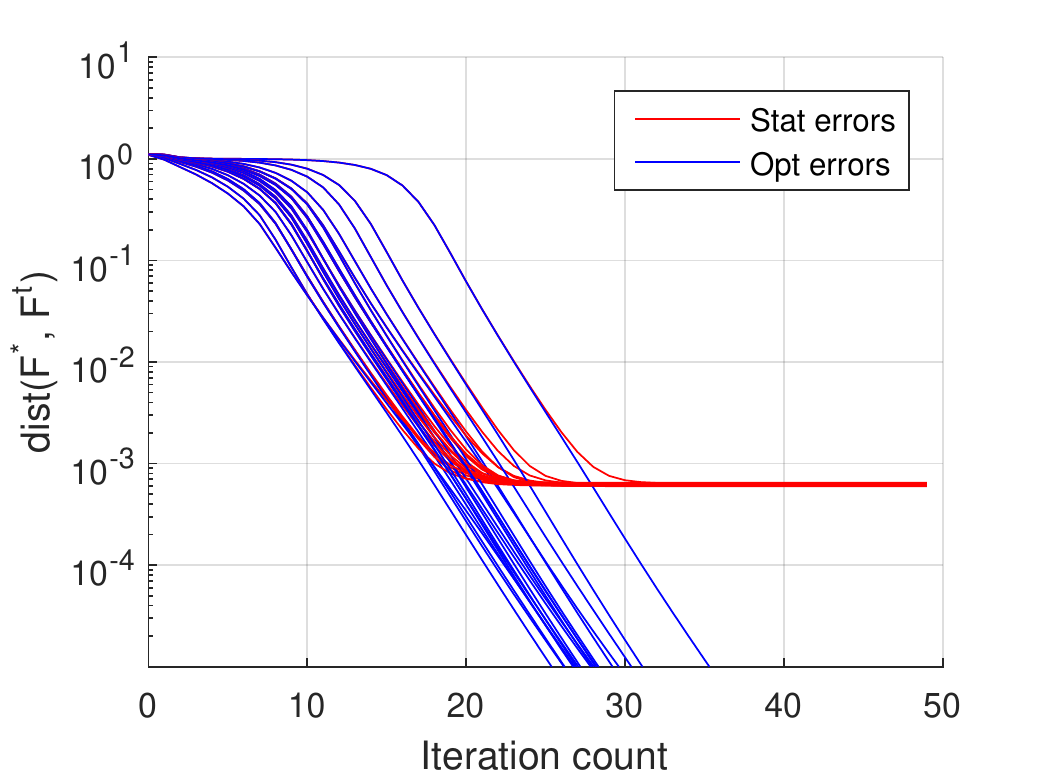} & 
\includegraphics[scale=0.5, clip, trim = 0 0 0 0]{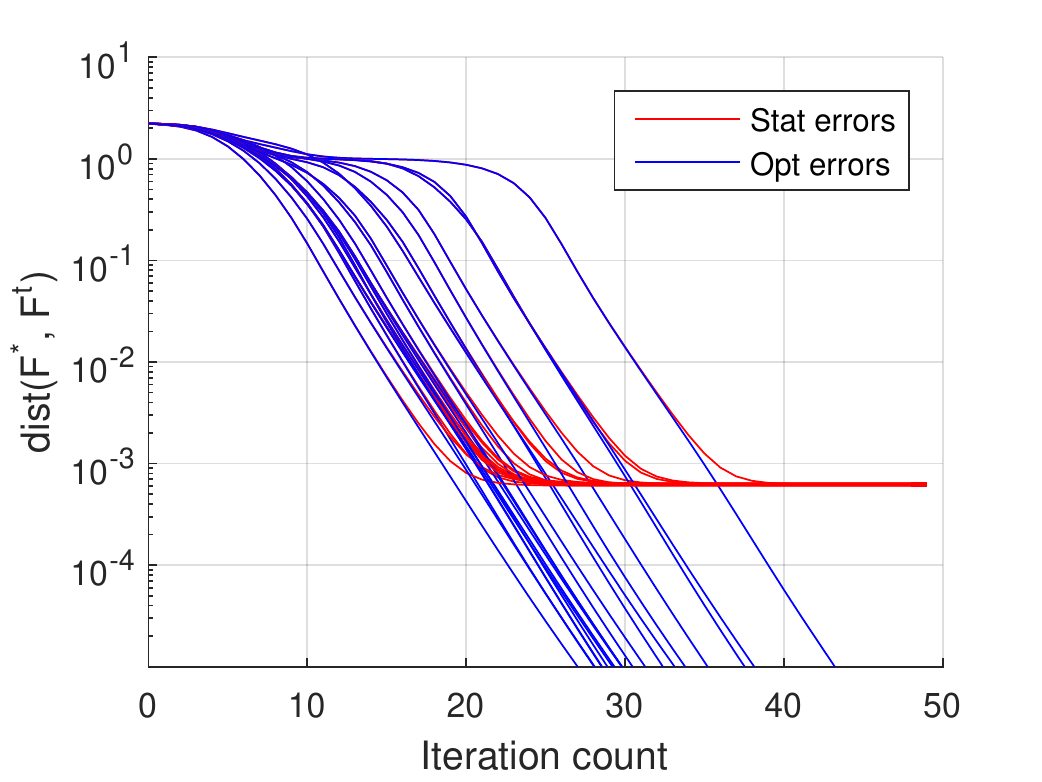} &
\includegraphics[scale=0.5, clip, trim = 0 0 0 0]{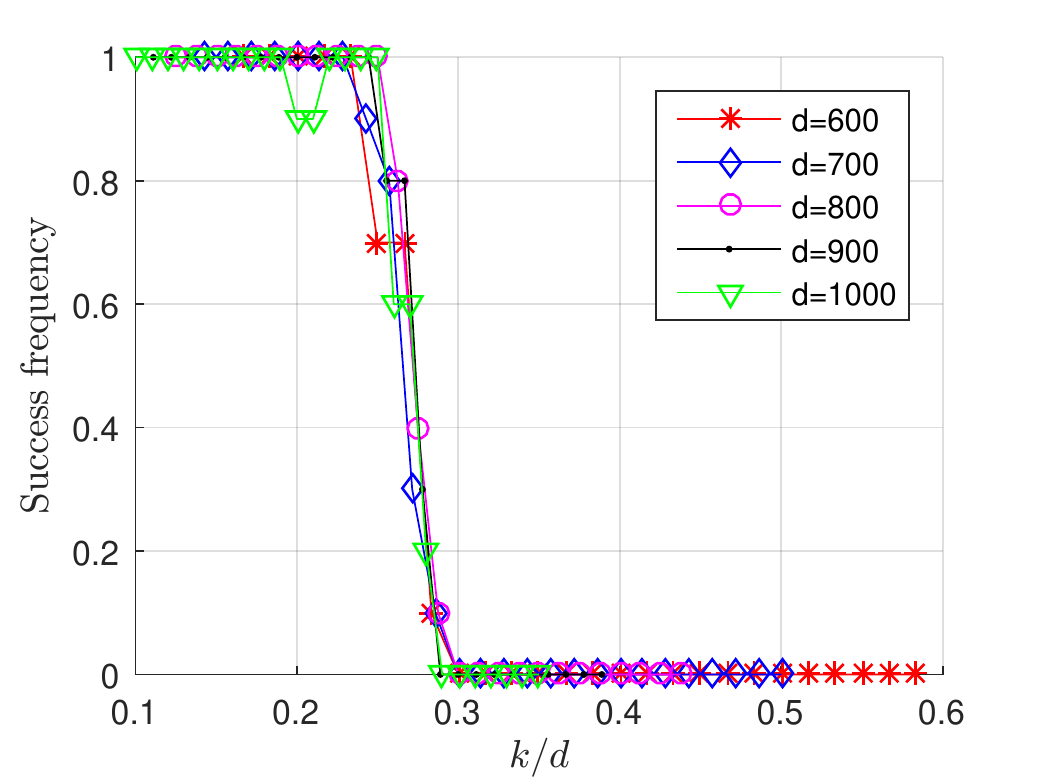} \\
(a) & (b) & (c)
\end{tabular} 
\caption{Matrix decomposition: plots of optimization error $
  \textup{d} (\thetait{t}, \thetait{T} ) $ and statistical error $
  \dist{\thetait{t}} $ versus the iteration number $ t$, using (a)
  SVD-based initialization and (b) random initialization. The
  simulation is performed using $ \usedim = 600 $, $ \rdim = 5 $, $
  \kdim =100 $ and $ \noisestd = 0.1 \cdot \frac{\rdim}{\usedim} $.
  Panel (c): plots of the probability of successful exact recovery of
  $ \thetastar $ versus $ \frac{\kdim}{\usedim} $, for different
  values of $ (\usedim, \kdim) $ using SVD-based initialization. We
  declare exact recovery if $ \dist{\widehat{\pltheta}} \le 2\times
  10^{-3} $, and each point represents frequency of exact recovery
  over $ 20 $ random instances.  The simulation is performed using $
  \rdim=6 $ and $ \noisestd=0 $.  }
	\label{FigMatDecomp}
\end{figure}


\paragraph{Initialization:}

We describe how to get a good initial matrix $\thetait{0}$ in the
noiseless setting $\noisevar= 0$. Suppose $\opnorm{\thetastar} = 1$.
Let $\widebar{\Yout}$ be obtained from hard-thresholding $\Yout$ at
the level $\frac{\inco\rdim}{\usedim}$; that is, for each element
$(i,j)$,
\begin{align*}
\widebar{\Yout}_{ij} & =
\begin{cases}
\Yout_{ij}, & \text{if \ensuremath{\abs{\Yout_{ij}} \le \frac{\inco\rdim}{\usedim}},} \\
\frac{\inco\rdim}{\usedim}\sgn(\Yout_{ij}), & \text{if \ensuremath{\abs{\Yout_{ij}} > \frac{\inco\rdim}{\usedim}}.}
\end{cases}
\end{align*}
We then set $\thetait{0}$ to be the $\usedim\times \rdim$ matrix with
columns being the top-$\rdim$ singular vector of $\bar{\Yout}$
projected onto the set $ \GenCons $.  In
Appendix~\ref{AppProofLemLSini}, we prove that under these conditions,
we have
\begin{align}
\label{EqnLSini}
\dist{\thetait{0}} \le \frac{4\inco\rdim\sqrt{\rdim}\kdim}{\usedim}.
\end{align}
Therefore, the requirement in Corollary~\ref{CorLS} is satisfied if
$\frac{\inco\rdim\sqrt{\rdim}\kdim}{\usedim} \le c_{1}$ for a
universal constant $c_{1}$ that is sufficiently small. The condition
$\frac{\inco\rdim\sqrt{\rdim}\kdim}{\usedim} \le c_{1}$ is sub-optimal
by a factor of $\sqrt{\rdim}$.

\paragraph{Computation:}

To compute the gradient $\nabla_{\PlTheta} \EmpLoss(\PlTheta)=
\PlTheta + \ProjConsS(\Yout - \PlTheta) - \Yout$, we need to project
each row of $\Yout - \PlTheta$ to the $\ell_{1}$ balls, which can be
done efficiently~\cite{duchi2008L1project}.  As discussed in
Section~\ref{SecMC}, the projection $\ProjGenCons$ can be computed by
row-wise clipping.\\

\paragraph{Simulations:}

We performed experiments under the same general set-up as the matrix
completion (see the discussion surrounding
equation~\eqref{EqnSineApprox}).  The matrix $ \thetastar $ is random
orthonormal; the sparse matrix $ \Sstar $ has $ \kdim \times \usedim $
non-zero entries whose locations are sampled uniformly without
replacement and whose values are independently and uniformly sampled
from the interval $[0, 10 \cdot \frac{\rdim}{\usedim}] $; the noise
matrix $ \noisevar $ has i.i.d.\ zero-mean entries with standard
deviation $ \noisestd $. The initial matrix $ \thetait{0} $ is
obtained using the SVD-based procedure described in
Section~\ref{SecMatrixDecomposition}.  The step size is fixed at $
\stepit{t} \equiv 1$.  Panels (a) and (b) of Figure~\ref{FigMatDecomp}
confirms the predicted geometric convergence.

In terms of the estimation error, in the noiseless setting with
$\noisestd = 0$, the output $ \widebar{\pltheta} $ equals $ \thetastar
$ with high probability provided that $ \frac{\rdim \kdim}{\usedim}
\precsim 1 $.  Therefore, with $ \rdim $ fixed, exact recovery of $
\thetastar $ can be achieved with probability close to one as soon as
$ \frac{\kdim}{\usedim} $ is below a constant threshold.  Panel (c) of
Figure~\ref{FigMatDecomp} confirms this prediction.


\section{Proofs of general theorems}
\label{SecProofThmGeneral}

This section is devoted to the proofs of our general theorems on
projected gradient convergence, namely Theorem~\ref{ThmGeneralLip} in
the Lipschitz case, and Theorem~\ref{ThmGeneralSmooth} under stronger
smoothness conditions.  Throughout these proofs, we make use of the
convenient shorthand $\GradTilit{t} \defn \nabla
\LossTil(\thetait{t})$ for the gradient of $\LossTil$ at step $t$.  We
also define the difference matrix $\diffit{t} \defn \thetait{t + 1} -
\thetait{t}$, as well as the parameters $\starnorm \defn
\opnorm{\thetastar}$ and $ \sigma_\rdim \defn \sigma_\rdim(\thetastar)
$.


\subsection{Proof of Theorem~\ref{ThmGeneralLip}}
\label{SecProofThmGeneralLip}

Our proof proceeds via induction on the event
\begin{align}
\label{EqnDefnEvent}
\Event_t \defn \big \{ \dist{ \thetait{s} } \leq \underbrace{(1 -
  \inisig) \sigma_\rdim}_{\smallrad} \quad \mbox{for all $s \in \{
  0,1,\ldots, t\}$} \big \}.
\end{align}
For the base case $t = 0$, note that $\Event_0$ holds by the
assumptions of the theorem.  Assuming that $\Event_t$ holds, it
suffices to show $\dist{\thetait{t + 1}} \le \smallrad$, which then
implies that $\Event_{t+1}$ holds.

We require the following auxiliary result:
\begin{lem}
\label{LemThetaStar}
For any matrix $\pltheta \in \real^{\usedim\times \rdim}$ such that
$\dist{\pltheta} < \sigma_{\rdim}(\thetastar)$, the optimization problem
$\min_{A \in \eclass} \fronorm{A - \pltheta}$ has a unique optimum
$\PISTAR{\pltheta}$ such that (i) the matrix $\pltheta^{\top}
\tdthetastar \in \real^{\rdim\times \rdim}$ is positive definite; and
(ii) the matrix $(\pltheta - \tdthetastar)^{\top} \tdthetastar$ is
symmetric.
\end{lem}
\noindent See Section~\ref{SecProofLemThetaStar} for the proof of
this claim. 
In view of Lemma~\ref{LemThetaStar}, the matrix $\PISTAR{\thetait{s}}
\defn \arg \min_{A \in \eclass} \fronorm{A- \thetait{s}}$ is
uniquely defined for each time step $s \in \{0, 1, \ldots,t \}$.\\

The projected gradient descent update can be decomposed into the two
steps
\begin{align}
\preprojit{s + 1} = \thetait{s}- \stepit{s} \nabla
\LossTil(\thetait{s}) \quad \mbox{and} \quad \thetait{s + 1} =
\ProjGenCons \big( \preprojit{s + 1} \big). \label{EqnGD_decompose}
\end{align}
For each $ s \in \{0, 1, \ldots, t\} $,  the local descent condition~\eqref{EqnLtilCurveCond} implies that
\begin{align*}
\trinprod{\nabla \LossTil(\thetait{s})}{ \thetait{s} - \thetastarit{s} } 
\geq \curv \fronorm{\thetait{s} - \thetastarit{s}}^{2} - \curv \epsnum^{2}.
\end{align*}
On the other hand, from the decomposition~\eqref{EqnGD_decompose}, we have $\nabla
\LossTil(\thetait{s}) = \frac{\thetait{s} - \preprojit{s+1}}{\stepit{s}}$, and hence the above inequality implies that
\begin{align*}
\curv \fronorm{\thetait{s} - \thetastarit{s}}^{2} - \curv \epsnum^{2} 
& \leq \frac{1}{\stepit{s}}\trinprod{\thetait{s} - \preprojit{s + 1}}{\thetait{s} - \thetastarit{s} } \\ 
& = \frac{1}{2\stepit{s}} \big( \fronorm{\thetait{s} - \thetastarit{s} }^{2} +
\fronorm{\thetait{s} - \preprojit{s + 1}}^{2} - \fronorm{\preprojit{s + 1} - \thetastarit{s} }^{2} \big) \\
& = \frac{1}{2\stepit{s}} \big( \fronorm{\thetait{s} - \thetastarit{s} }^{2} - \fronorm{\preprojit{s + 1} - \thetastarit{s} }^{2} \big) + \stepit{s} \fronorm{\nabla \LossTil(\thetait{s})}^2.
\end{align*}
Due to the $\ThetaStar$-faithfulness and convexity assumption on $ \GenCons $, we are guaranteed that
$\thetastarit{s} \in \GenCons$, and hence
\begin{align*}
\fronorm{\thetait{s + 1} - \thetastarit{s}} \leq \fronorm{\preprojit{s + 1} - \thetastarit{s}}
\end{align*}
since Euclidean projection onto a convex set is non-expansive~\cite{Bertsekas_nonlin}.	Moreover, by the Lipschitz condition~(\ref{EqnLtilLipCon}) on $\LossTil$, we have $\fronorm{\nabla \LossTil(\thetait{s})}^2 \leq \LipCon^2 \starnorm^2$.  Combining
the pieces, we find that
\begin{align}
\label{EqnRio}
\curv \fronorm{\thetait{s} - \thetastarit{s} }^{2} - \curv\epsnum^{2} 
& \le \frac{1}{2\stepit{s}} \big( \fronorm{\thetait{s} - \thetastarit{s} }^{2} - \fronorm{\thetait{s + 1} - \thetastarit{s} }^{2} \big) + \stepit{s} \LipCon^{2} \starnorm^{2}.
\end{align}
Introducing the shorthand $\gamma \defn \frac{20 \cond^2\LipCon^{2}}{\curv^{2}}-1$, we then make the step size choice $\stepit{s} = \frac{1}{\curv(s + 1 + \gamma)}$.  Substituting into our bound~\eqref{EqnRio} and rearranging yields
\begin{align*}
\frac{\curv(s + 1 + \gamma)}{2} \fronorm{\thetait{s + 1} - \thetastarit{s}}^{2} 
& \le \frac{\curv(s-1 + \gamma)}{2} \fronorm{\thetait{s} - \thetastarit{s} }^{2} + \frac{1}{\curv(s + 1 + \gamma)} \LipCon^{2} \starnorm^{2} + \curv\epsnum^{2}.
\end{align*}
Multiplying both sides by $s + \gamma$ and using the fact that $\fronorm{\thetait{s + 1} - \thetastarit{s+1} } \le \fronorm{\thetait{s + 1} - \thetastarit{s}} $  yields
\begin{align*}
\frac{\curv(s + \gamma)(s + 1 + \gamma)}{2} \fronorm{\thetait{s + 1} -  \thetastarit{s+1}}^{2}
\le & \frac{\curv(s + \gamma)(s-1 + \gamma)}{2} \fronorm{\thetait{s} - \thetastarit{s}}^{2} +
\frac{(s + \gamma)}{\curv(s + 1 + \gamma)}\LipCon^{2}\starnorm^{2} +
\curv(s + \gamma)\epsnum^{2}. \\
\leq & \frac{\curv(s + \gamma)(s-1 + \gamma)}{2} \fronorm{\thetait{s} - \thetastarit{s}}^{2} +
\frac{1}{\curv}\LipCon^{2}\starnorm^{2} + \curv(s + \gamma)\epsnum^{2}.
\end{align*}
Summing the above inequality over $s= 0,\ldots t$ yields
\begin{align}
\label{EqnRioCentro}
\frac{\curv(t + \gamma)(t + 1 + \gamma)}{2} \fronorm{\thetait{t + 1} - \thetastarit{t+1}}^{2} 
& \le \frac{\curv\gamma^{2}}{2} \fronorm{\thetait{0} - \thetastarit{0}}^{2} + \frac{t + 1}{\curv}\LipCon^{2}\starnorm^{2} + \curv(t + 1)(t + \gamma)\epsnum^{2}.
\end{align}

Now observe that the assumptions $ \inisig \le \frac{1}{2} $ and $\curv\le \LipCon$  imply that $ \gamma \ge \frac{4\cond^2 \LipCon^2 }{(1-\inisig)^2 \curv^2  } \ge 1 $.
These inequalities, combined with the facts that $\fronorm{\thetait{0} - \thetastarit{0}} \le (1 - \inisig)\sigma_\rdim$ by assumption and $ \starnorm/\cond = \sigma_\rdim $, when applied
to the bound~\eqref{EqnRioCentro}, yield
\begin{align}
\fronorm{\thetait{t + 1} - \thetastarit{t + 1}}^{2} 
& \leq
\frac{\gamma^{2}}{(t + \gamma)(t + 1 + \gamma)}(1 -
\inisig)^{2}\sigma_\rdim^{2} + \frac{(t + 1)\gamma/2}{(t + \gamma)(t + 1
  + \gamma)}(1 - \inisig)^{2}\frac{\starnorm^{2}}{\cond^2} + \frac{2(t + 1)}{t + 1 +
  \gamma}\epsnum^{2} \nonumber \\
\label{EqnLipIntermediate}
& = \frac{\gamma^{2} + (t + 1) \gamma/2}{(t + \gamma)(t + 1 +
  \gamma)}(1 - \inisig)^{2}\sigma_\rdim^{2} + \frac{2(t + 1)}{t + 1 +
  \gamma}\epsnum^{2}.
\end{align}
This bound, together with the assumed bound $\epsnum\le \frac{(1 -
  \inisig)\sigma_\rdim}{2}$ yields
\begin{align*}
\fronorm{\thetait{t + 1} - \thetastarit{t + 1}}^{2} & \le
\frac{\gamma^{2} + (t + 1)\gamma/2 + (t + 1)(t + \gamma)/2}{(t +
  \gamma)(t + 1 + \gamma)}(1 - \inisig)^{2}\sigma_\rdim^{2} \le (1 -
\inisig)^{2}\sigma_\rdim^{2} \; = \; \smallrad^2
\end{align*}
whence $\dist{\thetait{t + 1}} \le \smallrad$, thereby proving the
induction hypothesis for $t + 1$.

Moreover, since $\gamma \geq 1$, the
inequality~\eqref{EqnLipIntermediate} implies that
\begin{align*}
\fronorm{\thetait{t + 1} - \thetastar_{t + 1}}^{2} 
\le \frac{\gamma}{t + \gamma}(1 - \inisig)^{2}\sigma_\rdim^{2} + 2\epsnum^{2} 
\le \frac{20\LipCon^{2} \starnorm^{2}}{(t + 1)\curv^{2}} + 4\epsnum^{2},
\end{align*}
thereby establishing the bound~\eqref{EqnLipGuarantee} stated in the theorem.


\subsubsection{Proof of Lemma~\ref{LemThetaStar}}
\label{SecProofLemThetaStar}

We use the shorthand  $\sigma_k = \sigma_k(\thetastar)$ for $ k=1,\ldots,\rdim $.  Since $\dist{\pltheta} =\min \limits_{A\in \eclass}
\fronorm{\pltheta - A} < \sigma_\rdim$, there must exist a matrix $\thetastar_{0}
\in \eclass$ such that
\begin{align}
\label{EqnCafeZinho}
\opnorm{\pltheta - \thetastar_{0}} \le \fronorm{\pltheta -
  \thetastar_{0}} < \sigma_\rdim.
\end{align}
It follows that the
matrix $\pltheta$ must have full column rank with all its singular
values contained in the interval $\big[\sigma_\rdim - \dist{\pltheta}, \sigma_1+
\dist{\pltheta} \big]$. Let $\I_{\rdim}$ denote the $\rdim$-dimensional identity matrix, and the rank-$\rdim $ SVD of $ \thetastar_{0} $ be $ \thetastar_{0} = VSR^\top $, where $ V\in \real^{\usedim\times \rdim} $, $ S=\text{diag}( \sigma_1, \ldots, \sigma_\rdim ) $ and $ R\in\real^{\rdim\times\rdim} $ is orthonormal.  
Any unit vector in $ \real^\rdim $ has the form $  Rw $ for some unit vector $ w\in \real^\rdim $, and 
$ \pltheta^\top\thetastar_{0} Rw
= (\thetastar_{0} + \pltheta - \thetastar_{0})^\top VSw
= (RS  + (\pltheta - \thetastar_{0})^\top V) Sw  $.
We therefore have the bound
\begin{align*}
\twonorm{\pltheta^\top\thetastar_{0} Rw}
\ge \sigmin{RS  + (\pltheta - \thetastar_{0})^\top V} \twonorm{Sw} 
&\ge (\sigmin{RS} - \dist{\pltheta} )  \sigma_\rdim
= (\sigma_\rdim - \dist{\pltheta} )  \sigma_\rdim.
\end{align*}
It follows that the $\rdim$-dimensional matrix $U \defn
\pltheta^{\top}\thetastar_{0}$  satisfies $ \sigma_\rdim(U) \ge (\sigma_\rdim - \dist{\pltheta} )  \sigma_\rdim >0 $ and is thus invertible. A similar argument shows that $ \sigma_1(U) \le (\sigma_1+ \dist{\pltheta} )  \sigma_1  $.

Defining the matrix $\tdthetastar \defn
\thetastar_{0}U^{\top}(UU^{\top})^{-1/2}$, it is easy to verify
$\tdthetastar \in \eclass$. Observe that
\begin{align*}
\pltheta^{\top}\tdthetastar & =
\pltheta^{\top}\thetastar_{0}U^{\top}(UU^{\top})^{-1/2} =
UU^{\top}(UU^{\top})^{-1/2} = (UU^{\top})^{1/2},
\end{align*}
which is symmetric and positive definite since $U$ has strictly
positive singular values.  It is then clear that the matrix $(\pltheta
- \tdthetastar)^{\top}\tdthetastar$ is symmetric.

Any matrix $A \in \eclass$ can be written as $A=
\tdthetastar\Xi$ for some orthonormal matrix $\Xi \in \real^{\rdim\times\rdim}$, whence
\begin{align}
\label{EqnStar}
\trace(\pltheta^{\top}A)= \trace \big( (UU^{\top})^{1/2} \Xi \big) \leq \opnorm{\Xi}
\sum_{i= 1}^{\rdim} \sigma_{i}(U) = \sum_{i= 1}^{\rdim}\sigma_{i}(U).
\end{align}
Noting that
\begin{align*}\fronorm{\pltheta - \tdthetastar}^{2} \; = \;
\fronorm{\pltheta}^{2} + \fronorm{\tdthetastar}^{2} - 2\trace
\big( (UU^{\top})^{1/2} \big) & =
\fronorm{\pltheta}^{2} + \fronorm{A}^{2} - 2\sum_{i=
  1}^{\rdim}\sigma_{i}(U),
\end{align*}
we thus have the bound
\begin{align*}
\fronorm{\pltheta - \tdthetastar}^{2} & \leq \fronorm{\pltheta}^{2} +
\fronorm{A}^{2} - 2\trace \big( \pltheta^{\top} A \big) \; = \;
\fronorm{\pltheta - A}^{2}.
\end{align*}
Since $A$ was arbitrary in $\eclass$. we conclude that $\tdthetastar$
is a constrained minimizer of $\fronorm{\pltheta - \thetastar}$ over
$\eclass$.

Finally, we claim that the inequality in~\eqref{EqnStar} is strict if
$\Xi\neq\I_{\rdim}$, so that $\tdthetastar$ is in fact the unique
minimizer, as claimed.	To establish this strictness, suppose that the
SVD of $(UU^{\top})^{1/2}$ is given by $(UU^{\top})^{1/2} = R' \Sigma {R'}^{\top}$, where $\Sigma =
\text{diag} \big( \sigma_{1}(U),\ldots,\sigma_{\rdim}(U) \big)$ and
$R' $ is an $\rdim\times \rdim$ orthonormal matrix. Then $\trace\big( (UU^{\top})^{1/2}\Xi \big)=
\trace(\Sigma {R'}^{\top}\Xi R')= \trace(\Sigma\Xi')$, where $\Xi' \defn
{R'}^{\top}\Xi R'$ is also orthonormal. If $\Xi\neq\I_{\rdim}$, then
$\Xi'\neq\I_{\rdim}$ and therefore
\begin{align*}
\trace\big( (UU^{\top})^{1/2}\Xi \big) 
\; = \; \sum_{i= 1}^{r}\sigma_{i}(U)\Xi'_{ii} 
\; < \;
\sum_{i= 1}^{r}\sigma_{i}(U)\twonorm{\Xi_{i\cdot}'} \; = \; \sum_{i=
  1}^{r}\sigma_{i}(U),
\end{align*}
where the inequality follows from the facts that $\sigma_{i}(U) > 0$
for each index $i \in [\rdim]$, and since $\Xi'\neq\I_{\rdim}$, we
must have $\abs{\Xi'_{ii}} < 1 = \twonorm{\Xi_{i\cdot}'}$ for some
index $i \in [\rdim]$.


\subsection{Proof of Theorem~\ref{ThmGeneralSmooth}}
\label{SecProofThmGeneralSmooth}

We will in fact prove a slightly stronger form of the theorem, where the $ \LipCon $-Lipschitz condition is replaced by the following relaxed condition:
\begin{equation}
\abs{\trinprod{\nabla_{\PlTheta} \LossTil( \OUTER{\pltheta})}{\pltheta - \pltheta'}} \le \LipCon \big( \opnorm{\thetastar}^{2} + \opnorm{\thetastar} \fronorm{\pltheta - \pltheta'} \big),  \label{EqnLtilWeakLipCon}
\end{equation}
valid for all $ \pltheta \in \GenCons \cap
\Ballstar{(1 - \inisig^{2}) \sigma_{\rdim}(\thetastar)} $ and $ \pltheta' \in
\GenCons $.
We use the step size choice $\stepit{t} \defn
c_{\inisig} \frac{\curv}{\smoo} = \frac{\inisig^{4}(1 -
  \inisig)^{8}}{C} \frac{\curv}{\cond^6\smoo^{2}}$, where $C = 272$.

As in the proof of Theorem~\ref{ThmGeneralLip}, we will show that
$\dist{\thetait{t}} \le  (1 - \inisig)\starnorm$ for all iterates $t =
0, 1, 2, \ldots$.  We do so via induction on the event $\Event_t$
previously defined in equation~\eqref{EqnDefnEvent}.  As before, the
base event $\Event_0$ holds by the theorem's conditions.  The
induction step is based on assuming that $\Event_t$ holds, and then
showing that $\Event_{t+1}$ also holds.  As before, it suffices to
establish the bound $\dist{\thetait{t + 1}} \le (1 -
\inisig)\starnorm$. The proof is divided  into several steps below. \\

\paragraph{Showing that $\dist{\thetait{t + 1}}  \le \smallrad = (1-\inisig^2) \starnorm$:}

\noindent In the first step, we establish a slightly weaker bound on $ \dist{\thetait{t + 1}}  $. 

In view of Lemma~\ref{LemThetaStar}, the matrix $\PISTAR{\thetait{s}}
\defn \arg \min_{A \in \eclass} \fronorm{A- \thetait{s}}$ is
uniquely defined for each time step $s \in \{0, 1, \ldots,t \}$.
Recall that the iterate $\thetait{t + 1}$ is an optimal solution to
the convex optimization problem~\eqref{EqnGenProjGradOpt}.
Consequently, the first-order conditions for optimality imply that
\begin{align}
\label{EqnOptCond}
\trinprod{\GradTilit{t} + \frac{1}{\stepit{t}}\diffit{t}}{ -
  \thetait{t + 1} + \pltheta} \ge 0, \quad \forall\pltheta \in
\GenCons.
\end{align}
Applying this condition with $\pltheta = \thetait{t}$ and using the relaxed Lipschitz condition~\eqref{EqnWeakLipCon} with the assumption $ \LipCon = \smoo $, we obtain
\begin{align*}
\fronorm{\diffit{t}}^{2} & \le \stepit{t}\trinprod{\GradTilit{t}}{ - \diffit{t}} 
\le  \stepit{t}\smoo(\starnorm^{2} + \starnorm\fronorm{\diffit{t}}) 
= \frac{\inisig^{8}(1 - \inisig)^{2}\curv}{C\cond^6\smoo}(\starnorm^{2}  + \starnorm\fronorm{\diffit{t}}).
\end{align*}
With the constant $C = 272$ and the fact that $ \max\{\inisig, 1/\cond, \curv/\smoo\} \le 1 $, this inequality implies that
$\fronorm{\diffit{t}} \le \inisig(1 - \inisig)\starnorm$ and hence
that
\begin{align}
\label{EqnWeakBallBound}
\dist{\thetait{t + 1}} \le \fronorm{\thetait{t + 1} - \thetastarit{t}}
& \le \fronorm{\thetait{t} - \thetastarit{t}} + \fronorm{\diffit{t}} 
\le  (1 - \inisig^{2})\starnorm.
\end{align}
Note that we have not yet completed the induction step, but the bound~(\ref{EqnWeakBallBound}) is useful below.

\paragraph{Establishing a recursive bound:}

With $ \thetait{t+1} $ satisfying the bound~(\ref{EqnWeakBallBound}), Lemma~\ref{LemThetaStar} guarantees that the matrix $\thetastarit{t +
  1} \defn \arg\min_{A \in \eclass} \fronorm{\thetastar - \thetait{t +
    1}}$ is well-defined.  Our analysis involves the matrix
differences $\d_{s'}^{s} \defn \thetait{s} - \thetastarit{s'}$,
defined by pairs $s, s' \in \{ 0,1,\ldots,t + 1 \}$.  The following
inequality is central to our analysis:
\begin{align}
\label{EqnAccor}
\frac{1}{\stepit{t}}\trinprod{\diffit{t}}{ - \d_{t}^{t + 1}} 
& \ge \frac{\curv}{2} \fronorm{\d_{t + 1}^{t + 1}}^{2} - \frac{136\cond^6 \smoo^{2} \fronorm{\diffit{t}}^{2}}{\curv\inisig^{8}} - 3\curv\epsnum^{2}.
\end{align}
We return to prove it shortly; taking it as given for the moment, it
follows that
\begin{align*}
\fronorm{\d_{t + 1}^{t + 1}}^{2} \le \fronorm{\d_{t}^{t + 1}}^{2} & = \fronorm{\thetait{t + 1} - \thetastar_{t}}^{2} \\
& = \fronorm{ \thetait{t} - \thetait{t + 1} + \thetait{t + 1} - \thetastar_{t}}^{2} - \fronorm{\diffit{t}}^{2} + 2\trinprod{\diffit{t}}{\thetait{t + 1} - \thetastar_{t}} \\
& \le \fronorm{\d_{t}^{t}}^{2} - \fronorm{\diffit{t}}^{2} + \stepit{t}\cdot \Bigl( - \curv \fronorm{\d_{t + 1}^{t + 1}}^{2} + \frac{272\cond^6\smoo^{2}}{\curv\inisig^{8}}
\fronorm{\diffit{t}}^{2} + 6\curv\epsnum^{2} \Bigr) \\
& = \fronorm{\d_{t}^{t}}^{2} - \fronorm{\diffit{t}}^{2} - \frac{\curv^{2}\inisig^{8}(1-\inisig)^2}{C\cond^6 \smoo^{2}} \fronorm{\d_{t + 1}^{t + 1}}^{2} + \frac{272(1-\inisig)^2}{C} \fronorm{\diffit{t}}^{2} + \frac{6\curv^{2}\inisig^{8}(1-\inisig)^2}{C\cond^6 \smoo^{2}}\epsnum^{2} \\
& \le \fronorm{\d_{t}^{t}}^{2} - \frac{\curv^{2}\inisig^{8}}{C\cond^6\smoo^{2}} \fronorm{\d_{t + 1}^{t + 1}}^{2} + \frac{6\curv^{2}\inisig^{8}}{C\cond^6\smoo^{2}}\epsnum^{2}.
\end{align*}
where we used the step size choice $\stepit{t} = \frac{\curv\inisig^{8} (1-\inisig)^2}{C\cond^6\smoo^{2}}$ and the assumption $C \geq 272$ in the last two lines. Rearranging this inequality yields the recursive bound
\begin{align}
\label{EqnDecay}
\fronorm{\d_{t + 1}^{t + 1}}^{2} 
& \le \Big( 1 + \frac{\curv^{2}\inisig^{8}}{C\cond^6\smoo^{2}} \Big)^{-1} \Big(
\fronorm{\d_{t}^{t}}^{2} + \frac{6\curv^{2}\inisig^{8}}{C\cond^6\smoo^{2}}\epsnum^{2} \Big)
\le \Big( 1 - \frac{\curv^{2}\inisig^{8}}{2C\cond^6\smoo^{2}} \Big) \Big( \fronorm{\d_{t}^{t}}^{2} + \frac{6\curv^{2}\inisig^{8}}{C\cond^6\smoo^{2}}\epsnum^{2} \Big).
\end{align}

\paragraph{Completing the induction and proof:}

Since $ \fronorm{\d_t^t} \le (1-\inisig) \starnorm $ by induction hypothesis and $ \epsnum \le \frac{1-\inisig}{4} \starnorm $ by assumption, the inequality~(\ref{EqnDecay}) above implies
\[
\fronorm{\d_{t + 1}^{t + 1}}^{2} \le  \Big( 1 - \frac{\curv^{2}\inisig^{8}}{2C\cond^6\smoo^{2}} \Big) \Big( 1 + \frac{\curv^{2}\inisig^{8}}{2C\cond^6\smoo^{2}} \Big)(1 - \inisig)^{2}\psi^{2} \le (1 - \inisig)^{2}\psi^{2},
\]
which completes the induction step.  Moreover, by applying the inequality~(\ref{EqnDecay}) recursively, we find that
\begin{align*}
\fronorm{\d_{t}^{t}}^{2} 
\le \Big( 1 - \frac{\curv^{2}\inisig^{8}}{2C\cond^6\smoo^{2}} \Big)^{t}
\fronorm{\d_{0}^{0}}^{2} + \frac{2C\cond^6\smoo^{2}}{\curv^{2}\inisig^{8}} \cdot \frac{6\curv^{2}\inisig^{8}}{C\cond^6\smoo^{2}}\epsnum^{2}
\le \Big( 1 - \frac{\curv^{2}\inisig^{8}}{2C\cond^6\smoo^{2}} \Big)^{t}
\fronorm{\d_{0}^{0}}^{2} + (4\epsnum)^{2},
\end{align*}
thereby completing the proof of the theorem.

\paragraph{Proof of inequality~\eqref{EqnAccor}:}
It remains to prove the intermediate claim~\eqref{EqnAccor}.  With $ \dist{\thetait{t+1}} \le \smallrad = (1-\inisig^2) \psi$ as established in~(\ref{EqnWeakBallBound}),  the local descent condition~\eqref{EqnLtilCurveCond} yields
\begin{align*}
\trinprod{\GradTilit{t + 1}}{\d_{t}^{t + 1}} 
 & \ge \curv \fronorm{\d_{t + 1}^{t + 1}}^{2} - \frac{\smoo^2}{\curv} \fronorm{\thetastarit{t + 1} - \thetastarit{t} }^2 - \curv \epsnum^2.
\end{align*}
In order to proceed, we need a second technical lemma. Recall that $ \cond = \frac{\sigma_1(\thetastar)}{\sigma_\rdim(\thetastar) }$ is the condition number of $ \thetastar $.
\begin{lem}
\label{LemThetaStarIt}
Under the conditions of Theorem~\ref{ThmGeneralSmooth}, we have
\begin{align*}
\fronorm{\thetastarit{t + 1} - \thetastarit{t}} \le
\frac{10 \cond^3 \fronorm{\diffit{t}}}{\inisig^{4}}.
\end{align*}
\end{lem}
\noindent See Section~\ref{SecProofLemThetaStarIt} for the proof of
this claim. \\

\noindent Applying Lemma~\ref{LemThetaStarIt}, we find that
\begin{align}
\trinprod{\GradTilit{t + 1}}{\d_{t}^{t + 1}}  
& \geq \curv \fronorm{\d_{t + 1}^{t +1}}^{2} - \frac{100\cond^6 \smoo^{2}}{\curv\inisig^{8}} \fronorm{\diffit{t}}^{2} - \curv\epsnum^{2}.
\label{EqnCurveConseq}
\end{align}

On the other hand, the smoothness condition~\eqref{EqnLtilSmoothCond} yields that
\begin{align*}
\trinprod{\GradTilit{t} - \GradTilit{t + 1}}{\d_{t}^{t + 1}} 
& \geq - \smoo\fronorm{\diffit{t}} \fronorm{\d_{t}^{t + 1}} - \curv\epsnum\fronorm{\d_{t}^{t + 1}}.
\end{align*}
Together with Lemma~\ref{LemThetaStarIt}, we obtain
\begin{align}
& \trinprod{\GradTilit{t} - \GradTilit{t + 1}}{\d_{t}^{t + 1}} \\
\geq &  - \Big( \smoo\fronorm{\diffit{t}} + \curv\epsnum \Big)
\Big( \fronorm{\d_{t + 1}^{t + 1}} +
\frac{10\cond^3 \fronorm{\diffit{t}}}{\inisig^{4}} \Big).\nonumber \\
= &  - \smoo \fronorm{\diffit{t}} \fronorm{\d_{t + 1}^{t + 1}} - \curv\epsnum\fronorm{\d_{t + 1}^{t + 1}} - \frac{10\cond^3 \smoo}{\inisig^{4}} \fronorm{\diffit{t}}^{2} - \frac{10\cond^3 \curv}{\inisig^{4}}\epsnum\fronorm{\diffit{t}}  \nonumber \\
\overset{(i)}{\ge} &  - \Big( \frac{\curv}{4} \fronorm{\d_{t + 1}^{t + 1}}^{2} + \frac{\smoo^{2}}{\curv} \fronorm{\diffit{t}}^{2} \Big)
- \Big( \frac{\curv}{4} \fronorm{\d_{t + 1}^{t + 1}}^{2} + \curv\epsnum^{2} \Big) 
- \frac{10\cond^3 \smoo}{\inisig^{4}} \fronorm{\diffit{t}}^{2} 
- \Big( \frac{25\cond^6 \curv}{\inisig^{8}} \fronorm{\diffit{t}}^{2} + \curv\epsnum^{2} \Big)\nonumber \\
\label{EqnSmoothConseq}
\geq &  - \frac{\curv}{2} \fronorm{\d_{t + 1}^{t + 1}}^{2} - 2\curv\epsnum^{2} - \frac{36\cond^6 \smoo^{2}}{\curv\inisig^{8}}
\fronorm{\diffit{t}}^{2},
\end{align}
where the step~(i) follows from the AM-GM inequality.

Finally, the $\ThetaStar$-faithfulness of $\GenCons$ ensures that
$\thetastarit{t} \in \GenCons$, so that we may apply the
bound~\eqref{EqnOptCond} with $\pltheta = \thetastarit{t}$, thereby
obtaining
\begin{equation}
\label{EqnOptCondStar}
\trinprod{\GradTilit{t} + \frac{1}{\stepit{t}}\diffit{t}}{ -
  \d_{t}^{t + 1}} \ge 0.
\end{equation}
Adding together inequalities~\eqref{EqnCurveConseq},
~\eqref{EqnSmoothConseq} and~\eqref{EqnOptCondStar} yields the
claim~\eqref{EqnAccor}.


\subsubsection{Proof of Lemma~\ref{LemThetaStarIt}}
\label{SecProofLemThetaStarIt}

By dividing through by $ \starnorm $ we may assume that $ \starnorm=1 $, so $ \sigma_1(\thetastar) = \cond $, where we recall that $ \cond $ is the condition number of $ \thetastar $. Define the matrices $U_{s} \defn (\thetait{s})^{\top}\thetastarit{t}$
for $s \in \{t,t + 1\}$, and recall that we have shown
\begin{align*}
\max \Big \{ \fronorm{\thetait{t} -
  \thetastarit{t}},\fronorm{\thetait{t + 1} - \thetastarit{t}} \Big \}
& \leq 1 - \inisig^2.
\end{align*}
The same argument as in the proof of Lemma~\ref{LemThetaStar} from the
previous section show that the singular values of $U_{s}$ are in
the interval $[\inisig^2,\cond + (1 - \inisig^2)]$, and we have the expression 
$\thetastarit{s} \defn
\thetastarit{t}U_{s}^{\top}(U_{s}U_{s}^{\top})^{-1/2}$ for $s \in \{
t,t + 1 \}$. Since $U_{t} = U_{t}^{\top} =
(U_{t}U_{t}^{\top})^{1/2}$, we have
\begin{align*}
\fronorm{U_{t + 1}-U_{t}} = \fronorm{(\thetait{t + 1} -
  \thetait{t})^{\top}\thetastarit{t}} \le \sigma_1 \fronorm{\diffit{t}}.
\end{align*}
By applying a known perturbation bound for matrix square roots
(\cite[Lemma 15]{gao2014cca}), we find that
\begin{align*}
\fronorm{(U_{t}U_{t}^{\top})^{1/2}-(U_{t + 1}U_{t + 1}^{\top})^{1/2}}
& \le \frac{ \fronorm{U_{t}U_{t}^{\top}-U_{t + 1}U_{t + 1}^{\top}}}
	{\sigma_{\min}((U_{t}U_{t}^{\top})^{1/2}) +
	  \sigma_{\min}((U_{t + 1}U_{t + 1}^{\top})^{1/2})}\\
 & \leq \frac{1}{2\inisig^2 } \fronorm{U_{t}U_{t}^{\top}-U_{t + 1}U_{t +
    1}^{\top}} \\
& = \frac{1}{2\inisig^2 } \fronorm{(U_{t}-U_{t + 1})U_{t}^{\top} + U_{t +
    1}(U_{t}-U_{t + 1})^{\top}} \\
& \le \frac{2 \cond^2}{\inisig^2} \fronorm{\diffit{t}}.
\end{align*}
Moreover, we have
\begin{align*}
\fronorm{U_{t}^{\top} \left[(U_{t + 1} U_{t + 1}^{\top})^{-1/2} -
    (U_{t} U_{t}^{\top})^{-1/2} \right ] } & = \fronorm{(U_{t + 1}U_{t
    + 1}^{\top})^{-1/2} \left[(U_t U_t^{\top})^{1/2}-(U_{t+1}
    U_{t+1}^{\top})^{1/2} \right ] } \\
& \leq \frac{1}{\inisig^2 } \fronorm{(U_{t}U_{t}^{\top})^{1/2} - (U_{t +
    1} U_{t+1}^{\top})^{1/2}} \\
& \leq \frac{2 \cond^2 \fronorm{\diffit{t}}}
	{\inisig^{4} }.
\end{align*}
Putting together the pieces, it follows that
\begin{align*}
\fronorm{\thetastar_{t + 1} - \thetastar_{t}} & =
\fronorm{\thetastarit{t}U_{t + 1}^{\top}(U_{t + 1}U_{t +
    1}^{\top})^{-1/2} - \thetastarit{t}
  U_{t}^{\top}(U_{t}U_{t}^{\top})^{-1/2}} \\
& = \fronorm{\thetastarit{t} ( U_{t + 1} - U_{t} )^{\top}(U_{t +
    1}U_{t + 1}^{\top})^{-1/2} + \thetastarit{t}U_{t}^{\top}
  \left[(U_{t + 1}U_{t+1}^{\top})^{-1/2} - (U_{t}U_{t}^{\top})^{-1/2}
    \right ] } \\
& \leq  \opnorm{\thetastarit{t}} \fronorm{U_{t + 1}-U_{t}}
\opnorm{(U_{t + 1}U_{t + 1}^{\top})^{-1/2}} + \opnorm{\thetastarit{t}}
\fronorm{U_{t}^{\top} \left[(U_{t + 1}U_{t +
      1}^{\top})^{-1/2}-(U_{t}U_{t}^{\top})^{-1/2}\right]} \\
& \leq \cond \cdot \cond \fronorm{\diffit{t}} \cdot \frac{1}{\inisig^2} +
\cond \cdot \frac{2 \cond^2 \fronorm{\diffit{t}}}{\inisig^{4}} \\
& \leq \frac{10 \cond^3 \fronorm{\diffit{t}}}{\inisig^{4}},
\end{align*}
as claimed.



\section{Proofs of corollaries}
\label{SecProofCor}

In this section, we prove the corollaries by applying our general
theory.


The general theorems in Section~\ref{SecMain} are stated in terms of
the loss function $ \LossTil $ of the factor variable $ \pltheta
$. Sometimes it is convenient to work with the original loss function
$ \EmpLoss $ of the $ \usedim\times\usedim $ variable $ \PlTheta
$. These two loss functions are related by $ \LossTil(\pltheta ) =
\EmpLoss(\OUTER{\pltheta}) $ and $ \nabla_{\pltheta}
\LossTil(\pltheta) = \big[ \nabla_{\PlTheta}
  \EmpLoss(\OUTER{\pltheta}) + (\nabla_{\PlTheta}
  \EmpLoss(\OUTER{\pltheta}) )^\top \big] \pltheta $, and the
convergence results can be restated in terms of $ \EmpLoss $. We do so
below for the result in Theorem~\ref{ThmGeneralSmooth}.\\

The following conditions for $ \EmpLoss $ are the counterparts of the corresponding conditions for $ \LossTil $.
\begin{dfn}[Local descent condition for $ \EmpLoss $] 
\label{DefCurveCon} For some curvature
parameter $\curv$, statistical tolerance $\epsnum$ and radius $\smallrad$,
we say that the cost function $\EmpLoss$ satisfies a \emph{local
descent condition} with parameters $(\curv,\epsnum,\smallrad)$
over $\GenCons$ if for each $\pltheta \in \GenCons \cap \Ballstar{\smallrad}$,
there exists $\tdthetastar \in \arg\min_{\thetastar \in \eclass}\fronorm{\thetastar - \pltheta}$
such that
\begin{align}
\trinprod{\nabla_{\PlTheta}\EmpLoss(\OUTER{\pltheta})}{\OUTER{\pltheta} - \OUTER{\tdthetastar} + \OUTER{(\pltheta - \tdthetastar)}} 
& \ge 2 \curv \fronorm{\pltheta - \tdthetastar}^{2} - \frac{1}{4} \curv \epsnum \fronorm{\pltheta - \tdthetastar}. 
\label{EqnCurveCond}
\end{align}
\end{dfn}
\begin{dfn}[Relaxed Local Lipschitz condition for $ \EmpLoss $]  
\label{DefLipCon} 
For some Lipschitz constant $\LipCon$ and radius $\smallrad$, we say that the $\EmpLoss$
satisfies a \emph{relaxed local Lipschitz condition} with parameter $(\LipCon,\smallrad)$
over $\GenCons$ if for each $\pltheta \in \GenCons \cap \Ball_{2}(\smallrad;\thetastar)$ and $ \pltheta' \in \GenCons $,
\begin{align}
\abs{ \trinprod{\nabla_{\PlTheta} \EmpLoss(\OUTER{\pltheta})}{(\pltheta-\pltheta') \otimes \pltheta} } \le \frac{1}{2} \LipCon \big( \opnorm{\thetastar}^2 + \opnorm{\thetastar} \fronorm{\pltheta - \pltheta'} \big).
\label{EqnWeakLipCon}
\end{align}
\end{dfn}
\noindent Of course, this relaxed Lipschitz condition for $ \EmpLoss $ is implied by a Lipschitz condition of the form
\begin{align}
\fronorm{\nabla_{\PlTheta}\EmpLoss(\OUTER{\pltheta})\pltheta} & \leq \frac{1}{2} \LipCon \opnorm{\thetastar}, \quad \forall \pltheta \in \GenCons \cap \Ball_{2}(\smallrad;\thetastar).
\label{EqnLipCon}
\end{align}
\begin{dfn}[Local smoothness condition for $ \EmpLoss$] 
\label{DefSmoothCon} 
For some curvature and smoothness parameters $\curv$ and $\smoo$, statistical tolerance
$\epsnum$ and radius $\smallrad$, we say that the loss function
$\EmpLoss$ satisfies a\emph{ local smoothness condition} with parameters
$(\curv,\smoo,\epsnum,\smallrad)$ over $\GenCons$ if for each $\pltheta,\pltheta',\pltheta'' \in \GenCons \cap \Ballstar{\smallrad}$
and any $\thetastar \in \eclass$, 
\begin{subequations}
\begin{align}
\label{EqnSmoothCond1}
\abs{\trinprod{\nabla_{\PlTheta} \EmpLoss( \OUTER{\pltheta} ) - \nabla_{\PlTheta} \EmpLoss( \OUTER{\pltheta'} )}{\pltheta'\otimes(\pltheta - \thetastar)}} 
& \le
\frac{1}{4} \big( \smoo\fronorm{\pltheta - \pltheta'} + \curv\epsnum \big) \fronorm{\pltheta  - \thetastar}, \\
\label{EqnSmoothCond2}
\abs{\trinprod{\nabla_{\PlTheta} \EmpLoss( \OUTER{\pltheta} )}{(\pltheta - \thetastar) \otimes (\pltheta' - \pltheta'')}} 
& \le
\frac{1}{4} \big( \smoo \fronorm{\pltheta' - \pltheta''} + \curv\epsnum \big) \fronorm{\pltheta - \thetastar}.
\end{align}
\end{subequations}
\end{dfn}
\vtiny

Now suppose that for some numbers $\curv$, $\smoo$, $\LipCon$,
$\epsnum$ and $\inisig$ with $0 < \curv \le \smoo = \LipCon$, $ 0<
\inisig <1 $ and $\epsnum \le \frac{1 - \inisig}{4}
\sigma_{\rdim}(\thetastar)$, the empirical loss $\EmpLoss$ satisfies
the local descent, relaxed Lipschitz and smoothness conditions in
Definitions~\ref{DefCurveCon}--\ref{DefSmoothCon} over $\GenCons$ with
parameters $\curv$, $\smoo$, $\LipCon$, $\epsnum$ and $\smallrad= (1 -
\inisig^{2}) \sigma_{\rdim}(\thetastar)$, that the set $\GenCons$ is
$\ThetaStar$-faithful and convex, and that the matrix $ \nabla
\EmpLoss(\PlTheta) $ is symmetric for any symmetric matrix $ \PlTheta
$. As we show in the proof of Theorem~\ref{ThmGeneralL}, the loss
function $ \LossTil $ then satisfies the corresponding conditions with
the same parameters. Consequently, we have the following result:
\begin{thm}
\label{ThmGeneralL}
 Under the previously stated conditions, the conclusion in
 Theorem~\ref{ThmGeneralSmooth} holds.
\end{thm}
\noindent See Section~\ref{SecProofThmGeneralL} for the proof of this
claim. \\

In remainder of this section, we verify the above conditions for each
of our examples. It is easy to see that in these examples the matrix $
\nabla \EmpLoss(\PlTheta) $ is indeed symmetric for any symmetric $
\PlTheta $, so it remains to verify the conditions in
Definitions~\ref{DefCurveCon}--\ref{DefSmoothCon} for $ \EmpLoss $ and
the $ \ThetaStar $-faithfulness of $ \GenCons $.

Recall that $ \sigma_i $ and $ \cond $ are the $ i $-th singular value
and the condition number of $\thetastar $, respectively. Throughout
this section, we let $\pltheta$ be an arbitrary matrix in $\GenCons
\cap \Ballstar{\smallrad}$ and $\PlTheta= \OUTER{\pltheta}$, where
$\GenCons$ and $\smallrad$ are specified for each of our examples. In
all these examples $ \smallrad < \sigma_\rdim $, so
Lemma~\ref{LemThetaStar} guarantees that we can write $\tdthetastar
\defn \arg\min_{A \in \eclass} \fronorm{A - \pltheta}$ and $\d \defn
\pltheta - \tdthetastar$, and $\d^{\top}\tdthetastar$ is a symmetric
matrix. Let $\thetastar$ be an arbitrary matrix in $\eclass$, and
recall that $\ThetaStar= \OUTER{\thetastar} = \OUTER{\tdthetastar}$.
Denote by $ C, c, c_1 $ etc.\ positive universal constants, whose
values could change from line to line.


\subsection{Proof of Corollary~\ref{CorMatrixSensing}}
\label{SecProofCorMatSensing}

We begin by proving our claims for the matrix sensing observation
model. By dividing through by $\sigma_\rdim$, we may assume without
loss of generality $\sigma_\rdim = 1$, so $ \cond = \sigma_1 $. Recall
that $\curv= 6\ripparam{4\rdim} $, $\LipCon= \smoo= 64 \cond^2$,
$\epsnum= \frac{2 \sqrt{\rdim} \cond
  \opnorm{\numobs^{-1}\Xmap^{*}(\smallnoisevar)}}{\ripparam{4\rdim}}$
and $\smallrad= 1-12\ripparam{4\rdim} $. It is a standard result that
RIP implies preservation of inner products between low rank matrices,
as summarized in the lemma below:
\begin{lem}
\label{LemMatSenInprod}
If $\Xmap$ satisfies a RIP-$\ripparam{4\rdim}$ condition, then
\begin{align*}
\frac{1}{\numobs} \Big|\inprod{\Xmap(A)}{\Xmap(B)} - \inprod{A}{B}
\Big| & \leq \ripparam{4\rdim} \fronorm A\fronorm B \quad \mbox{for
  all matrices \ensuremath{A,B \in \real^{\usedim\times \usedim}} of
  rank at most $2 \rdim$.}
\end{align*}
\end{lem}
\noindent For completeness, we provide a proof in
Appendix~\ref{SecProofLemMatSenInprod}. \\

\noindent Under the matrix sensing observation
model~\eqref{EqnMatrixSensing}, the gradient of $\EmpLoss$ takes the
form
\begin{align}
\label{EqnMatrixSensingGradient}
\nabla \EmpLoss(\OUTER{\pltheta}) & =
\frac{1}{\numobs}\Xmap^{*}\Xmap(\OUTER{\pltheta} - \OUTER{\thetastar})
- \frac{1}{\numobs}\Xmap^{*}(\smallnoisevar),
\end{align}
Below we verify the local descent, Lipschitz and smoothness
conditions.

\paragraph{Local descent:}

We have the decomposition $\trinprod{\nabla
  \EmpLoss(\PlTheta)}{\PlTheta - \ThetaStar + \d\otimes\d} = \Term_1 +
\Term_2$, where
\begin{align*}
\Term_1 & \defn \frac{1}{\numobs}\trinprod{\Xmap^{*}\Xmap(\PlTheta -
  \ThetaStar)}{\PlTheta - \ThetaStar + \d\otimes\d} \qquad \mbox{and}
\quad \Term_2 \defn -
\frac{1}{\numobs}\trinprod{\Xmap^{*}(\smallnoisevar)}{\PlTheta -
  \ThetaStar + \d\otimes\d}.
\end{align*}
Lemma~\ref{LemMatSenInprod} implies that
\begin{align*}
\Term_{1} & = \frac{1}{\numobs}\trinprod{\Xmap(\PlTheta -
  \ThetaStar)}{\Xmap(\PlTheta - \ThetaStar + \d\otimes\d)} \\
& \geq \trinprod{\PlTheta - \ThetaStar}{\PlTheta - \ThetaStar +
  \d\otimes\d} - \ripparam{4\rdim} \fronorm{\PlTheta - \ThetaStar}
\fronorm{\PlTheta - \ThetaStar + \d\otimes\d}.
\end{align*}
Since the matrix $\d^{\top}\tdthetastar$ is symmetric, some algebra
shows that
\begin{align*}
\trinprod{\PlTheta - \ThetaStar}{\PlTheta - \ThetaStar + \d\otimes\d}
& = \trinprod{\tdthetastar\otimes\d + \d\otimes\tdthetastar +
  \d\otimes\d}{\tdthetastar\otimes\d + \d\otimes\tdthetastar +
  2\d\otimes\d} \\ 
& = 2\fronorm{\tdthetastar\otimes\d}^{2} +
2\trinprod{\tdthetastar\otimes\d}{\d\otimes\d} +
2\fronorm{(\tdthetastar)^{\top}\d + \d^{\top}\d}^{2} \\
 & \geq 2 \fronorm{\tdthetastar\otimes\d} \big( \fronorm{\tdthetastar\otimes\d} - \fronorm{\d}^2 \big) 
\end{align*}
In addition, we have
\begin{align*}
\fronorm{\PlTheta - \ThetaStar} \fronorm{\PlTheta - \ThetaStar +
  \d\otimes\d} & = \fronorm{\tdthetastar\otimes\d +
  \d\otimes\tdthetastar + \d\otimes\d} \fronorm{\tdthetastar\otimes\d
  + \d\otimes\tdthetastar + 2\d\otimes\d}\\
& \leq \big( 2\fronorm{\tdthetastar\otimes\d} + \fronorm{\d}^2 \big) \big( 2\fronorm{\tdthetastar\otimes\d} + 2\fronorm{\d}^2 \big).
\end{align*}
It follows that 
\begin{align*}
T_1 &\geq \fronorm{\tdthetastar\otimes\d} \Big( (2-4\ripparam{4\rdim}) \fronorm{\tdthetastar\otimes\d}- (2+6\ripparam{4\rdim}) \fronorm{\d}^2 \Big) - 2\ripparam{4\rdim} \fronorm{\d} ^4 \\
& \geq  \fronorm{\d} \Big( (2-4\ripparam{4\rdim}) \fronorm{\d}- (2+6\ripparam{4\rdim}) (1-12\ripparam{4\rdim}) \fronorm{\d} \Big) - 2\ripparam{4\rdim} \fronorm{\d}^2 
 \geq   12\ripparam{4\rdim} \fronorm{\d}^2,
\end{align*}
where the second step uses the inequalities  $ \fronorm{\tdthetastar\otimes\d}  \ge  \sigma_\rdim \fronorm{\d} = \fronorm{\d}$ and $\fronorm{\d} \le \smallrad = 1 - 12\ripparam{4\rdim}$.

On the other hand, we have
\begin{align*}
\abs{\Term_{2}} 
%
\leq & \opnorm{\numobs^{-1}\Xmap^{*}(\smallnoisevar)}\cdot\sqrt{2\rdim} \fronorm{\PlTheta - \ThetaStar + \d\otimes\d} \\
\leq &
\opnorm{\numobs^{-1}\Xmap^{*}(\smallnoisevar)} \cdot  \sqrt{2 \rdim} \big(2\fronorm{\tdthetastar\otimes\d} + 2\fronorm{\d}^2 \big) 
\leq  \opnorm{\numobs^{-1}\Xmap^{*}(\smallnoisevar)} \cdot  6\sqrt{\rdim} \cond \fronorm{\d},
\end{align*}
where the last step uses the inequalities $ \fronorm{\tdthetastar\otimes\d} \le \sigma_1(\tdthetastar)\fronorm{\d} = \cond \frobnorm{\d}$ and $ \fronorm{\d} \le 1$. Combining this upper bound with our lower bound on $\Term_{1}$, we find that
\begin{align*}
\trinprod{\nabla \EmpLoss(\PlTheta)}{\PlTheta - \ThetaStar +
  \d\otimes\d} 
\geq 12\ripparam{4\rdim} \Big( \fronorm{\d}^{2} -
\frac{\sqrt{\rdim} \cond}{2\ripparam{4\rdim}}
\opnorm{\numobs^{-1}\Xmap^{*}(\smallnoisevar)} \fronorm{\d} \Big),
\end{align*}
thereby establishing the local descent condition~(\ref{EqnCurveCond}) for $ \EmpLoss $.


\paragraph{Local Lipschitz and smoothness:}

We have the following variational representation:
\begin{align*}
\fronorm{\nabla \LossTil(\pltheta)} & = \sup_{ \substack{\plphi \in
    \real^{\usedim \times \rdim} \\ \fronorm{\plphi} = 1}}
\trinprod{\nabla \EmpLoss(\OUTER{\pltheta})}{\plphi \otimes \pltheta}.
\end{align*}
Using the form of  the gradient $\nabla \EmpLoss$ given in
equation~\eqref{EqnMatrixSensingGradient}, we have
\begin{align*}
\nabla \EmpLoss(\OUTER{\pltheta}) - \nabla
\EmpLoss(\OUTER{\pltheta'})= \frac{1}{\numobs} \Xmap^{*}\Xmap \Big(
\OUTER{\pltheta} - \OUTER{\pltheta'} \Big).
\end{align*}
Note moreover that $0 \in \Ballstar{1}$, $\curv\le \LipCon= \smoo$ and
$\epsnum \leq 1$.  Using these facts, it can be verified that the
local Lipschitz and smoothness conditions~(\ref{EqnLipCon})--(\ref{EqnSmoothCond2}) for $\EmpLoss $  are implied by a bound of
the form
\begin{equation}
\label{EqnMatSenSmooth}
\abs{\trinprod{ \numobs ^{-1} \Xmap^{*}\Xmap(\OUTER{\pltheta} - \OUTER{\pltheta'})}{\plphi\otimes\plomega}} +
\abs{\trinprod{ \numobs ^{-1} \Xmap^{*}(\smallnoisevar)}{\plphi \otimes \plomega}}
\le \frac{1}{8\cond} \Big( \smoo\fronorm{\pltheta - \pltheta'} +
\curv\epsnum \Big) \fronorm{\plphi} \opnorm{\plomega},
\end{equation}
valid for all $\pltheta,\pltheta' \in \Ballstar{1}$ and for all
$\plphi,\plomega \in \real^{\usedim\times \rdim}$.

Let us prove the bound~\eqref{EqnMatSenSmooth}.
Lemma~\ref{LemMatSenInprod} guarantees that
\begin{align*}
\abs{\trinprod{ \numobs ^{-1} \Xmap^{*}\Xmap(\OUTER{\pltheta} - \OUTER{\pltheta'})}{\plphi\otimes\plomega}} 
& \leq (1 + \ripparam{4\rdim}) \fronorm{\OUTER{\pltheta} - \OUTER{\pltheta'}}\cdot\fronorm{\plphi} \opnorm{\plomega} \\
& \leq  (1 + \ripparam{4\rdim}) \big(\opnorm{\pltheta}+\opnorm{\pltheta'} \big) \fronorm{\pltheta - \pltheta'} \cdot\fronorm{\plphi} \opnorm{\plomega} \\
& \leq 8 \cond \fronorm{\pltheta - \pltheta'} \fronorm{\plphi}
\opnorm{\plomega},
\end{align*}
where the last step follows from the facts that
$\opnorm{\pltheta}  \leq \opnorm{\thetastar} + \dist{\pltheta} \leq 2\cond $ for all $\pltheta\in \Ballstar 1$ and $\ripparam{4\rdim} \leq 1$. 
We also have
\begin{align*}
\abs{\trinprod{ \numobs ^{-1} \Xmap^{*}(\smallnoisevar)}{\plphi\otimes\plomega}} \le
\opnorm{\numobs^{-1}\Xmap^{*}
  (\smallnoisevar)}\cdot\nucnorm{\plphi\otimes\plomega} \le
 \opnorm{\numobs^{-1}
  \Xmap^{*}(\smallnoisevar)} \cdot \sqrt{\rdim}\fronorm{\plphi} \opnorm{\plomega}.
\end{align*}
Combining these inequalities and recalling the values of
$\curv,\smoo,\epsnum$ yields the claim~\eqref{EqnMatSenSmooth}.


\subsection{Proof of Corollary~\ref{CorMC}}
\label{SecProofCorMC}

We now turn to the proof of our claims for the matrix completion
model.	By dividing through by $\opnorm{\thetastar}$ and using the equal eigenvalue assumption, we may assume
without loss of generality $\opnorm{\thetastar} = \sigma_\rdim(\thetastar) = 1$. We first show
that $\GenCons$ is \mbox{$\ThetaStar$-faithful.} Note that $\GenCons$
is the set of matrices with each row in the $\ell_{2}$ ball of radius
$\gamma \defn \sqrt{\frac{2\inco}{\usedim \rdim }} \opnorm{\thetait{0}}$.
Because $\thetait{0} \in \Ball_{2}(\frac{1}{5};\thetastar)$, we have
$\opnorm{\thetait{0}} \ge \frac{4}{5} \opnorm{\thetastar}$, whence
$\gamma\ge \sqrt{\frac{\inco}{\usedim \rdim}} \opnorm{\thetastar}$. Combined with the definition of the incoherence
parameter $\inco$, we see that any matrix $\thetastar \in \eclass$
satisfies the maximum row norm bound $\|\thetastar\|_{2, \infty} \leq
\gamma$, so that $\thetastar \in \GenCons$ as desired.

For future reference, we make note of a useful property satisfied by any matrices $ \pltheta  \in \GenCons \cap \Ballstar{\smallrad} $ and $\thetastar \in \eclass$. As a consequence of the clipping operation~$\ProjGenCons$, the row norms of the matrices $\pltheta$ and $\pltheta - \thetastar$ satisfy the bounds
\begin{subequations}
\begin{align}
\label{EqDtwoinf}
\twoinfnorm{\pltheta} & \le \sqrt{\frac{2\inco \rdim}{\usedim}}
\opnorm{\thetait{0}} \le 2\sqrt{\frac{\inco\rdim}{\usedim}}
\opnorm{\thetastar} \le 2\sqrt{\frac{\inco\rdim}{\usedim}} \quad
\text{and} \\
\label{EqDtwoinfB}
\twoinfnorm{\pltheta - \thetastar} & \le \twoinfnorm{\pltheta} +
\twoinfnorm{\thetastar}
\le 3\sqrt{\frac{\inco\rdim}{\usedim}}
\end{align}
\end{subequations}
where we use the inequality $\opnorm{\thetait{0}} \le
\opnorm{\thetastar} + \dist{\thetait{0}} \le \frac{6}{5}
\opnorm{\thetastar}$ and the normalization assumption $
\opnorm{\thetastar} = 1$. Inequality~\eqref{EqDtwoinfB} applies in
particular to the difference matrix \mbox{$\d \defn \pltheta -
  \tdthetastar$,} where we recall that the matrix $ \tdthetastar \defn
\arg \min_{A \in \eclass} \fronorm{A - \pltheta} $ is uniquely defined
thanks to Lemma~\ref{LemThetaStar}.

It remains to verify that the local descent, Lipschitz and
smoothness conditions are satisfied with high probability. Under the
matrix completion observation model~\eqref{EqnMC}, the gradient takes
the form $\nabla \EmpLoss(\PlTheta) = \frac{1}{\pobs}\ProjObs(\PlTheta
- \ThetaStar - \noisevar)$, where $\ProjObs$ is the projection
operator. We need two technical lemmas. The first lemma shows that the
projection operator $\ProjObs$ approximately preserves inner products
between matrices whose column or row spaces are equal to the column
space of $\thetastar$.
\begin{lem}
\label{LemMCInprod}
There is a universal constant $c$ such that for any $0 < \epsilon < 1$
and $\pobs \ge c\frac{\inco\rdim\log\usedim}{\epsilon^{2}\usedim}$,
uniformly for all $\plphi,\plomega \in \real^{\usedim\times \rdim}$,
we have
\begin{subequations}
\begin{align}
\label{EqnMCInprod1}
\abs{\pobs^{-1}\trinprod{\ProjObs\left(\thetastar\otimes\plphi\right)}{\ProjObs(\plomega\otimes\thetastar)}
  - \trinprod{\thetastar\otimes\plphi}{\plomega\otimes\thetastar}} &
\le \epsilon\opnorm{\thetastar}^{2} \fronorm{\plphi}
\fronorm{\plomega}, \\
\label{EqnMCInprod2}
\abs{\pobs^{-1}\trinprod{\ProjObs\left(\thetastar\otimes\plphi\right)}{\ProjObs(\thetastar\otimes\plomega)}
  - \trinprod{\thetastar\otimes\plphi}{\thetastar\otimes\plomega}} &
\le \epsilon\opnorm{\thetastar}^{2} \fronorm{\plphi}
\fronorm{\plomega},
\end{align}
\end{subequations}
with probability at least $1- 2 \usedim^{-3}$.
\end{lem}
\noindent See Appendix~\ref{SecProofLemMCInProd} for the proof of this
claim.\\

Our second lemma is useful for controlling the projection of ``small''
matrices to $\Obs$.
\begin{lem}
\label{LemMCProjSmall}
There is a universal constant $c > 0$ such that for any $\epsilon \in
(0,1)$ and $\pobs \geq \frac{C}{\epsilon^{2}} \big(
\frac{\inco^{2}\rdim^{2}}{\usedim} + \frac{\log\usedim}{\usedim}
\big)$, then unformly for all matrices $\PlZ \in \real^{\usedim\times
  \usedim}$, $\plomega \in \real^{\usedim\times \rdim}$ and $\plphi$
with $\twoinfnorm{\plphi} \leq 6 \sqrt{\frac{\inco\rdim}{\usedim}}$,
we have
\begin{subequations}
\begin{align}
\label{EqnMCProjSmall1}
\pobs^{-1} \fronorm{\ProjObs(\plphi\otimes\plphi)}^{2} & \le (1 +
\epsilon) \fronorm{\plphi}^{4} + \epsilon\fronorm{\plphi}^{2}, \\
\label{EqnMCProjSmall2}
\pobs^{-1} \fronorm{\ProjObs(\PlZ)\plphi}^{2} &
\le72\inco\rdim\fronorm{\ProjObs(\PlZ)}^{2}, \\
\label{EqnMCProjSmall3}
\pobs^{-1} \fronorm{\ProjObs(\phi\otimes\omega)}^{2} &
\le72\inco\rdim\fronorm{\plomega}^{2}
\end{align}
\end{subequations}
with probability at least $1 - 2 \usedim^{-4}$.
\end{lem}
\noindent
See Section~\ref{SecProofLemMCProjSmall} for the proof of this
claim. \\

\vsmall

For the remainder of the proof, we condition on the
intersection of the events in Lemmas~\ref{LemMCInprod}
and~\ref{LemMCProjSmall}.


\paragraph{Local descent:}
We have the decomposition $\trinprod{\nabla
  \EmpLoss(\PlTheta)}{\PlTheta - \ThetaStar + \OUTER{\d}} = \Term_1 -
\Term_2$, where
\begin{align*}
\Term_1 & \defn \frac{1}{\pobs}\trinprod{\ProjObs(\PlTheta -
  \ThetaStar)}{\ProjObs(\PlTheta - \ThetaStar + \OUTER{\d})}, \qquad
\mbox{and} \quad \Term_2 \defn
\frac{1}{\pobs}\trinprod{\ProjObs(\noisevar)}{\PlTheta - \ThetaStar +
  \OUTER{\d}}.
\end{align*}
Our strategy is to lower bound $\Term_1$ and upper bound $|\Term_2|$.
Beginning with $\Term_{1}$, we have
\begin{align*}
\Term_{1} & = \pobs^{-1}\trinprod{\ProjObs(\tdthetastar\otimes\d +
  \d\otimes\tdthetastar + \d\otimes\d)}{\ProjObs(\tdthetastar\otimes\d +
  \d\otimes\tdthetastar + 2\d\otimes\d)} \\
& \geq \pobs^{-1} \Big( \fronorm{\ProjObs(\tdthetastar\otimes\d +
  \d\otimes\tdthetastar)}-2\fronorm{\ProjObs(\d\otimes\d)} \Big) \Big(
\fronorm{\ProjObs(\tdthetastar\otimes\d + \d\otimes\tdthetastar)} -
\fronorm{\ProjObs(\d\otimes\d)} \Big).
\end{align*}
Recall that $\pobs\ge \frac{C}{\epsilon^{2}} \Big(
\frac{\inco^{2}\rdim^{2}}{\usedim} + \frac{\inco\rdim\usedim}{\usedim}
\Big)$ by assumption of the corollary. By Lemma~\ref{LemMCInprod}, we find that
\begin{align*}
\pobs^{-1} \fronorm{\ProjObs(\tdthetastar\otimes\d +
  \d\otimes\tdthetastar)}^{2} \ge (1 - \epsilon)
\fronorm{\tdthetastar\otimes\d + \d\otimes\tdthetastar}^{2} \ge2(1 -
\epsilon) \fronorm{\d}^{2},
\end{align*}
where the last step follow from the inequality 
\[
\fronorm{\tdthetastar\otimes\d + \d\otimes\tdthetastar}^{2} 
= 2\fronorm{\tdthetastar \otimes \d }^2 + 2\trinprod{\tdthetastar\otimes\d}{\d\otimes\tdthetastar} = 2\fronorm{\tdthetastar \otimes \d }^2 + 2\fronorm{\d^\top \tdthetastar}^2 \ge 2 \fronorm{\d}^2
\] thanks to the symmetry of the matrix $\d^{\top}\tdthetastar$ (cf. Lemma~\ref{LemThetaStar}).
Since $\twoinfnorm{\d} \leq 3
\sqrt{\frac{\inco \rdim}{\usedim}}$ and $\fronorm{\d} \le
\frac{3}{5}$, we can use the inequality~\eqref{EqnMCProjSmall1} from
Lemma~\ref{LemMCProjSmall} to get
\begin{align}
\pobs^{-1} \fronorm{\ProjObs(\OUTER{\d})}^{2} & \le (1 + \epsilon)
\fronorm{\d}^{4} + \epsilon\fronorm{\d}^{2} \le \frac{9}{25}(1 +
4\epsilon) \fronorm{\d}^{2}. \label{EqnMC_DDBound}
\end{align}
With the constant $\epsilon$ sufficiently small, we get that
$2\fronorm{\ProjObs(\OUTER{\d})} \le
\fronorm{\ProjObs(\thetastar\otimes\d + \d\otimes\thetastar)}$ and
\begin{subequations}
\begin{align}
\label{EqnT1Lower}
\Term_{1} & \geq \fronorm{\d}^{2} \big( \sqrt{2(1 - \epsilon)} -
\frac{6}{5}\sqrt{1 + 4\epsilon} \big) \big( \sqrt{2(1 - \epsilon)} -
\frac{3}{5}\sqrt{1 + 4\epsilon} \big)\ge \frac{4}{25}
\fronorm{\d}^{2}.
\end{align}
On the other hand, we have
\begin{align}
\label{EqnT2Upper}
\abs{\Term_{2}} \le \frac{1}{\pobs}
\opnorm{\ProjObs(\noisevar)}\cdot\sqrt{\rdim} \fronorm{\PlTheta -
  \ThetaStar + \OUTER{\d}} \le \frac{4\sqrt{\rdim}}{\pobs}
\opnorm{\ProjObs(\noisevar)}\cdot\fronorm{\d}.
\end{align}
\end{subequations}

Combining inequalities~\eqref{EqnT1Lower} and~\eqref{EqnT2Upper} with
our original decomposition yields
\begin{align*}
\trinprod{\nabla \EmpLoss(\PlTheta)}{\PlTheta - \ThetaStar +
  \OUTER{\d}} \ge \frac{4}{25} \fronorm{\d}^{2} -
\frac{4\sqrt{\rdim}}{\pobs}
\opnorm{\ProjObs(\noisevar)}\cdot\fronorm{\d},
\end{align*}
showing that the local descent~(\ref{EqnCurveCond}) for $ \EmpLoss $ holds with $\curv= \frac{2}{25}$ and $\epsnum= \frac{100\sqrt{\rdim}}{\pobs}
\opnorm{\ProjObs(\noisevar)}.$

\paragraph{Local Lipschitz and smoothness:}

Observe that $\curv\le \LipCon= \smoo$ and $\max\{\smallrad, \epsnum \}\leq 1$, 
$\twoinfnorm{\pltheta - \pltheta'} \leq 6
\sqrt{\frac{\inco\rdim}{\usedim}}$ for all $\pltheta,\pltheta' \in
\GenCons$, and
\begin{align*}
\fronorm{\nabla \EmpLoss(\pltheta) \pltheta} = \sup_{\plomega \in
  \real^{\usedim\times \rdim}, \frobnorm{\plomega}\le 1}\trinprod{\nabla
  \EmpLoss(\OUTER{\pltheta})}{\plomega\otimes\pltheta}.
\end{align*}
Using these facts, it follows that the Lipschitz and
smoothness conditions~(\ref{EqnLipCon})--(\ref{EqnSmoothCond2}) for $\EmpLoss $ can be verified by showing that
\begin{align}
\label{EqnMCSmoothness}
\Biggr| \frac{\trinprod{\ProjObs(\OUTER{\pltheta} -
    \OUTER{\pltheta'})}{\plomega\otimes\plphi}}{\pobs} \Biggr| +
\Biggr |
\frac{\trinprod{\ProjObs(\noisevar)}{\plomega\otimes\plphi}}{\pobs}
\Biggr| \le \frac{1}{8} \big( \smoo\fronorm{\pltheta - \pltheta'} +
\curv\epsnum\opnorm{\plphi} \big) \fronorm{\plomega},
\end{align}
valid for all $\pltheta,\pltheta' \in \GenCons \cap \Ballstar{1}$, and
for all matrices $\plphi, \plomega \in \real^{\usedim\times \rdim}$
such that $\twoinfnorm{\plphi} \le6\sqrt{\frac{\inco\rdim}{\usedim}}$.

Let us now verify the bound~\eqref{EqnMCSmoothness}.  For an arbitrary
$\thetastar \in \eclass$, define the matrices $\diff= \pltheta' -
\pltheta$, $\d_{1} = \pltheta - \thetastar$ and $\d_{2} = \pltheta' -
\thetastar$, and observe that
\begin{align*}
\frac{\fronorm{\ProjObs(\OUTER{\pltheta} -
    \OUTER{\pltheta'})}}{\sqrt{\pobs}} & =
\frac{\fronorm{\ProjObs(\thetastar\otimes\diff +
    \diff\otimes\thetastar + \diff\otimes\d_{2} +
    \d_{1}\otimes\diff)}}{\sqrt{\pobs}} \\
& \le \underbrace{
  \frac{\fronorm{\ProjObs(\thetastar\otimes\diff)}}{\sqrt{\pobs}}}_{\Term_1}
+ \underbrace{ \frac{\fronorm{\ProjObs(\diff \otimes
      \d_2)}}{\sqrt{\pobs}}}_{\Term_2} + \underbrace{
  \frac{\fronorm{\ProjObs(\d_{1}\otimes\diff)}}{\sqrt{\pobs}}}_{\Term_3}.
\end{align*}
Lemma~\ref{LemMCInprod} implies that $\Term_{1} \le (1 + \epsilon)
\fronorm{\thetastar\otimes\diff + \diff\otimes\thetastar} \le 2(1 +
\epsilon) \fronorm{\diff}$, whereas inequality~\eqref{EqnMCProjSmall3}
from Lemma~\ref{LemMCProjSmall} ensures that
$\max\{\Term_{2},\Term_{3}\} \le6\sqrt{2\inco\rdim} \fronorm{\diff}$.
Combining these bounds yields
\begin{align}
\label{EqnMCObsLambda}
\frac{\fronorm{\ProjObs(\OUTER{\pltheta} -
    \OUTER{\pltheta'})}}{\sqrt{\pobs}} & \leq 14 \sqrt{2\inco\rdim}
\fronorm{\pltheta - \pltheta'}.
\end{align}
On the other hand, using the inequality~(\ref{EqnMCProjSmall2}) from
Lemma~\ref{LemMCProjSmall}, we have
\begin{align*}
\abs{\pobs^{-1}\trinprod{\ProjObs(\OUTER{\pltheta} - \OUTER{\pltheta}')}{\plomega\otimes\plphi}} & \le \pobs^{-1/2}\cdot\fronorm{\pobs^{-1/2}\ProjObs(\OUTER{\pltheta} - \OUTER{\pltheta'})\plphi} \fronorm{\plomega} \\
 & \le6\sqrt{\frac{2\inco\rdim}{\pobs}} \fronorm{\ProjObs(\OUTER{\pltheta} - \OUTER{\pltheta'})} \fronorm{\plomega}.
\end{align*}
Combining with the earlier inequality~\eqref{EqnMCObsLambda} yields
\begin{align*}
\abs{\pobs^{-1} \trinprod{\ProjObs(\OUTER{\pltheta} -
    \OUTER{\pltheta}')}{\plomega\otimes\plphi}} \leq 168 \inco \rdim
\fronorm{\pltheta - \pltheta'} \fronorm{\plomega}.
\end{align*}
Finally, observe that
\begin{align*}
\abs{\pobs^{-1}\trinprod{\ProjObs(\noisevar)}{\plomega \otimes
    \plphi}} \leq \frac{1}{\pobs} \opnorm{ \ProjObs(\noisevar)}
\nucnorm{\plomega \otimes \plphi} \leq \frac{\sqrt{\rdim}}{\pobs}
\opnorm{\ProjObs(\noisevar)} \fronorm{\plomega} \opnorm{\plphi}.
\end{align*}
Combining the last two inequalities establishes the
claim~\eqref{EqnMCSmoothness}, thereby completing the proof of
Corollary~\ref{CorMC}.


\subsection{Proof of Corollary~\ref{CorSparsePCA}}
\label{SecProofCorSparsePCA}

We now prove our claims for the sparse PCA model.  Define the sampling
noise matrix $\noisemat \defn \SigHat - \CovMat$, corresponding to the
deviation between the sample and population covariance matrices.
Recall that $\CovMat= \snr \, \big(\OUTER{\thetastar}\big) +
\I_{\usedim}$ with $\opnorm{\thetastar} = 1$. Under the spiked
covariance model~\eqref{EqnSpikedCovariance}, we have \mbox{$\nabla
  \EmpLoss(\Theta)= - \SigHat= -(\CovMat + \noisemat)$.}  Let $\Tset$
index the non-zero rows of $\thetastar$; observe that the choice of
$\Tset$ does not depend on the choice of $\thetastar$ in $\eclass$.

In light of Remark~\ref{RemThetaStarRowNorm}, we have $\eclass
\subseteq \GenCons$, which guarantees the $\ThetaStar$-faithfulness
condition.  It remains to verify the local descent, Lipschitz and
smoothness conditions.

\paragraph{Local descent:}
For a given matrix $\pltheta$, let $\tdthetastar$ be its projection
onto $\eclass$, and deifne $\d = \pltheta - \tdthetastar$.  Since
$\OUTER{\tdthetastar} = \OUTER{\pltheta^*}$, we have
\begin{align*}
\trinprod{\nabla \EmpLoss(\PlTheta)}{\PlTheta - \ThetaStar +
  \OUTER{\d}} & = \underbrace{- \trinprod{\CovMat}{\PlTheta -
    \ThetaStar + \OUTER{\d}}}_{\Term_1} \underbrace{-
  \trinprod{\noisemat}{\PlTheta - \ThetaStar + \OUTER{\d}}}_{\Term_2}.
\end{align*}
The remainder of the proof consists of lower bounding $T_1$ and $T_2$.\\

Beginning with $\Term_1$, observe that
\begin{align*}
\Term_{1} & = - \trinprod{\snr\OUTER{\tdthetastar} +
  \I}{\tdthetastar\otimes\d + \d\otimes\tdthetastar + 2\d\otimes\d} \\
& = -2 \snr \big( \inprod{\tdthetastar}{\d} +
\fronorm{\d^{\top}\tdthetastar}^{2}
\big)-2\trinprod{\tdthetastar}{\d}-2\fronorm{\d}^{2}.
\end{align*}
By Lemma~\ref{LemThetaStar}, the matrix $\pltheta^{\top}\tdthetastar$
is positive semidefinite, and has operator norm bounded as
$\opnorm{\pltheta^{\top}\tdthetastar} \le \opnorm{\pltheta}
\opnorm{\tdthetastar} \le 1$, so that the matrix $ -
\d^{\top}\tdthetastar= \I_{\rdim} - \pltheta^{\top}\tdthetastar$ is
also positive semidefinite. We therefore have the bound
\begin{align*}
\fronorm{\d^{\top}\tdthetastar}^{2} \le
\nucnorm{\d^{\top}\tdthetastar} \opnorm{\d^{\top}\tdthetastar} = -
\trace(\d^{\top}\tdthetastar)\cdot\opnorm{\d^{\top}\tdthetastar}.
\end{align*}
Combined with the bound $\opnorm{\pltheta} \le 1$ and the
orthonormality of $\tdthetastar$, we find that
\begin{align*}
- \inprod{\tdthetastar}{\d} & = - \inprod{\tdthetastar}{\pltheta} +
\fronorm{\tdthetastar}^{2} \ge - \inprod{\tdthetastar}{\pltheta} +
\frac{1}{2} \fronorm{\tdthetastar}^{2} + \frac{1}{2}
\fronorm{\pltheta}^{2} = \frac{1}{2} \fronorm{\d}^{2}.
\end{align*}
It follows that
\begin{align*}
\Term_{1} & \ge \snr\fronorm{\d}^{2} \big( 1 -
\opnorm{\d^{\top}\tdthetastar} \big) +
\fronorm{\d}^{2}-2\fronorm{\d}^{2} \ge \fronorm{\d}^{2} \big(
\snr\inisig^{2}-1 \big),
\end{align*}
where in the last inequality we use $\opnorm{\d^{\top}\tdthetastar}
\le \opnorm{\d} \leq 1 - \inisig^{2}$.	Combined with the assumption
$\snr \geq \frac{2}{\inisig^{2}}$, it thus follows that $\Term_{1}
\geq \frac{\snr\inisig^{2}}{2} \fronorm{\d}^{2}$.

In order to bound $\Term_{2}$, we require control on how the matrix
$\noisemat$ behaves when acting on matrices in the set $\Ckset{\kdim}
\defn \{\plu \in \real^{\usedim\times \rdim} \,\mid\,
\twoonenorm{\plu} \le \sqrt{k} \fronorm{\plu}\}$.
\begin{lem}
\label{LemREofW}
There is a universal constant $c > 0$ such that
\begin{align}
\label{EqnREofW}
\abs{\trinprod{\noisemat}{\plu\otimes\plv}} \le c \; \max
\Big\{\sqrt{\frac{\kdim\log\usedim}{\numobs}},\frac{\kdim\log\usedim}{\numobs}
\Big\}(\snr + 1) \fronorm{\plu} \fronorm{\plv}, \quad \mbox{for all
  $\plu, \plv \in \Ckset{\kdim}$}
\end{align}
with probability at least $1-2\usedim^{-4}$.
\end{lem}
\noindent See Appendix~\ref{SecProofLemREofW} for the proof of this
lemma; it is based on variants of techniques from Lemma~12 of Loh and
Wainwright~\cite{loh2012nonconvex}.

Now observe that the row sparsity of $\thetastar$ implies $\thetastar
\in \Ckset{\kdim}$, $\forall\thetastar \in \eclass$. Recall that
$\Tset$ is the row support set of $\thetastar$, with $\TsetComp$
denoting its complement. Since $\pltheta \in \GenCons$, we are
guaranteed that $\twoonenorm{\pltheta} \leq\twoonenorm{\thetastar}$,
which implies the cone inequality $\twoonenorm{\d_{\TsetComp}}
\leq\twoonenorm{\d_{\Tset}}$.  By assumption $|\Tset|\leq\kdim$,
whence $\twoonenorm{\d} \leq\sqrt{\kdim} \fronorm{\d}$.  It follows
that $\d \in \Ckset{\kdim}$. Applying Lemma~\ref{LemREofW}, we find
that with probability at least $1-8\usedim^{-4}$,
\begin{align*}
|\Term_{2}| \; = \; \big| 2\trinprod{\noisemat}{\tdthetastar\otimes\d
  + 2 \d \otimes \d} \big| & \leq
2\abs{(\tdthetastar)^{\top}\noisemat\d} + 2\abs{\d^{\top}\noisemat\d}
\\
& \leq 4 c \max
\Big\{\sqrt{\frac{\kdim\log\usedim}{\numobs}},\frac{\kdim\log\usedim}{\numobs}
\Big\}(\snr + 1)\sqrt{\rdim} \fronorm{\d}.
\end{align*}
Combining the bounds for $\Term_{1}$ and $\Term_{2}$ proves that the
local descent condition~(\ref{EqnCurveCond}) for $ \EmpLoss $ is satisfied.

\paragraph{Local Lipschitz}

Let us verify the relaxed Lipschitz condition~\eqref{EqnWeakLipCon}.
Observe that for all matrices $\pltheta \in \GenCons \cap
\Ballstar{\smallrad}$, $\pltheta' \in \GenCons$ and $\thetastar \in
\eclass$, we have $\pltheta= \thetastar + (\pltheta - \thetastar)$ and
\begin{align*}
\pltheta - \pltheta'= (\pltheta - \thetastar)-(\pltheta' -
  \thetastar).
\end{align*}
Following the argument above one can show that the three matrices
$\thetastar,\pltheta - \thetastar,\pltheta' - \thetastar$ all belong
to the set $\Ckset{\kdim}$.  Consequently, Lemma~\ref{LemREofW}
guarantees that with probability at least $1-8\usedim^{-4}$,
\begin{align*}
\abs{\trinprod{\noisemat}{(\pltheta - \pltheta')\otimes\pltheta}} &
\leq 2 c \max
\Big\{\sqrt{\frac{\kdim\log\usedim}{\numobs}},\frac{\kdim\log\usedim}{\numobs}
\Big\}(\snr + 1)\cdot \big( \fronorm{\thetastar} + \fronorm{\pltheta -
  \thetastar} \big) \big( \fronorm{\pltheta - \thetastar} +
\fronorm{\pltheta' - \thetastar} \big) \\
& \leq \underbrace{12 c \max
  \Big\{\sqrt{\frac{\kdim\log\usedim}{\numobs}},\frac{\kdim\log\usedim}{\numobs}
  \Big\}(\snr + 1)\cdot\sqrt{\rdim}}_{ \frac{3}{4}\curv\epsnum},
\end{align*}
where in the last inequality we use $\fronorm{\pltheta' - \thetastar}
\le \sqrt{\rdim}(\opnorm{\pltheta'} + \opnorm{\thetastar})\le
2\sqrt{\rdim}.$ It follows that
\begin{align*}
\abs{\trinprod{\nabla \EmpLoss(\pltheta)}{(\pltheta - \pltheta')\otimes\pltheta}} & \leq
\abs{\trinprod{\CovMat}{(\pltheta - \pltheta')\otimes\pltheta}} +
\abs{\trinprod{\noisemat}{(\pltheta - \pltheta')\otimes\pltheta}} \\
& \leq \opnorm{\CovMat}\sqrt{\rdim} \fronorm{\pltheta - \pltheta'} \opnorm{\pltheta} + \frac{3}{4}\curv\epsnum \\
& \leq 2(\snr + 1)\sqrt{\rdim} \fronorm{\pltheta - \pltheta'} + \frac{3}{4} \curv\epsnum \\
& \leq 2(\snr + 1)\sqrt{\rdim}(1 + \fronorm{\pltheta - \pltheta'}),
\end{align*}
where the last inequality follows from $\epsnum\le 1$ and
$\curv\le 2(\snr + 1)\sqrt{\rdim}$. Thus, we have established the
relaxed Lipschitz condition~\eqref{EqnWeakLipCon} for $ \EmpLoss $.

\paragraph{Local smoothness: }

Since $\nabla \EmpLoss( \OUTER{\pltheta} ) - \nabla \EmpLoss( \OUTER{\pltheta'} )= 0$, the
first smoothness condition~\eqref{EqnSmoothCond1} for $ \EmpLoss $is satisfied
trivially.  On the other hand, we have
\begin{align*}
\abs{\trinprod{\nabla \EmpLoss(\PlTheta)}{(\pltheta -
    \thetastar)\otimes(\pltheta' - \pltheta'')}} & \leq
\abs{\trinprod{\CovMat}{(\pltheta - \thetastar)\otimes(\pltheta' -
    \pltheta'')}} + \abs{\trinprod{\noisemat}{(\pltheta -
    \thetastar)\otimes(\pltheta' - \pltheta'')}} \\
& \leq (\snr + 1)\sqrt{\rdim} \fronorm{\pltheta - \thetastar}
\fronorm{\pltheta' - \pltheta''} + \abs{\trinprod{\noisemat}{(\pltheta
    - \thetastar)\otimes(\pltheta' - \pltheta'')}}.
\end{align*}
Following the same argument above, we can show that $\pltheta -
\thetastar,\pltheta' - \thetastar,\pltheta'' - \thetastar \in
\Ckset{\kdim}$, whence Lemma~\ref{LemREofW} guarantees that
\begin{align*}
\abs{\trinprod{\noisemat}{(\pltheta - \thetastar)\otimes(\pltheta' -
    \pltheta'')}} \leq & c (\snr + 1)\max
\Big\{\sqrt{\frac{\kdim\log\usedim}{\numobs}},\frac{\kdim\log\usedim}{\numobs}
\Big\}\cdot\fronorm{\pltheta - \thetastar} \big( \fronorm{\pltheta' -
  \thetastar} + \fronorm{\pltheta'' - \thetastar} \big) \\
\stackrel{(i)}{\leq} & 2 c (\snr + 1)\sqrt{\rdim}\max
\Big\{\sqrt{\frac{\kdim\log\usedim}{\numobs}},\frac{\kdim\log\usedim}{\numobs}
\Big\}\cdot\fronorm{\pltheta - \thetastar}\\
= &\frac{1}{8}\curv\epsnum,
\end{align*}
with probability at least $1-8\usedim^{-4}$.  Here step (i) follows
from the inequality $\max \{ \fronorm{\pltheta' - \thetastar},
\fronorm{\pltheta'' - \thetastar} \} \le \sqrt{\rdim}$, valid for any
pair $\pltheta',\pltheta'' \in \Ballstar{\smallrad}$. We conclude that
\begin{align*}
\abs{\trinprod{\nabla \EmpLoss(\PlTheta)}{(\pltheta -
    \thetastar)\otimes(\pltheta' - \pltheta'')}} \le (\snr +
1)\sqrt{\rdim} \fronorm{\pltheta - \thetastar} \fronorm{\pltheta' -
  \pltheta''} + \frac{1}{8}\curv\epsnum\fronorm{\pltheta - \thetastar},
\end{align*}
thereby establishing the second smoothness condition~\eqref{EqnSmoothCond2} for $ \EmpLoss $.


\subsection{Proof of Corollary~\ref{CorClustering}}
\label{SecProofCorClustering}

We now prove our claims for the planted densest subgraph model.  Since
$\eclass$ is a two-element set, it follows that for any vector
$\pltheta \in \GenCons \cap \Ballstar{\smallrad}$, the projection
$\tdthetastar \defn \arg \min \limits_{\thetastar \in \eclass}
\fronorm{\thetastar - \pltheta}$ is always equal to the cluster
membership vector.  The $\ThetaStar$-invariance of the set $\GenCons$
thus follows. Under the planted densest subgraph model, the
expectation of the shifted adjacency matrix $\Shift$ has the
expression
\begin{align*}
\ShiftBar \defn & \Exs[\Shift]= \frac{\pin - \qout}{2}
\Big\{2\OUTER{\thetastar} - \OUTER{\onevec} \Big\},
\end{align*}
where $\onevec \in \real^\usedim$ denotes a vector of all ones.  The
noise matrix $\noisemat \defn \Shift - \ShiftBar$ has i.i.d. zero mean
entries with variance bounded by $\pin$. The gradient of $\EmpLoss$ is
given by $\nabla \EmpLoss(\OUTER{\pltheta})= -2\Shift\pltheta$. Below
we verify the local descent, Lipschitz and smoothness conditions.

\paragraph{local descent:}

We have the decomposition
\begin{align*}
\trinprod{\nabla \EmpLoss(\PlTheta)}{\PlTheta - \ThetaStar +
  \OUTER{\d}} & = \underbrace{-2
  \trinprod{\ShiftBar}{\tdthetastar\otimes\d + \d\otimes\d}}_{\Term_1}
\quad \underbrace{- 2\trinprod{\noisemat}{\tdthetastar\otimes\d +
    \d\otimes\d}}_{\Term_2}.
\end{align*}
We proceed by lower bounding $\Term_1$ and upper bounding $|\Term_2|$. \\

\vtiny

Beginning with the term $\Term_1$, for any feasible deviation $\d$,
the two matrices $-\ShiftBar$ and $\thetastar\d^{\top}$ have the same
sign on each entry, whence
\begin{align*}
-2\trinprod{\ShiftBar}{\thetastar\otimes\d} & = (\pin -
\qout)\onenorm{\thetastar\otimes\d} \ge (\pin -
\qout)\csize\onenorm{\d}.
\end{align*}
On the other hand, the bounds $\fronorm{\d} \leq\frac{2}{5}
\opnorm{\thetastar} = \frac{2\sqrt{\csize}}{5}$ and
$\onenorm{\pltheta} \le \csize= \onenorm{\thetastar}$ imply that
$\onenorm{\d} \le 2\sqrt{\csize} \fronorm{\d} \le \frac{4}{5}\csize$,
from which it follows that
\begin{align*}
\abs{2\trinprod{\ShiftBar}{\OUTER{\d}}} & \le
2\infnorm{\ShiftBar}\onenorm{\OUTER{\d}} = (\pin - \qout)
\onenorm{\d}^{2} \leq\frac{4}{5}(\pin - \qout)\csize\onenorm{\d}.
\end{align*}
Putting together the piecesr, we obtain the lower bound $\Term_{1} \ge
\frac{1}{5}(\pin - \qout)\csize\onenorm{\d}.$

Now turning to term $\Term_2$, by Bernstein's inequality and
Proposition~\ref{PropBandeira},
there is a universal constant  $c_{0} > 0$ such that
\begin{align*}
\infnorm{W\thetastar} \le c_{0}\sqrt{\pin\csize\log\usedim} \qquad
\text{and} \qquad \opnorm W\le c_{0}\sqrt{\pin\usedim +
  \pin\csize\log\usedim},
\end{align*}
with probability at least $1 - \usedim^{-3}$.  On this event, the term
$\Term_{2}$ can be bounded as
\begin{align*}
\abs{\Term_{2}} & \le 2 \big( \infnorm{W\thetastar}\onenorm{\d} +
\opnorm W\fronorm{\d}^{2} \big) \\
& \overset{}{\le}2c_{0} \big( \sqrt{\pin\csize\log\usedim}\onenorm{\d}
+ \sqrt{\pin\usedim + \pin\csize\log\usedim} \fronorm{\d}^{2} \big) \\
& \overset{(ii)}{\le} \frac{1}{10}\csize(\pin - \qout)\onenorm{\d},
\end{align*}
where the step (ii) follows from the clustering
condition~\eqref{EqnClusteringCond}, as well as the upper bound
$\fronorm{\d}^{2} \le \infnorm{d} \onenorm{\d} \leq \onenorm{d}$,
using the fact that $\infnorm{\d} \le1$.  Combining the bounds for
$\Term_{1}$ and $\Term_{2}$, we conclude that
\begin{align*}
\trinprod{\nabla \EmpLoss(\PlTheta)}{\PlTheta - \ThetaStar +
  \OUTER{\d}} \ge \frac{1}{10}\csize(\pin - \qout)\onenorm{\d} \ge
\frac{1}{10}\csize(\pin - \qout) \fronorm{\d}^{2},
\end{align*}
thereby establishing the local descent condition~\eqref{EqnCurveCond} for $ \EmpLoss $.


\paragraph{Local Lipschitz and smoothness:}

Since $\nabla \EmpLoss(\PlTheta) - \nabla \EmpLoss(\PlTheta')= 0$, the
first smoothness condition~\eqref{EqnSmoothCond1} is satisfied
trivially.  It remains to verify the second smoothness
condition~(\ref{EqnSmoothCond2}) and the relaxed Lipschitz
condition~(\ref{EqnWeakLipCon}). For each
$\pltheta,\pltheta',\pltheta'' \in \GenCons$ and $\plomega \in
\real^{\usedim\times \rdim}$, observe that
\begin{align*}
\abs{\trinprod{\nabla \EmpLoss(\PlTheta)}{(\pltheta -
    \plomega)\otimes(\pltheta' - \pltheta'')}} \leq
\underbrace{\abs{\trinprod{\ShiftBar}{(\pltheta -
      \plomega)\otimes(\pltheta' - \pltheta'')}}}_{\Term_1} +
\underbrace{\abs{\trinprod{\noisemat}{(\pltheta -
      \plomega)\otimes(\pltheta' - \pltheta'')}}}_{\Term_2}.
\end{align*}
Note that the matrices $\pltheta',\pltheta'' \in \GenCons$ satisfy the
constraint $\sum_{i}\pltheta'= \sum_{i}\pltheta''= \csize$, which
implies that $(\OUTER{\onevec})(\pltheta' - \pltheta'')= 0$. It
follows that $\Term_1$ can be upper bounded as
\begin{align*}
\Term_{1} & = \frac{\pin - \qout}{2} \; \Big|
\trinprod{2\OUTER{\thetastar} - \OUTER{\onevec}}{(\pltheta -
  \plomega)\otimes(\pltheta' - \pltheta'')} \Big| \\
& = (\pin - \qout) \; \Big | \trinprod{\OUTER{\thetastar}}{(\pltheta -
  \plomega)\otimes(\pltheta' - \pltheta'')} \Big| \\
& \leq (\pin - \qout) \fronorm{\thetastar}^{2} \fronorm{\pltheta -
  \plomega} \fronorm{\pltheta' - \pltheta''} \\
& = 2 (\pin - \qout)\csize\fronorm{\pltheta - \plomega}
\fronorm{\pltheta' - \pltheta''}.
\end{align*}
Similarly, the second term can be upper bounded as
\begin{align*}
\Term_{2} \le & \opnorm{\noisemat} \fronorm{\pltheta - \plomega}
\fronorm{\pltheta' - \pltheta''} \\ \overset{(i)}{\le} &
c_{0}\sqrt{\pin\usedim + \pin\csize\log\usedim} \fronorm{\pltheta -
  \plomega} \fronorm{\pltheta' - \pltheta''} \\
 \overset{(ii)}{\le} & \frac{1}{64}\csize(\pin - \qout)
 \fronorm{\pltheta - \plomega} \fronorm{\pltheta' - \pltheta''},
\end{align*}
where inequality (i) holds with probability at least $1 -
\usedim^{-3}$ as proved above, and inequality (ii) holds under the
clustering condition~\eqref{EqnClusteringCond}.  Combining the bounds
for $\Term_{1}$ and $\Term_{2}$ with the choice $\smoo= 12(\pin -
\qout)\csize$, we conclude that
\begin{align*}
\abs{\trinprod{\nabla \EmpLoss(\PlTheta)}{(\pltheta -
    \plomega)\otimes(\pltheta' - \pltheta'')}} \le \frac{\smoo}{4}
\fronorm{\pltheta - \plomega} \fronorm{\pltheta' - \pltheta''}, \qquad
\forall\pltheta,\pltheta',\pltheta'' \in \GenCons, \, \plomega \in
\real^{\usedim\times \rdim}.
\end{align*}
For an arbitrary $\thetastar \in \eclass$, setting $\plomega=
\thetastar$ in this inequality establishes the smoothness
condition~\eqref{EqnSmoothCond2}.  On the other hand, setting
$\plomega= 0$ and noting that $\fronorm{\pltheta}^{2} \le
\onenorm{\pltheta}\infnorm{\pltheta} \le \csize=
\opnorm{\thetastar}^{2}$, we obtain the relaxed Lipschitz
condition~\eqref{EqnWeakLipCon} for $ \EmpLoss $.


\subsection{Proof of Corollary~\ref{CorOB}}
\label{SecProofCorOB}

We now prove our claims for the one-bit matrix completion model.  By
assumption, the initial matrix $\thetait{0}$ belongs to the set
$\Ballstar{\frac{1}{5}} \cap \GenCons$, where the set $\GenCons$ is
was previously involved in our analysis of ordinary matrix completion
(see Section~\ref{SecProofCorMC}). Therefore, following the argument
in Remark~\ref{RemThetaStarRowNorm}, we can show that
$\twoinfnorm{\pltheta} \le 2\sqrt{\frac{\inco\rdim}{\usedim}}$,
$\twoinfnorm{\pltheta - \thetastar}
\le3\sqrt{\frac{\inco\rdim}{\usedim}}$ for all $\thetastar \in
\eclass$ and $\pltheta \in \Ballstar{\smallnoisevar} \cap \GenCons$,
and that $\GenCons$ is $\ThetaStar$-faithful. Consequently, we have
for all relevant matrices $\ThetaStar$ and $\PlTheta$, we are guaranteed
that
\begin{align*}
\max \Big \{ \infnorm{\ThetaStar/\noisestd}, \;
\infnorm{\PlTheta/\noisestd} \Big \} & \leq 4 \;
\underbrace{\frac{\mu\rdim}{\usedim\noisestd}}_{\obsnr}.
\end{align*}
Now define a (random function) $H: \real^{\usedim \times \usedim}
\rightarrow \real^{\usedim \times \usedim}$ with entries $[H(x)]_{ij}
\defn \frac{f'(x)\left( - \Yout_{ij} +
  2f(x)-1\right)}{f(x)\left(1-f(x)\right)}$. With this notation, we
have
\begin{align}
\label{EqnFactorGrad}
\nabla \EmpLoss(\PlTheta) & = \frac{1}{\noisestd} \ProjObs
\left[H(\PlTheta/\sigma)\right].
\end{align}
For future reference, we claim that each component $[H(\cdot)]_{ij}$
is bounded and $4\obupper{4\obsnr}$-Lipschitz over $[-4\nu,4\nu]$.
This property follows because $f$ satisfies the bounds in
equation~\eqref{EqnFlatness} and $\abs{\Yout_{ij}} \leq 1$, so that
$\abs{[H(x)]_{ij}} \le 2\obupper{4\obsnr}$ surely. Moreover, the
derivative can be bounded as
\begin{align*}
\abs{[H'(x)]_{ij}} = \Biggr| \frac{f''(x)(\Yout_{ij}-2f(x) +
  1)-2f'(x)^{2}}{f(x)(1-f(x))} - \frac{f'(x)(\Yout_{ij}-2f(x) +
  1)}{f^{2}(x)(1-f(x))^{2}}f'(x)(1-2f(x)) \Biggr| \leq 4
\obupper{4\obsnr}
\end{align*}
which certifies the Lipschitz property.\\

\vtiny

\noindent With this set-up, we are now prepared to establish the local
descent, Lipschitz and smoothness conditions for $ \EmpLoss $.

\paragraph{Local descent:}
Let us introduce the shorthand $\PopGrad(\PlTheta) \defn \Exs [ \nabla
  \EmpLoss(\PlTheta) ] = \nabla \Exs[\EmpLoss(\PlTheta)]$.  We begin
by splitting the gradient into two terms, corresponding to the
expectation and the zero-mean deviation, thereby obtaining
\begin{align*}
\trinprod{\nabla \EmpLoss(\PlTheta)}{\PlTheta - \ThetaStar +
  \OUTER{\d}} & = \underbrace{\trinprod{\PopGrad(\PlTheta)}{\PlTheta -
    \ThetaStar + \OUTER{\d}}}_{\Term_1} + \underbrace{\trinprod{\nabla
    \EmpLoss(\PlTheta) - \PopGrad(\PlTheta)}{\PlTheta - \ThetaStar +
    \OUTER{\d}}}_{\Term_2}.
\end{align*}

\vtiny

\noindent \emph{Controlling $\Term_1$:} For the expectation term
$\Term_{1}$, we first note that $\Exs[\EmpLoss(\PlTheta)]$ is convex
in $\PlTheta$ and has the form of expected negative log likelihood,
whence
\begin{align*}
\trinprod{\PopGrad(\PlTheta)}{\PlTheta - \ThetaStar} & \geq \Exs
	 [\EmpLoss(\PlTheta)] - \Exs[\EmpLoss(\ThetaStar)] \\
& = \pobs \; D \big( f(\ThetaStar/\noisestd)\Vert
	 f(\PlTheta/\noisestd) \big)\\
 & \geq \pobs \; d_{H}^{2} \big( f(\ThetaStar/\noisestd)\Vert
	 f(\PlTheta/\noisestd) \big),
\end{align*}
where $D(\cdot)$ and $d_{H}^{2}(\cdot)$ denote the KL and Hellinger
distances, respectively. To proceed, we use a known lower bound (Lemma
2 in the paper~\cite{davenport2014_onebitMC}) on the Hellinger
distance: for matrices $\PlX,\PlX' \in \real^{\usedim\times \usedim}$
such that $\infnorm X,\infnorm{X'} \le a$, we have
\begin{align}
\label{EqnHellingerLower}
d_{H}^{2} \big( f(X),f(X') \big) & \geq \frac{\fronorm{X-X'}^{2}}{8
  \oblower a}.
\end{align}
Since $\infnorm{\PlTheta/\noisestd},\infnorm{\ThetaStar/\noisestd}
\le4\obsnr$, applying the lower bound~\eqref{EqnHellingerLower} with
$a = 4 \obsnr$ yields the lower bound
$\trinprod{\PopGrad(\PlTheta)}{\PlTheta - \ThetaStar} \geq
\frac{\pobs}{8\noisestd^{2}\oblower{4\obsnr}} \fronorm{\PlTheta -
  \ThetaStar}^{2}$.  Furthermore, since $\thetastar$ is orthonormal,
$\opnorm{\d} \le \frac{1}{16}$ and $\d^{\top}\tdthetastar$ is
symmetric, we have
\begin{align*}
\trinprod{\PopGrad(\PlTheta)}{\PlTheta - \ThetaStar} \ge
\frac{\pobs}{16\noisestd^{2}\oblower{4\obsnr}} \fronorm{\d}^{2},
\end{align*}
On the other hand, using the expression~\eqref{EqnFactorGrad} for the
gradient and the assumptions in equation~\eqref{EqnFlatness}, we find
that
\begin{align*}
\fronorm{\PopGrad(\PlTheta)} & = \frac{\pobs}{\noisestd} \cdot
\fronorm{\frac{f'(\PlTheta/\noisestd)\circ(f(\PlTheta/\noisestd)-f(\ThetaStar/\noisestd))}{f(\PlTheta/\noisestd)
    \circ (1 - f(\PlTheta/\noisestd))}} \leq
\frac{\pobs\sqrt{\obupper{4\obsnr}}}{\noisestd} \cdot
\fronorm{f(\PlTheta/\noisestd)-f(\ThetaStar/\noisestd)}.
\end{align*}
Inequality~\eqref{EqnFlatness} combined with the bound $\sup_{z \in
  \real} |f(z)(1-f(z))| \leq 1$ ensures that $f$ is
\mbox{$\sqrt{\obupper{4\obsnr}}$-Lipschitz} over the interval
$[-4\obsnr,4\obsnr]$. It follows that
\begin{align}
\label{EqnOBThetaGradBound}
\fronorm{\PopGrad(\PlTheta)} & \le
\frac{\pobs\obupper{4\obsnr}}{\noisestd^{2}} \fronorm{\PlTheta -
  \ThetaStar} \le \frac{3\pobs\obupper{4\obsnr}}{\noisestd^{2}}
\fronorm{\d},
\end{align}
where the last step uses the upper bound $\opnorm{\d} \le
1$. Combining these bounds, we obtain the following lower bound on the
expectation:
\begin{align*}
\Term_{1} = \trinprod{\PopGrad(\PlTheta)}{\PlTheta - \ThetaStar +
  \OUTER{\d}} & \geq \trinprod{\PopGrad(\PlTheta)}{\PlTheta -
  \ThetaStar} - \fronorm{\PopGrad(\PlTheta)} \fronorm{\d}^{2} \\
 & \geq \frac{\pobs}{16\noisestd^{2}\oblower{4\obsnr}}
\fronorm{\d}^{2} - \frac{3\pobs\obupper{4\obsnr}}{\noisestd^{2}}
\opnorm{\d}\cdot\fronorm{\d}^{2}\\
& \ge \frac{\pobs}{32\noisestd^{2}\oblower{4\obsnr}} \fronorm{\d}^{2}
\end{align*}
where the last step uses the bound $\opnorm{\d} \le
\frac{1}{96\oblower{4\obsnr}\obupper{4\obsnr}}$. \\

\vtiny

\noindent \emph{Controlling $\Term_2$:} We now turn to analysis of the
deviation term $\Term_2$.  Using the symmetry of the matrix $\nabla
\EmpLoss(\PlTheta) - \PopGrad(\PlTheta)$, we may rewrite it as
\begin{align*}
\Term_{2} = 2\inprod{\nabla \EmpLoss(\PlTheta) -
  \PopGrad(\PlTheta)}{\d\otimes\pltheta}
\end{align*}
We control this quantity via two auxiliary lemmas.  For each $B \in
(0,1)$, define the annular set
\begin{align*}
\Gamma(B) & \defn \{ \d \,\mid\, \frac{B}{2} < \fronorm{\d} \le
B,\twoinfnorm{\d}^{2} \leq 4\obsnr\noisestd\}, \qquad \mbox{and the
  event} \\
\Event(B) & \defn \Big \{ \sup_{\d \in \Gamma(B)}
\frac{\abs{2\inprod{\nabla \EmpLoss(\OUTER{\pltheta}) -
      \PopGrad(\OUTER{\pltheta})}{\d\otimes\pltheta}}}{\fronorm{\d}} > \frac{1}{4}
\alpha \epsnum \Big \}.
\end{align*}
Our first lemma controls the probability of this ``bad event'':
\begin{lem}
\label{LemOBPeel}
For any $B \in (0,1)$, we have $\mprob[\Event(B)] \leq 2\usedim^{-14}$.
\end{lem}
Our next lemma gives controls over the small ball $\Gamma_{0}$ around
the origin.  In particular, we define the set
\begin{align*}
\Gamma_{0}& \defn \{\d \,\mid\, \fronorm{\d} \le 2^{ -
  \usedim^{2}},\twoinfnorm{\d}^{2} \le4\obsnr\noisestd\}, \;
\mbox{and} \; \Event_0 \defn \Big \{ \sup_{\d \in \Gamma_{0}}
\frac{\abs{2\inprod{\nabla \EmpLoss(\OUTER{\pltheta}) -
      \PopGrad(\OUTER{\pltheta})}{\d\otimes\pltheta}}}{\fronorm{\d}} >
\frac{1}{4} \alpha\epsnum \Big \}.
\end{align*}
\begin{lem}
\label{LemOBSmallRad}
We have $\mprob[\Event_0] \leq \usedim^{-12}$.
\end{lem}
\noindent See Appendices~\ref{SecProofLemOBPeel}
and~\ref{SecProofLemOBSmallRad} for the proofs of these two auxiliary
results.\\

Taking the lemmas as given, the union bound then guarantees that
\begin{align*}
\mprob \left[\sup_{\fronorm{\d} \leq 1,\twoinfnorm{\d}^{2} \leq 4 \obsnr \noisestd} \frac{\abs{2\inprod{\nabla \EmpLoss(\OUTER{\pltheta}) - \PopGrad(\OUTER{\pltheta})}{\d\otimes\pltheta}}}{\fronorm{\d}}
> \frac{1}{4} \alpha\epsnum\right] & \leq \mprob[\Event_0] + \sum_{i= 0}^{\usedim^{2}}\mprob[\Event(2^{-i})] \\
& \leq \usedim^{2}\cdot2\usedim^{-14} + \usedim^{-12} =
  3\usedim^{-12},
\end{align*}
which implies that $\abs{\Term_{2}} \le \frac{1}{4} \curv\epsnum\fronorm{\d}$ with
probability at least $1-3\usedim^{-12}$, Conditioned on this event, we
can combine the bounds for $\Term_{1}$ and $\Term_{2}$ to conclude
that
\begin{align*}
\trinprod{\nabla \EmpLoss(\PlTheta)}{\PlTheta - \ThetaStar + \OUTER{\d}} \ge \Term_{1} - \abs{\Term_{2}} 
\ge 2\curv\fronorm{\d}^{2} - \frac{1}{4} \curv\epsnum\fronorm{\d},
\end{align*}
valid for all matrices $\d$ such that $\fronorm{\d} \le \smallrad =
\max\{\frac{1}{16},\frac{1}{96\oblower{4\obsnr}\obupper{4\obsnr}}\}$,
thereby establishing the local descent condition~\eqref{EqnCurveCond} for $ \EmpLoss $.

\paragraph{Local Lipschitz condition~\eqref{EqnLipCon} and smoothness condition~\eqref{EqnSmoothCond2}:}

We begin by making note of the bounds
\begin{align*}
\twoinfnorm{\pltheta' - \pltheta''} \leq 4
\sqrt{\frac{\inco\rdim}{\usedim}}, \quad \mbox{and} \quad
\fronorm{\pltheta' - \pltheta''} \le 2,
\end{align*}
valid for all matrices $\pltheta',\pltheta'' \in \GenCons \cap
\Ballstar{\smallrad}$.	Moreover, we have $\curv\le \LipCon= \smoo$
and $\epsnum \leq 1$. Using these facts, it follows that the local Lipschitz and
smoothness conditions~\eqref{EqnLipCon} and~\eqref{EqnSmoothCond2} for $ \EmpLoss $ are implied a bound of the
form
\begin{align}
\label{EqnOBSmoothness}
 \fronorm{\nabla \EmpLoss(\PlTheta)\plphi} \le \frac{1}{8} \big(
 \smoo\fronorm{\pltheta - \thetastar} + \curv\epsnum\fronorm{\plphi}
 \big),
\end{align}
valid for all matrices $\pltheta \in \GenCons \cap
\Ballstar{\smallrad}$, and $\thetastar \in \eclass$, and matrices
$\plphi$ such that \mbox{$\twoinfnorm{\plphi} \le
  4\sqrt{\frac{\inco\rdim}{\usedim}}$.}  Accordingly, the remainder of
our effort is devoted to establishing the
bound~\eqref{EqnOBSmoothness}. We first condition on the event in
Lemma~\ref{LemMCProjSmall}, which occurs with probability at least
$1-2\usedim^{-4}$. We then make note of the decomposition
\begin{align*}
\fronorm{\nabla \EmpLoss(\PlTheta)\plphi} & = \frac{1}{\noisestd}
\fronorm{\ProjObs[h(\PlTheta/\noisestd)]\plphi} \; \leq
\underbrace{\frac{2}{\noisestd} \fronorm{\ProjObs
    \big[h(\PlTheta/\noisestd)-h(\ThetaStar/\noisestd)
      \big]\plphi}}_{\Term_1} + \underbrace{\frac{2}{\noisestd}
  \fronorm{\ProjObs[h(\ThetaStar/\noisestd)]\plphi}}_{\Term_2},
\end{align*}
and we bound each of these two terms separately.

\noindent \emph{Bounding term $\Term_1$:} We have
\begin{align*}
\Term_{1} & = \frac{2}{\noisestd}\sup_{\plomega \in
  \real^{\usedim\times \rdim}:\fronorm{\plomega} =
  1}\abs{\trinprod{\ProjObs
    \big[H(\PlTheta/\noisestd)-H(\ThetaStar/\noisestd)
      \big]}{\plomega\otimes\plphi}} \\
& \leq \frac{2}{\noisestd} \fronorm{\ProjObs \big(
  H(\PlTheta/\noisestd)-H(\ThetaStar/\noisestd)
  \big)}\cdot\sup_{\plomega \in \real^{\usedim\times
    \rdim}:\fronorm{\plomega} = 1}
\fronorm{\ProjObs(\plomega\otimes\plphi)}.
\end{align*}
Recall that each component $[H(\cdot)]_{ij}$ is almostly surely
$4\obupper{4\obsnr}$-Lipschitz over the interval $[-4\obsnr,4\obsnr]$,
and that $\infnorm{\PlTheta/\noisestd},\infnorm{\ThetaStar/\noisestd}
\leq 4 \obsnr$.  Combining these facts yields the bound
\begin{align*}
\fronorm{\ProjObs \big( H(\PlTheta/\noisestd)-H(\ThetaStar/\noisestd)
  \big)} & \leq \frac{4\obupper{4\obsnr}}{\noisestd} \fronorm{\ProjObs
  \big( \PlTheta - \ThetaStar \big)} \\
& \leq \frac{4\obupper{4\obsnr}}{\noisestd} \Big( 2\fronorm{\ProjObs
  \big( \thetastar\otimes(\pltheta - \thetastar) \big)} +
\fronorm{\ProjObs \big( \OUTER{(\pltheta - \thetastar)} \big)} \Big)
\\
& \leq \frac{72\obupper{4\obsnr}\sqrt{2\pobs\inco\rdim}}{\noisestd}
\fronorm{\d},
\end{align*}
where the last inequality follows from the
inequality~(\ref{EqnMCProjSmall3}) in Lemma~\ref{LemMCProjSmall}
combined with the fact that $\max \big \{
\twoinfnorm{\thetastar},\twoinfnorm{\pltheta - \thetastar} \big \}
\leq 4 \sqrt{\frac{\inco\rdim}{\usedim}}$.  Applying
inequality~\eqref{EqnMCProjSmall3} a second time yields
\begin{align*}
\sup_{\plomega \in \real^{\usedim\times \rdim}:\fronorm{\plomega} = 1}
\fronorm{\ProjObs(\plomega\otimes\plphi)}
\le6\sqrt{2\pobs\inco\rdim}\sup_{\plomega \in \real^{\usedim\times
    \rdim}:\fronorm{\plomega} = 1} \fronorm{\plomega} =
6\sqrt{2\pobs\inco\rdim}.
\end{align*}
Putting together the pieces yields the upper bound
\begin{subequations}
\begin{align}
\label{EqnAppleOne}
\Term_{1} \le \frac{48\obupper{4\obsnr}\pobs\inco\rdim}{\noisestd^{2}}
\fronorm{\pltheta - \thetastar} \le \frac{\smoo}{8} \fronorm{\pltheta
  - \thetastar}.
\end{align}

\vtiny

\noindent \emph{Bounding term $\Term_2$:} Turning to the term
$\Term_{2}$, we observe that conditioned on $\Yout$, the matrix
$H(\ThetaStar/\noisestd)$ is a deterministic quantity with entries
bounded uniformly as $\abs{[H(\ThetaStar/\noisestd)]_{oj}} \le
\obupper{4\obsnr}$ for all indices $i,j \in [\usedim]$.  By applying
Lemma~\ref{LemCensorGaussMat} and integrating out the conditioning, we
find that the second term is bounded as
\begin{align}
\label{EqnAppleTwo}
\Term_{2} & \le \frac{2}{\noisestd}
\opnorm{\ProjObs[h(\ThetaStar/\noisestd)]}\cdot\sqrt{\rdim}
\fronorm{\plphi} \leq \frac{4 \sqrt{\pobs\rdim\usedim}
  \obupper{4\obsnr}}{\noisestd}\cdot\fronorm{\plphi} \leq
\frac{1}{8}\curv\epsnum\fronorm{\plphi},
\end{align}
\end{subequations}
where these bounds hold with probability at least $1 - \usedim^{-4}$.\\

\vtiny

\noindent Finally, combining our bounds for $\Term_{1}$ and
$\Term_{2}$ in equations~\eqref{EqnAppleOne} and~\eqref{EqnAppleTwo}
yields the desired bound~\eqref{EqnOBSmoothness}.

\paragraph{Local smoothness condition~(\ref{EqnSmoothCond1}):}

We begin by conditioning on the event in Lemma~\ref{LemMCProjSmall},
which holds with probability at least $1-2\usedim^{-4}$. For each pair
of matrices $\pltheta',\pltheta''$ in the set $\GenCons \cap
\Ballstar{\smallrad}$ and any matrix $\thetastar \in \eclass$, we then
have
\begin{align*}
\abs{\trinprod{\nabla \EmpLoss(\PlTheta) - \nabla
    \EmpLoss(\PlTheta')}{\pltheta'\otimes(\pltheta - \thetastar)}} & =
\frac{1}{\noisestd}\abs{\inprod{\ProjObs \big(
    H(\PlTheta/\noisestd)-h(\PlTheta'/\noisestd)
    \big)}{\pltheta'\otimes(\pltheta - \thetastar)}} \\
& \leq \frac{1}{\noisestd} \fronorm{\ProjObs \big(
  H(\PlTheta/\noisestd)-h(\PlTheta'/\noisestd) \big)}
\fronorm{\ProjObs(\pltheta'\otimes(\pltheta - \thetastar))} \\
& \leq \frac{4\obupper{4\obsnr}}{\noisestd^{2}} \fronorm{\ProjObs
  \big( \PlTheta - \PlTheta' \big)}
\fronorm{\ProjObs(\pltheta'\otimes(\pltheta - \thetastar))},
\end{align*}
where the last inequality follows from the fact that the function $H$
is element-wise $4\obupper{4\obsnr}$-Lipschitz.
Recall equation~\eqref{EqnMCObsLambda} from
Section~\ref{SecProofCorMC}, which ensures that
\begin{align*}
\fronorm{\ProjObs \big( \PlTheta - \PlTheta' \big)} &
\le14\sqrt{2\pobs\inco\rdim} \fronorm{\pltheta - \pltheta'}.
\end{align*}
Combined with inequality~\eqref{EqnMCProjSmall3} from
Lemma~\ref{LemMCProjSmall}, we find that
$\fronorm{\ProjObs(\pltheta'\otimes(\pltheta - \thetastar))} \le 6
\sqrt{2\pobs\inco\rdim} \fronorm{\pltheta - \thetastar}$, from which
it follows that
\begin{align*}
\abs{\trinprod{\nabla \EmpLoss(\PlTheta) - \nabla
    \EmpLoss(\PlTheta')}{\pltheta' \otimes (\pltheta - \thetastar)}} &
\leq \frac{672 \obupper{4 \obsnr}}{\noisestd^{2}}\cdot \pobs \inco
\rdim\fronorm{\pltheta - \pltheta'} \fronorm{\pltheta - \thetastar} \\
& \leq \frac{1}{4} \smoo\fronorm{\pltheta - \pltheta'} \fronorm{\pltheta -
  \thetastar},
\end{align*}
thereby establishing the local smoothness
condition~\eqref{EqnSmoothCond1}.


\subsection{Proof of Corollary~\ref{CorLS}}
\label{SecProofCorLS}

We now prove our claims for the matrix decomposition problem.  By
dividing through $\opnorm{\thetastar}$, we may assume without loss of
generality that $\opnorm{\thetastar} = 1$.  The set $\GenCons$ and the
values of $\smallrad$ and $\dist{\thetait{0}}$ are the same as used in
the proof of matrix completion in Section~\ref{SecProofCorMC}, so we
make use of the results therein.  In particular, we showed there that
the set $\GenCons$ is $\ThetaStar$-faithful.

Given the observation matrix
$\Yout= \ThetaStar + \SStar + \noisevar$, the
gradient takes the form
\begin{align*}
\nabla \EmpLoss(\PlTheta)= (\PlTheta - \ThetaStar) + (\PlS(\PlTheta) -
\SStar) - \noisevar,
\end{align*}
where $\PlS(\PlTheta) \defn \ProjConsS(\Yout - \PlTheta).$ Below we
verify the local descent, Lipschitz and smoothness conditions.

\paragraph{local descent:}

Expanding $\nabla \EmpLoss(\PlTheta)$, the quantity
$\trinprod{\nabla \EmpLoss(\PlTheta)}{\PlTheta - \ThetaStar + \OUTER{\d}}$
can be decomposed into the sum
\begin{align*}
\underbrace{\fronorm{\PlTheta - \ThetaStar}^{2}}_{\Term_1} +
\underbrace{\trinprod{\PlTheta - \ThetaStar}{\OUTER{\d}}}_{\Term_2} +
\underbrace{\trinprod{\PlS(\PlTheta) - \SStar}{\PlTheta - \ThetaStar +
    \OUTER{\d}}}_{\Term_3} + \underbrace{\trinprod{ -
    \noisevar}{\PlTheta - \ThetaStar + \OUTER{\d}}}_{\Term_4}.
\end{align*}
Note that $\PlTheta - \ThetaStar= \tdthetastar\otimes\d +
\d\otimes\tdthetastar + \OUTER{\d}.$ By Lemma~\ref{LemThetaStar}, the
matrix $\d^{\top}\tdthetastar$ is symmetric, so expanding the
Frobenius norm shows that $\Term_{1} \ge2\fronorm{\d}^{2}.$ Since
$\opnorm{\d} \le \frac{3}{5}$, we have
\[
\abs{\Term_{2}} \le 2\opnorm{\d} \fronorm{\d}^{2} + \opnorm{\d}^{2} \fronorm{\d}^{2} \le \frac{39}{25} \fronorm{\d}^{2}.
\]
 With $\dS \defn \PlS(\PlTheta) - \SStar$ and $e_{j}$ being the $j$-th
standard basis, we find that
\begin{align*}
\abs{\Term_{3}} = \abs{2\trinprod{\dS}{\d\otimes\pltheta}} \le & 2\fronorm{\pltheta^{\top}\dS} \fronorm{\d} \\
=  & 2\sqrt{\sum_{j= 1}^{\usedim}\twonorm{\pltheta^{\top}\dS e_{j}}^{2}} \fronorm{\d} \le 2\sqrt{\sum_{j= 1}^{\usedim}\twoinfnorm{\pltheta}^{2}\onenorm{\dS e_{j}}^{2}} \fronorm{\d},
\end{align*}
where we use the symmetry of $\dS$ in the first equality.
Inequality~\eqref{EqDtwoinf} ensures that \mbox{$\twoinfnorm{\pltheta}
  \le 2\sqrt{\frac{\inco\rdim}{\usedim}}$.} Moreover, for each $\PlS
\in \ConsS$, we have the inequalities $\onenorm{\dS e_{j}} \le
2\sqrt{k}\twonorm{\dS e_{j}},j\in[\usedim]$ thanks to the row-wise
$\ell_{1}$ constraints and the $\kdim$-sparsity of the columns of
$\SStar$. It follows that
\begin{align*}
\abs{\Term_{3}} \le &
4\sqrt{\frac{\inco\rdim\kdim}{\usedim}}\sqrt{\sum_{j=
    1}^{\usedim}\twonorm{\dS e_{j}}^{2}} \fronorm{\d} =
4\sqrt{\frac{\inco\rdim\kdim}{\usedim}} \fronorm{\dS} \fronorm{\d}.
\end{align*}
Under the assumption $\frac{\inco\rdim\kdim}{\usedim} \le c_{1}$ of
the corollary, we obtain $\abs{\Term_{3}} \le \frac{2}{25}
\fronorm{\dS} \fronorm{\d}.$ But $\SStar \in \ConsS$ and the
projection $\ProjConsS$ is non-expansive, whence
\begin{equation}
\label{EqnDS}
\fronorm{\dS} = \fronorm{\ProjConsS(\Yout - \PlTheta) - \SStar} \le
\fronorm{(\Yout - \PlTheta) - \SStar} = \fronorm{\PlTheta -
  \ThetaStar} \le3\fronorm{\d},
\end{equation}
and so we have shown that $\abs{\Term_{3}} \le \frac{6}{25}
\fronorm{\d}^{2}$.  Finally, we have
\begin{align*}
\abs{\Term_{4}} \le \fronorm{\noisevar} \fronorm{\PlTheta - \ThetaStar + \OUTER{\d}} \le \frac{16}{5} \fronorm{\noisevar} \fronorm{\d}.
\end{align*}
Putting together the bounds for $\Term_{1}$, $\Term_2$ and $\Term_3$
and $\Term_{4}$, we conclude that
\begin{align*}
\trinprod{\nabla \EmpLoss(\PlTheta)}{\PlTheta - \ThetaStar +
  \OUTER{\d}} \ge \frac{1}{5} \fronorm{\d}^{2} - \frac{16}{5}
\fronorm{\noisevar} \fronorm{\d},
\end{align*}
thereby proving the local descent condition~\eqref{EqnCurveCond} for $ \EmpLoss $.


\paragraph{Local Lipschitz condition:}

Using the inequality~(\ref{EqnDS}) above and the assumption
$\fronorm{\noisevar} \le \frac{2}{5}$, we have
\begin{align*}
\fronorm{\nabla \EmpLoss(\PlTheta)\pltheta} & = \fronorm{(\PlTheta -
  \ThetaStar + \PlS(\PlTheta) - \SStar - \noisevar)\pltheta} \\
& \leq \big( \fronorm{\PlTheta - \ThetaStar} + \fronorm{\dS} +
\fronorm{\noisevar} \big) \opnorm{\pltheta} \\
& \leq (6\fronorm{\d} + 1) \opnorm{\pltheta} \\
& \le 8,
\end{align*}
where the last inequality follows from $\opnorm{\d} \le \fronorm{\d}
\le \frac{3}{5}$.  Therefore, $ \EmpLoss $ satisfies the local
Lipschitz condition in~\eqref{EqnLipCon}.

\paragraph{Local smoothness:}

Observe that
\begin{align*}
\abs{\trinprod{\nabla \EmpLoss(\PlTheta) - \nabla
    \EmpLoss(\PlTheta')}{\pltheta'\otimes(\pltheta - \thetastar)}} & =
\abs{\trinprod{\PlTheta - \PlTheta' + \PlS(\PlTheta) -
    \PlS(\PlTheta')}{\pltheta'\otimes(\pltheta - \thetastar)}} \\
&
\le \big( \fronorm{\PlTheta - \PlTheta'} + \fronorm{\PlS(\PlTheta) -
  \PlS(\PlTheta')} \big) \fronorm{\pltheta - \thetastar}
\opnorm{\pltheta'}.
\end{align*}
The non-expansiveness of the projection $\ProjConsS$ ensures that
\begin{align*}
\fronorm{\PlS(\PlTheta) - \PlS(\PlTheta')} = \fronorm{\ProjConsS(\Yout
  - \PlTheta) - \ProjConsS(\Yout - \PlTheta')} \le \fronorm{\PlTheta -
  \PlTheta'} \le \frac{16}{5} \fronorm{\pltheta - \pltheta'},.
\end{align*}
where we use $\pltheta,\pltheta' \in \Ballstar{\frac{3}{5}}$. It
follows that
\begin{align*}
\abs{\trinprod{\nabla \EmpLoss(\PlTheta) - \nabla
    \EmpLoss(\PlTheta')}{\pltheta'\otimes(\pltheta - \thetastar)}} \le
12 \fronorm{\pltheta - \pltheta'} \fronorm{\pltheta - \thetastar},
\end{align*}
proving the first smoothness condition~\eqref{EqnSmoothCond1}.

Similarly, combining inequality~\eqref{EqnDS} with the bound
$\fronorm{\noisevar} \leq \frac{2}{5}$ implies that
\begin{align*}
\abs{\trinprod{\nabla \EmpLoss(\PlTheta)}{(\pltheta -
    \thetastar)\otimes(\pltheta' - \pltheta'')}} & =
\abs{\trinprod{\PlTheta - \ThetaStar + \PlS(\PlTheta) - \SStar -
    \noisevar}{(\pltheta - \thetastar)\otimes(\pltheta' -
    \pltheta'')}} \\
& \le \big( \fronorm{\PlTheta - \ThetaStar} + \fronorm{\PlS(\PlTheta)
  - \SStar} + \fronorm{\noisevar} \big) \fronorm{\pltheta -
  \thetastar} \fronorm{\pltheta' - \pltheta''} \\
 & \leq 7 \fronorm{\pltheta - \thetastar} \fronorm{\pltheta' -
  \pltheta''},
\end{align*}
thereby verifying the second smoothness
condition~\eqref{EqnSmoothCond2}.


\section{Discussion}
\label{SecConclusion}

In this paper, we have laid out a general framework for analyzing the
behavior of projected gradient descent for solving low-rank
optimization problems in the factorized space.  We have illustrated
the consequences of our general theory for a number of concrete
models, including matrix regression, structured PCA, matrix
completion, matrix decomposition and graph clustering.


\subsection*{Acknowledgements}
This work was partially supported by ONR-MURI grant DOD-002888, AFOSR
grant FA9550-14-1-0016, NSF grant CIF-31712-23800, and ONR MURI grant
N00014-11-1-0688.



\appendix



\section{Proof of Theorem~\ref{ThmGeneralL}}
\label{SecProofThmGeneralL}

Recall that in Section~\ref{SecProofThmGeneralSmooth} we proved
Theorem~\ref{ThmGeneralSmooth} under the assumption that $ \LossTil $
satisfies the \emph{relaxed} local Lipschitz
condition~(\ref{EqnLtilWeakLipCon}) as well as the local descent
condition~(\ref{EqnLtilCurveCond}) and smoothness
conditions~(\ref{EqnLtilSmoothCond}). We establish
Theorem~\ref{ThmGeneralL} by showing that these conditions for $
\LossTil $ are implied by the corresponding conditions for $ \EmpLoss
$ in Definitions~\ref{DefCurveCon}--\ref{DefSmoothCon} with the same
parameters $ \curv$, $\smoo$, $\LipCon$, $\epsnum $ and $ \smallrad$
with $\smallrad < \sigma_\rdim (\thetastar) $. \\

Let $ \pltheta$ be an arbitrary matrix in $ \GenCons\cap\Ballstar{\smallrad} $ and $ \thetastar $ an arbitrary member in $ \eclass $. Since $ \smallrad< \sigma_\rdim(\thetastar)$, Lemma~\ref{LemThetaStar} guarantees that $ \tdthetastar = \arg\min_{A \in \eclass} \fronorm{A - \pltheta}$ is uniquely defined. We use the shorthand $ \Grad \defn \nabla_{\PlTheta} \EmpLoss $, $ \GradTil \defn \nabla_{\pltheta} \LossTil $, $ \PISTAR{\d} \defn \pltheta - \tdthetastar $ and $ \d = \pltheta - \thetastar $.

\paragraph{Local descent condition:}

Observe that 
\begin{align*}
\trinprod{\GradTil (\pltheta)}{\pltheta - \thetastar}
=  & \trinprod{ \Grad (\OUTER{\pltheta}) }{\d \otimes\pltheta + \pltheta\otimes\d} \\
=  & \trinprod{ \Grad (\OUTER{\pltheta}) }{\OUTER{\pltheta} - \OUTER{\thetastar}} + \trinprod{ \Grad (\OUTER{\pltheta}) }{\d \otimes\d } \\
=  & \trinprod{ \Grad (\OUTER{\pltheta}) }{ \OUTER{\pltheta} - \OUTER{\tdthetastar} } + \trinprod{ \Grad (\OUTER{\pltheta}) }{\PISTAR{\d} \otimes\PISTAR{\d} }  \\
&   + \trinprod{ \Grad (\OUTER{\pltheta}) }{(\tdthetastar-\thetastar)\otimes\d}  +  \trinprod{ \Grad (\OUTER{\pltheta}) }{\PISTAR{\d} \otimes(\tdthetastar-\thetastar)},
\end{align*}
where the last step follows from $ \OUTER{\thetastar} = \OUTER{\tdthetastar} $ and $ \d = \PISTAR{\d} +\tdthetastar - \thetastar $.
We then apply the local descent condition~(\ref{EqnCurveCond}) for $ \EmpLoss $ to the first two terms above, and the local smoothness condition~(\ref{EqnSmoothCond2}) for $ \EmpLoss $ to the last two terms. Doing so yields
\begin{align*}
\trinprod{\GradTil (\pltheta)}{\pltheta - \thetastar} 
\ge & 2\curv\fronorm{\PISTAR{\d}}^{2} - \frac{\smoo}{4}\fronorm{\tdthetastar-\thetastar } \big( \fronorm{\d} + \fronorm{\PISTAR{\d}} \big) - \frac{\curv}{2}\epsnum\fronorm{\PISTAR{\d}}  - \frac{\curv}{2}\epsnum\fronorm{\d}\\
\ge & 2\curv\fronorm{\PISTAR{\d}}^{2} - \frac{\smoo}{2}\fronorm{\tdthetastar-\thetastar } \fronorm{\PISTAR{\d}}  - \frac{\smoo}{2}\fronorm{\tdthetastar-\thetastar }^2 
- \curv\epsnum\fronorm{\PISTAR{\d}} 
- \frac{\curv}{2}\epsnum\fronorm{\tdthetastar-\thetastar }\\
\ge & 2\curv\fronorm{\PISTAR{\d}}^{2} - \Big( \frac{1}{2}\curv \fronorm{\PISTAR{\d}}^2 + \frac{\smoo^2}{4\curv}\fronorm{\tdthetastar-\thetastar }^2  \Big)
- \frac{\smoo}{2}\fronorm{\tdthetastar-\thetastar }^2 \\
&
- \Big( \frac{1}{2}\curv \fronorm{\PISTAR{\d}}^2  + \frac{1}{2}\curv\epsnum^2 \Big) 
- \Big( \frac{1}{2}\curv\epsnum^2 + \frac{\curv}{4}  \fronorm{\tdthetastar-\thetastar }^2 \Big)\\
\ge  &  \curv\fronorm{\PISTAR{\d}}^{2} -  \frac{\smoo^2}{\curv} \fronorm{\tdthetastar-\thetastar } - \curv\epsnum^2, 
\end{align*}
where the last three steps follow from $ \fronorm{\d} \le \fronorm{\PISTAR{\d}} + \fronorm{\tdthetastar - \thetastar } $, the AM-GM inequality and the upper bound $ \curv \le \smoo $, respectively. This proves the local descent condition for~$ \LossTil $ in (\ref{EqnLtilCurveCond}).

\paragraph{Relaxed local Lipschitz condition:}

Let $ \pltheta' $ be an arbitrary matrix in $ \GenCons $.  With  $\Grad(\OUTER{\pltheta})$ symmetric and $ \EmpLoss $ satisfying the relaxed Lipschitz condition in~(\ref{EqnWeakLipCon}), we have
\begin{align*}
\abs{ \trinprod{\GradTil(\pltheta) }{ \pltheta - \pltheta' } }= \abs{ \trinprod{2\Grad(\OUTER{\pltheta})}{ (\pltheta - \pltheta') \otimes \pltheta} } 
\le \LipCon \big( \opnorm{\thetastar}^2 + \fronorm{\thetastar} \fronorm{\pltheta - \pltheta'}\big),
\end{align*}
which proves the relaxed local Lipschitz condition for $ \LossTil $ in~(\ref{EqnLtilWeakLipCon}).
 
\paragraph{Local smoothness condition:} 
Let $ \pltheta' $ be an arbitrary matrix in $ \GenCons \cap \Ballstar{\smallrad}$. The smoothness conditions~(\ref{EqnSmoothCond1}) and (\ref{EqnSmoothCond2}) yields that
\begin{align*}
\abs{ \trinprod{\GradTil(\pltheta) - \GradTil(\pltheta')}{\pltheta-\thetastar} }
& = \abs{ 2\trinprod{\Grad(\OUTER{\pltheta}) \pltheta - \Grad(\OUTER{\pltheta'}) \pltheta' }{\pltheta-\thetastar} }\\
& = \abs{ 2\trinprod{\Grad(\OUTER{\pltheta}) \pltheta' - \Grad(\OUTER{\pltheta'})\pltheta'}{\pltheta-\thetastar} + 2\trinprod{\Grad(\OUTER{\pltheta})(\pltheta - \pltheta')}{\pltheta-\thetastar} } \\
 & \leq \smoo\fronorm{\pltheta-\pltheta'} \fronorm{\pltheta-\thetastar} + \curv\epsnum\fronorm{\pltheta-\thetastar},
\end{align*}
which establishes the smoothness condition for $ \LossTil $ in~(\ref{EqnLtilSmoothCond}). \\


\section{Technical lemmas for Corollary~\ref{CorMatrixSensing}}

In this appendix, we prove the technical lemmas involved in the proof
of Corollary~\ref{CorMatrixSensing} on the matrix sensing model.

\subsection{Proof of Lemma~\ref{LemMatSenInprod}}
\label{SecProofLemMatSenInprod}

By the bilinearity of the inner product, we may assume without loss of
generality that \mbox{$\fronorm {A}= \fronorm {B}= 1$.}  Since the
matrices $A \pm B$ have rank at most~$4\rdim$, the RIP with
$\ripparam{4\rdim}$ ensures that
\begin{align*}
(1 - \ripparam{4\rdim}) \fronorm{A \pm B}^{2} \le \frac{1}{\numobs}
  \fronorm{\Xmap (A\pm B)}^{2} \le (1 + \ripparam{4\rdim}) \fronorm{A
    \pm B}^{2}.
\end{align*}
It follows that
\begin{align*}
\frac{1}{\numobs} \trinprod{\Xmap (A)} {\Xmap (B )} & =
\frac{1}{4\numobs} \Big( \fronorm{\Xmap(A + B)}^{2} -
\fronorm{\Xmap(A-B)}^{2} \Big) \\ \le & \frac{1}{4} \Big( (1 +
\ripparam{4\rdim}) \fronorm{A + B}^{2} - (1 - \ripparam{4\rdim})
\fronorm{A-B}^{2} \Big) \\
& = \trinprod{A}{B} + \frac{1}{2} \ripparam{4\rdim} \big( \fronorm
A^{2} + \fronorm B^{2} \big) \\
& = \trinprod AB + \ripparam{4\rdim} \fronorm{A} \fronorm{B}.
\end{align*}
It follows from a similar argument that $\frac{1}{\numobs} \trinprod{\Xmap (A)} {\Xmap (B)} \ge \trinprod{A}{B} - \ripparam{4\rdim} \fronorm{A} \fronorm{B}.$


\section{Technical lemmas for Corollary~\ref{CorMC}}

In this appendix, we prove the technical lemmas involved in the
proof of Corollary~\ref{CorMC} on the matrix completion model.

\subsection{Proof of Lemma~\ref{LemMCInprod}}
\label{SecProofLemMCInProd}

Define a  subspace $\Tspace \subseteq \real^{\usedim \times \usedim}$
of $\usedim$-dimensional matrices as follows
\begin{align*}
\Tspace & \defn \Big \{X \mid X = (\thetastar\otimes\plu) + ( \plv
\otimes \thetastar) \, \quad \mbox{for some $\plu,\plv \in
  \real^{\usedim\times \rdim}$} \Big \},
\end{align*}
and let $\ProjTspace$ be the Euclidean projection onto $\Tspace$.
Since $\thetastar$ is $4\inco$-incoherent, a known result in exact
matrix completion~\cite{candes2009exact} guarantees that as long as
$\pobs\ge c \frac{\inco\rdim\log\usedim}{\epsilon^{2}\usedim}$ for a
suficiently large universal constant $c$, then
\begin{align*}
\fronorm{(\ProjTspace\ProjObs\ProjTspace - \pobs\ProjTspace)X} \le
\epsilon\pobs\fronorm X, \quad \mbox{for all \ensuremath{X
    \in \Tspace}}
\end{align*}
with probability at least $1-2\usedim^{-3}$.  Noting that the matrices
$\thetastar \otimes \plphi \pm \plomega \otimes \thetastar$ belong to
the subspace $\Tspace$, we can apply the above inequality to obtain
\begin{align*}
(1 - \epsilon) \pobs \fronorm{\thetastar\otimes\plphi \pm
    \plomega\otimes\thetastar} \le
  \fronorm{\ProjTspace\ProjObs\ProjTspace \big(
    \thetastar\otimes\plphi \pm \plomega\otimes\thetastar \big)} \le
  (1 + \epsilon) \pobs \fronorm{\thetastar\otimes\plphi \pm
    \plomega\otimes\thetastar}.
\end{align*}
The rest of the proof is similar to that of
Lemma~\ref{LemMatSenInprod}.  In particular, by the bilinearity of the
inner product, we may assume $\fronorm{\thetastar\otimes\plphi} =
\fronorm{\plomega\otimes\thetastar} = 1$. Using the above
inequalities, we find that
\begin{align*}
&
  \trinprod{\ProjObs\left(\thetastar\otimes\plphi\right)}{\ProjObs\left(\plomega\otimes\thetastar\right)}
  =
  \trinprod{\ProjObs\ProjTspace(\thetastar\otimes\plphi)}{\ProjObs\ProjTspace(\plomega\otimes\thetastar)}
  \\
= & \frac{1}{4}
\left(\fronorm{\ProjObs\ProjTspace(\thetastar\otimes\plphi +
  \plomega\otimes\thetastar)}^{2} -
\fronorm{\ProjObs\ProjTspace(\thetastar\otimes\plphi -
  \plomega\otimes\thetastar)}^{2}\right) \\
\leq & \frac{1}{4} \big( \fronorm{\thetastar\otimes\plphi +
  \plomega\otimes\thetastar}
\fronorm{\ProjTspace\ProjObs\ProjTspace(\thetastar\otimes\plphi +
  \plomega\otimes\thetastar)} - \fronorm{\thetastar\otimes\plphi -
  \plomega\otimes\thetastar}
\fronorm{\ProjTspace\ProjObs\ProjTspace(\thetastar\otimes\plphi +
  \plomega\otimes\thetastar)} \big) \\
\leq & \frac{1}{4} \Big( (1 + \epsilon) \pobs
\fronorm{\thetastar\otimes\plphi + \plomega\otimes\thetastar}^{2} - (1
- \epsilon) \pobs \fronorm{\thetastar\otimes\plphi -
  \plomega\otimes\thetastar}^{2} \Big) \\
= & \pobs
\trinprod{\thetastar\otimes\plphi}{\plomega\otimes\thetastar} +
\epsilon\pobs = \pobs
\trinprod{\thetastar\otimes\plphi}{\plomega\otimes\thetastar} +
\epsilon \pobs \fronorm{\thetastar\otimes\plphi}
\fronorm{\plomega\otimes\thetastar}.
\end{align*}
This proves the first inequality in the lemma. The second inequality
can be proved in the same fashion by noting that the matrices
$(\thetastar \otimes \plphi) \pm (\thetastar\otimes\plomega)$ also
belong to the subspace $\Tspace$.


\subsection{Proof of Lemma~\ref{LemMCProjSmall}}
\label{SecProofLemMCProjSmall}

We need the following result on random
graphs~\cite{feige2005spectral}, which involves some universal
constants $c_1$ and $c_2$.
\begin{lem}
\label{LemRandomGraph}
If $\pobs \geq c_{1} \frac{\log\usedim}{\epsilon^{2}\usedim}$, then
with probability at least $1 - \frac{1}{2}d^{-4}$,
\begin{align}
\label{EqnRandomGraph}
\sum_{(i,j) \in \Obs}\plu_{i}\plv_{j} \le (1 +
\epsilon)\pobs\onenorm{\plu}\onenorm{\plv} +
c_{2}\sqrt{\pobs\usedim}\twonorm{\plu}\twonorm{\plv}, \quad
\forall\plu,\plv \in \real^{d}.
\end{align}
\end{lem}
From the bound~\eqref{EqnRandomGraph} and the assumption $\pobs \ge
\frac{C}{\epsilon^{2}} \big( \frac{\inco^{2}\rdim^{2}}{\usedim} +
\frac{\log\usedim}{\usedim} \big)$, we find that with probability at
least~$1 - \frac{1}{2}d^{-4}$,
\begin{align*}
\pobs^{-1} \fronorm{\ProjObs(\plphi\otimes\plphi)}^{2} & \le
\pobs^{-1} \sum_{(i,j) \in \Obs}\twonorm{\plphi_{i\cdot}}^{2}
\twonorm{\plphi_{j\cdot}}^{2} \leq (1 + \epsilon) \Big( \sum_{i}
\twonorm{\plphi_{i\cdot}}^{2} \Big)^{2} +
C_{2}\sqrt{\frac{\usedim}{\pobs}}\sum_{i}
\twonorm{\plphi_{i\cdot}}^{4} \\
& = (1 + \epsilon) \fronorm{\plphi}^{4}
+ C_{2}\sqrt{\frac{\usedim}{\pobs}} \fronorm{\plphi}^{2}
\twoinfnorm{\plphi}^{2}\\
& \leq (1 + \epsilon) \fronorm{\plphi}^{4} +
\epsilon\fronorm{\plphi}^{2},
\end{align*}
where the last step follows from the bound $\twoinfnorm{\plphi} \le
6\sqrt{\frac{\inco\rdim}{\usedim}}$.  We have thus established the
first inequality~\eqref{EqnMCProjSmall1} in the lemma statement.

Let $ \Obs_i \defn \{ j \,\mid\, (i,j) \in \Obs \} $. When $\pobs\ge C \frac{\log\usedim}{\usedim}$ for $C$ sufficiently
large, the event $\max_{i}\abs{\Obs_{i}} \le 2\pobs\usedim$ holds with
probability at least~$1 - \usedim^{-4}$. On this event, we have
\begin{align*}
\pobs^{-1} \fronorm{\ProjObs(Z)\plphi}^{2} & = \pobs^{-1}\sum_{i= 1}^{\usedim}\sum_{k= 1}^{\rdim} \Big( \sum_{j \in \Obs_{i}}(\ProjObs(Z^{\top})e_{i})_{j}\cdot\plphi_{jk} \Big)^{2} \\
 & \le \pobs^{-1}\sum_{i= 1}^{\usedim}\sum_{k= 1}^{\rdim}\twonorm{\ProjObs(Z^{\top})e_{i}}^{2}\sum_{j \in \Obs_{i}}\plphi_{jk}^{2} \\
 & = \pobs^{-1} \sum_{i= 1}^{\usedim} \twonorm{\ProjObs(Z^{\top})e_{i}}^{2} \sum_{j \in \Obs_{i}} \twonorm{\plphi_{j\cdot}}^{2} \\
 & \le \pobs^{-1} \sum_{i= 1}^{\usedim} \twonorm{\ProjObs(Z^{\top})e_{i}}^{2} \cdot \big( \max_{i}\abs{\Obs_{i}} \big)\twoinfnorm{\plphi}^{2}
  \le \fronorm{\ProjObs(Z)}^{2} \cdot 2 \usedim \twoinfnorm{\plphi}^{2}.
\end{align*}
But $\twoinfnorm{\plphi} \le6\sqrt{\frac{\inco\rdim}{\usedim}}$ by
assumption, so the second inequality~(\ref{EqnMCProjSmall2}) in
the lemma follows.

To establish the third inequality~\eqref{EqnMCProjSmall3} in the
lemma, observe that conditioned on the event $\{
\max_{i}\abs{\Obs_{i}} \le 2\pobs\usedim \}$, we have
\begin{align*}
\pobs^{-1} \fronorm{\ProjObs(\plphi\otimes\plomega)}^{2} \; = \;
\pobs^{-1}\sum_{j=1}^{\usedim}
\twonorm{\ProjObs(\plphi\otimes\plomega)e_{j}}^{2} & \leq
\pobs^{-1}\sum_{j= 1}^{\usedim}
\twonorm{\plomega_{j\cdot}}^{2}\sum_{\{ i \, \mid \, (i,j) \in \Omega
  \}} \twonorm{\plphi_{i\cdot}}^{2} \\
& \le \pobs^{-1}\sum_{j= 1}^{\usedim}\twonorm{\plomega_{j\cdot}}^{2}
\big( \max_{i}\abs{\Obs_{i}} \big)\twoinfnorm{\plphi}^{2} \\
& \le \pobs^{-1}\sum_{i=
  1}^{\usedim}\twonorm{\plomega_{j\cdot}}^{2}\cdot2\pobs\usedim\cdot36\frac{\inco\rdim}{\usedim}\\
& = 72 \inco\rdim\fronorm{\plomega}^{2}.
\end{align*}


\section{Proof of Lemma~\ref{LemREofW}}

\label{SecProofLemREofW}

By rescaling, it suffices to consider matrices $\plu$ and $\plv$ with
$\fronorm{\plu} = \fronorm{\plv} = 1$ and $\plu,\plv \in
\Balltwoone(\sqrt{\kdim})$. We need the following geometric result,
which is a simple generalization of Lemma~11 in the
paper~\cite{loh2012nonconvex}.	For completeness, we provide the proof
in Section~\ref{SecProofLemGEometric} to follow.
\begin{lem}
\label{LemGeometric}
For each integer $1 \leq \kdim \le \usedim$, we have
\begin{align}
\label{EqnGeometric}
\Balltwoone(\sqrt{\kdim}) \cap \Ballfro(1) \subseteq 3 \, \textup{cl}
\big\{\textup{conv}\{\Balltwozero(\kdim) \cap \Ballfro(1)\} \big\}.
\end{align}
\end{lem}
Based on this this lemma and continuity, it suffices to prove the
bound~\eqref{EqnREofW} for pairs of matrices \mbox{$\plu, \plv \in
  \textup{conv}\{\Balltwozero(\kdim) \cap \Ballfro(3)\}$.}  Any such
pair can be written as a weighted combination of the form $\plu =
\sum_{i} \alpha_{i} \plu_{i}$ and $\plv = \sum_{j}\beta_{j} \plv_{j}$,
with weights $\alpha_{i}, \beta_{j} \ge 0$ such that
$\sum_{i}\alpha_{i} = \sum_{j}\beta_{j} = 1$, and constituent matrices
$\plu_{i}, \plv_{j} \in \Balltwozero(\kdim) \cap \Ballfro(3)$ for each
$i,j$.	With this notation, observe that
\begin{align*}
\abs{ \trinprod{\noisemat}{\plu\otimes\plv}} \le
\sum_{i,j}\alpha_{i}\beta_{j}\abs{\trinprod{\noisemat}{\plu_{i}\otimes\plv_{j}}}
\le \Big( \sum_{i,j}\alpha_{i}\beta_{j}
\Big)\max_{i,j}\abs{\trinprod{\noisemat}{\plu_{i}\otimes\plv_{j}}} =
\max_{i,j}\abs{\trinprod{\noisemat}{\plu_{i}\otimes\plv_{j}}}.
\end{align*}
If we use $(\plu_{i})_{\cdot\ell}$ and $(\plv_{j})_{\cdot\ell}$
to denote the $\ell$-th column of $\plu_{i}$ and $\plv_{j}$, respectively,
then
\begin{align*}
\abs{\trinprod{\noisemat}{\plu_{i}\otimes\plv_{j}}} \le \sum_{\ell=
  1}^{\rdim}
\abs{\trinprod{\noisemat}{(\plu_{i})_{\cdot\ell}\otimes(\plv_{j})_{\cdot\ell}}}
& \le \Big( \sup_{\plx,\ply \in \Ball_{0}(\kdim)}
\frac{\abs{\plx^{\top}\noisemat\ply}}{\twonorm{\plx}\twonorm{\ply}}
\Big) \sum_{\ell= 1}^{\rdim} \twonorm{(\plu_{i})_{\cdot\ell}}
\twonorm{(\plv_{j})_{\cdot\ell}} \\
& \overset{(i)}{\le} 9 \Big( \sup_{\plx,\plz \in \Ball_{0}(\kdim)}
\frac{\abs{\plx^{\top}\noisemat\ply}}{\twonorm{\plx}\twonorm{\ply}}
\Big) = 9 \Big( \sup_{\plx,\ply \in \Ball_{0}(\kdim) \cap
  \Ball_{2}(1)} \abs{\plx^{\top}\noisemat\ply} \Big),
\end{align*}
where step (i) follows from the Cauchy-Schwarz inequality, and
$\fronorm{\plu_{i}} = \fronorm{\plv_{j}} \le3$.  It suffices to bound
the supremum in the last RHS by $18t$, where
\begin{align}
\label{EqnReofW_Tval}
t \defn c' \, (\snr + 1) \max \Big\{\sqrt{\frac{\kdim
    \log\usedim}{\numobs}}, \frac{\kdim\log\usedim}{\numobs} \Big\}
\end{align}
for a universal constant $c'$ to be specified later.

To proceed, we make use of a standard concentration result. Recall
that $X \in \real^{\numobs\times \usedim}$ is the matrix of
independent samples from a $\usedim$-dimensional Gaussian distribution
with zero mean and covariance $\CovMat= \snr (\OUTER{\thetastar}) +
\I_{\usedim}$.	By Lemma 15 in Loh and
Wainwright~\cite{loh2012nonconvex}, there is a universal constant $c >
0$ such that
\begin{align*}
\mprob \Big[\sup_{\plz \in \Ball_{0}(2\kdim) \cap
    \Ball_{2}(1)}\abs{\twonorm{X\plz}^{2}/\numobs -
    \plz^{\top}\CovMat\plz} \ge t \Big]\le 2\exp \Big( -cn\min
\big\{\frac{t^{2}}{(\snr + 1)^{2}},\frac{t}{\snr + 1} \big\} +
2\kdim\log\usedim \Big).
\end{align*}
Applying this inequality with $\plz= \frac{1}{6}(\plx\pm\ply)$ and our
previously specified~\eqref{EqnReofW_Tval} of $t$ with $c'=
\frac{8}{c}$, we find that with probability $1-2\usedim^{-4}$, we have
\begin{align*}
\frac{1}{36}\abs{\frac{1}{\numobs}\twonorm{X(\plx + \ply)}^{2}-(\plx +
  \ply)^{\top}\CovMat(\plx + \ply)} & \le t \quad \;\text{and} \quad
\;\frac{1}{36}\abs{\frac{1}{\numobs}\twonorm{X(\plx - \ply)}^{2}-(\plx
  - \ply)^{\top}\CovMat(\plx - \ply)} \le t
\end{align*}
for all $\plx,\ply \in \Ball_{0}(\kdim) \cap \Ball_{2}(1)$. On this event,
we have
\begin{align*}
\frac{4}{\numobs}\plx^{\top}X^{\top}X\ply & =
\frac{1}{\numobs}\twonorm{X(\plx + \ply)}^{2} -
\frac{1}{\numobs}\twonorm{X(\plx - \ply)}^{2} \\ & \le (\plx +
\ply)^{\top}\CovMat(\plx + \ply)-(\plx - \ply)^{\top}\CovMat(\plx -
\ply) + 72t \\ & = 4\plx^{\top}\CovMat\ply + 72t
\end{align*}
and similarly
$\frac{4}{n}\plx^{\top}X^{\top}X\ply\ge4\plx\CovMat\ply-72t$, whence
$\abs{\plx^{\top}\noisemat\ply} =
\abs{\frac{1}{\numobs}\plx^{\top}X^{\top}X\ply -
  \plx^{\top}\CovMat\ply} \le18t$.

\subsection{Proof of Lemma~\ref{LemGeometric}}
 \label{SecProofLemGEometric}

Let $A,B \subseteq
\real^{\usedim\times \rdim}$ be closed convex sets, with support
function given by $\phi_{A}(\plu)= \sup_{\pltheta\in A}
\trinprod{\pltheta}{\plu}$ and $\phi_{B}$ similarly defined. It is
well-known that $\phi_{A}(\plu) \leq \phi_{B}(\plu)$ for all $\plu \in
\real^{\usedim\times \rdim}$ if and only if $A \subseteq B$.
Accordingly, let us verify the first condition for the sets $A=
\Balltwoone(\sqrt{\kdim}) \cap \Ballfro(1)$ and $B= 3\text{cl}
\big\{\text{conv}\{\Balltwozero(\kdim) \cap \Ballfro(1)\} \big\}$.

For any $\plu \in \real^{\usedim\times \rdim}$, let $S \subseteq \{1,2,\ldots,\usedim\}$
be the subset that indexes the top $\lfloor\kdim\rfloor$ rows of
$\plu$ in $\ell_{2}$ norm. Then $\twoinfnorm{\plu_{S^{c}}} \le \twonorm{\plu_{j\cdot}}$
for all $j\in S$, whence
\begin{equation}
\twoinfnorm{\plu_{S^{c}}} \le \frac{1}{\lfloor\kdim\rfloor}\twoonenorm{\plu_{S}} \le \frac{1}{\sqrt{\lfloor\kdim\rfloor}} \fronorm{\plu_{S}}. \label{EqnGeoInequality}
\end{equation}
Therefore, we obtain
\begin{align*}
\phi_{A}(\plu)= \sup_{\pltheta\in A}\trinprod{\pltheta_{S}}{\plu_{S}}
+ \trinprod{\pltheta_{S^{c}}}{\plu_{S^{c}}} & \le
\sup_{\fronorm{\pltheta_{S}} \le1}\trinprod{\pltheta_{S}}{\plu_{S}} +
\sup_{\twoonenorm{\pltheta_{S^{c}}} \le
  \sqrt{\kdim}}\trinprod{\pltheta_{S^{c}}}{\plu_{S^{c}}} \\
& \leq \fronorm{\plu_{S}} + \sqrt{\kdim}\twoinfnorm{\plu_{S^{c}}} \\
& \overset{(i)}{\le} \Big( 1 +
\sqrt{\frac{\kdim}{\lfloor\kdim\rfloor}} \Big) \fronorm{\plu_{S}} \leq
3 \fronorm{\plu_{S}},
\end{align*}
where inequality (i) follows from the earlier
bound~\eqref{EqnGeoInequality}.  The claim then follows from the
observation that $\phi_{B}(\plu)= \sup \limits_{\pltheta\in
  B}\trinprod{\pltheta}{\plu} = 3 \max \limits_{\abs T=
  \lfloor\kdim\rfloor}\sup_{\fronorm{\pltheta_{T}}
  \le1}\trinprod{\pltheta_{T}}{\plu_{T}} = 3\fronorm{\plu_{S}}$.


\section{Technical lemmas for Corollary~\ref{CorOB}}
In this appendix, we provide the proofs of the technical lemmas required
for Corollary~\ref{CorOB} on one-bit matrix completion.

\subsection{Proof of Lemma~\ref{LemOBPeel}}
 \label{SecProofLemOBPeel}

In order to prove the lemma, we need to establish an upper tail bound
on the random variable
\begin{align*}
Z \defn \sup_{\d \in \Gamma(B)} \frac{\trinprod{
    \nabla \EmpLoss(\OUTER{\pltheta}) - \PopGrad(\OUTER{\pltheta})}
  {\d\otimes\pltheta} }{\fronorm{\d}}.
\end{align*}
Define the indicator variable $\Obsind_{ij} \defn \mathbb{I}\{(i,j)
\in \Obs\}$ for each $(i,j)$. Using the
expression~\eqref{EqnFactorGrad} for $\nabla \EmpLoss(\cdot)$ and the definition
of $\Gamma(B)$, we observe that
\begin{align}
\label{eqOBDeviationTerm1}
Z\le & \frac{2}{B\noisestd}\sup_{\d \in \Gamma(B)}\sum_{i,j}
\left[\Obsind_{ij}h_{ij}(\PlTheta_{ij}/\noisestd) -
  \Exs_{\Obsind,\Yout}\Obsind_{ij}h_{ij}(\PlTheta_{ij}/\noisestd)\right]\cdot(\d\otimes\theta)_{ij}.
\end{align}
By the usual symmetrization argument~\cite{LedTal91}, the expectation
$\Exs_{\Obsind, \Yout}[Z]$ is at most a factor of two times the
Rademacher-symmetrized version.  This is the supremum of a
sub-Gaussian process in terms of Rademacher variables, and so is
majorized by the expected supremum of the corresponding Gaussian
process~\cite{LedTal91} up to a universal constant $c$.  In
conjunction, these two steps yield the bound
\begin{align*}
\Exs_{\Obsind,\Yout}Z \leq c \Exs_{\Obsind,\Yout}\Exs_{g}\sup_{\d \in
  \Gamma(B)} \left\{ Z(\d)\right\} \defn & c
\Exs_{\Obsind,\Yout}\Exs_{g}\sup_{\d \in \Gamma(B)} \left\{ \frac{4}{B
  \noisestd} \sum_{i,j} g_{ij} \Obsind_{ij} h_{ij} (
\PlTheta_{ij}/\noisestd) \cdot (\d \otimes \theta)_{ij}\right\} ,
\end{align*}
where $\{g_{ij}\}$ are independent standard Gaussian variables. 

Our next step is to bound $\Exs_{g}\sup_{\d \in \Gamma(B)}Z(\d)$ using
the Sudakov-Fernique comparison inequality. For any $\d,\d' \in
\Gamma(B)$ and $\pltheta' \defn \thetastar + \d'$, $\PlTheta' \defn
\OUTER{\pltheta'}$, we have
\begin{align*}
\gamma(\d,\d') & \defn \Exs_{g}(Z(\d)-Z(\d'))^{2} \\
& = \frac{16}{(B\noisestd)^{2}}\Exs_{g}
\left[\sum_{i,j}g_{ij}\Obsind_{ij}
  \left(h_{ij}(\PlTheta_{ij}/\noisestd)\cdot(\d\otimes\pltheta)_{ij}-h_{ij}(\PlTheta_{ij}'/\noisestd)\cdot(\d'\otimes\pltheta')_{ij}\right)\right]^{2}
\\ 
 & = \frac{16}{(B\noisestd)^{2}}\sum_{i,j}\Obsind_{ij}^{2} \left\{
h_{ij}(\PlTheta_{ij}/\noisestd)\cdot(\d\otimes\pltheta)_{ij}-h_{ij}(\PlTheta_{ij}'/\noisestd)\cdot(\d'\otimes\pltheta')_{ij}\right\}
^{2} \\
& \leq \frac{32}{(B\noisestd)^{2}}\sum_{i,j}\Obsind_{ij}^{2} \left\{
h_{ij}(\Theta_{ij}/\noisestd)^{2}\cdot(\d\otimes\pltheta -
\d'\otimes\pltheta')_{ij}^{2} +
\left(h_{ij}(\Theta_{ij}/\noisestd)-h_{ij}(\PlTheta'_{ij}/\noisestd)\right)^{2}\cdot(\d'\otimes\pltheta')_{ij}^{2}\right\}.
\end{align*}
Recall that for each $(i,j)$ and over the interval
$[-4\obsnr,4\obsnr]$, the function $h_{ij}(\cdot)$ is surely bounded
by $4\obupper{4\obsnr}$ and $4\obupper{4\obsnr}$-Lipschitz. Moreover,
the Cauchy-Schwarz inequality implies that
\begin{align*}
\abs{\Theta_{ij}} = \abs{(\OUTER{\pltheta})_{ij}} \le
\twonorm{\pltheta{}_{i\cdot}}\twonorm{\pltheta{}_{j\cdot}} \le
2\sqrt{\frac{\inco\rdim}{\usedim}}\cdot2\sqrt{\frac{\inco\rdim}{\usedim}}
= 4\noisestd\obsnr.
\end{align*}
Note that the same bound holds for $\abs{\PlTheta'_{ij}}$,
$\abs{(\d\otimes\pltheta)_{ij}}$ and
$\abs{(\d'\otimes\pltheta')_{ij}}$.  It follows that
\begin{align*}
\gamma(\d,\d') & \le C^{2}
\frac{\obupper{4\obsnr}^{2}}{(B\noisestd)^{2}}\sum_{i,j}\Obsind_{ij}^{2}
\left\{ (\d\otimes\pltheta - \d'\otimes\pltheta')_{ij}^{2} +
(\PlTheta_{ij}/\noisestd -
\PlTheta_{ij}'/\noisestd)^{2}\cdot4\noisestd^{2}\obsnr^{2}\right\}
\\ 
& = C^{2}
\frac{\obupper{4\obsnr}^{2}}{(B\noisestd)^{2}}\sum_{i,j}\Obsind_{ij}^{2}
\left\{ (\d\otimes\pltheta - \d'\otimes\pltheta')_{ij}^{2} +
\left(2\d\otimes\pltheta_{ij} +
\d\otimes\d_{ij}-2\d'\otimes\pltheta'_{ij} -
\d'\otimes\d'_{ij}\right)^{2}4\obsnr^{2}\right\} \\ 
& \le C^{2}(1 + \obsnr)^{2}
\frac{\obupper{4\obsnr}^{2}}{(B\noisestd)^{2}}\sum_{i,j}\Obsind_{ij}^{2}
\left\{ (\d\otimes\pltheta - \d'\otimes\pltheta')_{ij}^{2} +
(\d\otimes\d - \d'\otimes\d')_{ij}^{2}\right\} .
\end{align*}
We compare $Z(\d)$ with an alternative stochastic process given by
\begin{align*}
\bar{Z}(\d) \defn C(1 + \obsnr)
\frac{\obupper{4\obsnr}}{B\noisestd}\sum_{i,j}
\left[g_{ij}\Obsind_{ij}(\OUTER{\d})_{ij} +
  g'_{ij}\Obsind_{ij}(\d\otimes\pltheta)_{ij}\right],
\end{align*}
where $\{g_{ij},g_{ij}'\}$ are independent standard Gaussian
variables.  Both $Z(\d)$ and $\bar{Z}(\d)$ are surely continuous
in~$\d$. Observe that by independence, we have
\begin{align*}
\bar{\gamma}(\d,\d'): & = \Exs_{g,g'}(\bar{Z}(\d) - \bar{Z}(\d'))^{2}
\\
& = C^{2}(1 + \obsnr)^{2}
\frac{\obupper{4\obsnr}^{2}}{(B\noisestd)^{2}}\sum_{i,j}\Obsind_{ij}^{2}
\left\{ (\OUTER{\d} - \OUTER{\d'})_{ij}^{2} + (\d\otimes\pltheta -
\d'\otimes\pltheta')_{ij}^{2}\right\} ,
\end{align*}
so we have $\gamma(\d,\d')\le \bar{\gamma}(\d,\d'),\forall\d,\d' \in
\Gamma(B)$.  By the Sudakov-Fernique comparison~\cite{LedTal91}, we
find that
\begin{align*}
\Exs_{\Obsind,\Yout}\Exs_{g}\sup_{\d \in \Gamma(B)}Z(\d) & \le
\Exs_{\Obsind,\Yout}\Exs_{g,g'}\sup_{\d \in \Gamma(B)}\bar{Z}(\d) \\ &
= C(1 + \obsnr)
\frac{\obupper{4\obsnr}}{B\noisestd}\cdot\Exs_{\Obsind,g,g'}\sup_{\d
  \in \Gamma(B)}\sum_{i,j} \left[g_{ij}\Obsind_{ij}(\OUTER{\d})_{ij} +
  g'_{ij}\Obsind_{ij}(\d\otimes\pltheta)_{ij}\right] \\ & \le C(1 +
\obsnr) \frac{\obupper{4\obsnr}}{B\noisestd} \left(\sup_{\d \in
  \Gamma(B)}\nucnorm{\OUTER{\d}} + \sup_{\d \in
  \Gamma(B)}\nucnorm{\d\otimes\pltheta}\right)\Exs_{\Obsind,g}
\opnorm{g\circ\Obsind},
\end{align*}
where the last inequality follows from the generalized Holder's
inequality and that $g$ and $g'$ are identically distributed. To
proceed, we use Lemma~\ref{LemCensorGaussMat} to get that
$\Exs_{\Obsind,g} \opnorm{g\circ\Obsind} \le c(\sqrt{\pobs\usedim} +
\log\usedim)\le 2c\sqrt{\pobs\usedim}$, where the last inequality
follows from the assumption $\pobs\ge
\frac{\log^{2}\usedim}{\usedim}$.  Moreover, for each $\d \in
\Gamma(B)$, the matrices $\OUTER{\d}$ and $\d\otimes\pltheta$ have
rank at most $\rdim$, so
\begin{align*}
\max\{\nucnorm{\OUTER{\d}},\nucnorm{\d\otimes\pltheta}\} \le
\sqrt{\rdim}\cdot\fronorm{\d}\max\{\opnorm{\d},\opnorm{\pltheta}\} \le
2\sqrt{\rdim}B.
\end{align*}
Putting together the pieces yields
\begin{align*}
\Exs_{\Obsind,\Yout}Z & \le \Exs_{\Obsind,\Yout}\Exs_{g}\sup_{\d \in
  \Gamma(B)}Z(\d)\le \frac{4C'}{\noisestd}\obupper{4\obsnr}(1 +
\obsnr)\sqrt{\pobs\usedim\rdim} \le \frac{1}{8}\alpha\epsnum.
\end{align*}

In order establish concentration of $Z$ around
$\Exs_{\Obsind,\Yout}Z$, we use a standard functional Hoeffding
inequality~\cite{Ledoux01}.  In particular, letting
$\{X_i\}_{i=1}^\numobs$ be independent random variables such that
$X_i$ takes values in $\mathcal{X}_i$, consider a random variable of
the form $Y \defn \sup \limits_{g \in \mathcal{G}} \sum_{i=1}^\numobs
g(X_i)$ where for each $g \in \mathcal{G}$, we have $\sup_{x \in
  \mathcal{X}_i} |g(x)| \leq b_i$.  Then we are guaranteed that
\begin{align}
\label{EqnFunctionalHoeffding}
\mprob \left[Y \geq \Exs[Y] + \tau\right]\le e^{ - \tau^{2}/16D^{2}}
\qquad \mbox{for all $\tau \geq 0$, and $D^{2} \defn
  \sum^\numobs_{i=1} b_{i}^{2}$.  }
\end{align}
Setting $\tau= \frac{1}{8}\alpha\epsnum$, we have
\begin{align*}
\sup_{\d \in \Gamma(B)}\sum_{i,j} \frac{4[\nabla \EmpLoss(\PlTheta) -
    \PopGrad(\PlTheta)]_{ij}^{2}(\d\otimes\pltheta)_{ij}^{2}}{\fronorm{\d}^{2}}
& \le C\frac{\obupper{4\obsnr}^{2}}{B^{2}\noisestd^{2}}\sup_{\d \in
  \Gamma(B)} \fronorm{\d\otimes\pltheta}^{2} \\ 
& \le C\frac{\obupper{4\obsnr}^{2}}{\noisestd^{2}B^{2}}\cdot\sup_{\d
  \in \Gamma(B)} \fronorm{\d}^{2} \opnorm{\pltheta}^{2} \le
\frac{\alpha^{2}\epsnum^{2}}{128\times14\log\usedim}\to D^{2}.
\end{align*}
Consequently, applying the bound~\eqref{EqnFunctionalHoeffding} with
these choices of $(\tau,D^{2})$,
\footnote{In particular, we apply it with the following setup: $X_{ij}
  = (e_{i}\otimes e_{j},\Obsind_{ij},\Yout_{ij})$ and
  $\mathcal{X}_{ij} = \{e_{i}\otimes e_{j}\}\times \{0,1\}\times
  \{-1,1\}$, where $e_{i}$ is the $i$-th standard basis vector in
  $\real^{\usedim}$; $\mathcal{F} = \{\zeta_{\d}:\d \in \Gamma(B)\}$
  with
\[
\zeta_{\d}(X_{ij})=
\frac{2}{\noisestd\fronorm{\d}}\trinprod{\frac{f'(\PlTheta/\noisestd)\left(\Yout_{ij}(e_{i}\otimes
    e_{j})-2f(\PlTheta/\noisestd) +
    1\right)}{f(\PlTheta/\noisestd)\left(1-f(\PlTheta/\noisestd)\right)}-p\frac{f'(\PlTheta/\noisestd)\left(2f(\ThetaStar/\noisestd)-2f(\PlTheta/\noisestd)\right)}{f(\PlTheta/\noisestd)\left(1-f(\PlTheta/\noisestd)\right)}}{\Obsind_{ij}(e_{i}\otimes
  e_{j})\circ(\d\otimes\pltheta)}.
\]
} we obtain $\mprob [ Z \geq \Exs_{\Obsind,\Yout}[Z] +
  \frac{1}{8}\curv \epsnum ] \le \usedim^{-14}$.  Combining with the
the expectation bound $\Exs_{\Obsind,\Yout} Z \le \frac{1}{8}
\curv \epsnum$, we find that
\begin{align*}
\mprob \Big[ Z = \sup_{\d \in \Gamma(B)} \trinprod{
    \nabla \EmpLoss(\OUTER{\pltheta}) - \PopGrad(\OUTER{\pltheta})}{\d} /
  \fronorm{\d} \ge \frac{1}{4} \curv \epsnum \Big] \le \usedim^{-14}.
\end{align*}
Following the same lines of argument we obtain a similar bound on the
lower tail:
\begin{align*}
\mprob \Big[ \inf_{\d \in \Gamma(B)} \trinprod{
    \nabla \EmpLoss(\OUTER{\pltheta}) - \PopGrad(\OUTER{\pltheta})} {\d} /
  \fronorm{\d} \le - \frac{1}{4} \curv \epsnum \Big] \le \usedim^{-14}.
\end{align*}
The proof the lemma is completed by applying the union bound.


\subsection{Proof of Lemma~\ref{LemOBSmallRad}}
\label{SecProofLemOBSmallRad}

Using the Cauchy-Schwarz inequality and the
expression~\eqref{EqnFactorGrad} for the gradient $\nabla \EmpLoss$,
we have
\begin{align*}
\sup_{\d \in \Gamma_{0}} \frac{ \abs{ \trinprod{ \nabla
      \EmpLoss(\OUTER{\pltheta}) - \PopGrad(\OUTER{\pltheta})}
    {\d\otimes\pltheta} } }{\fronorm{\d}} 
& \le \frac{1}{\noisestd}\sup_{\d \in \Gamma_{0}}
\fronorm{\left[\ProjObs h(\PlTheta/\noisestd) - \Exs\ProjObs
    h(\PlTheta/\noisestd)\right]\pltheta}  \leq \sum_{j=1}^3 T_3
\end{align*}
where $T_1 \defn
\frac{2}{\noisestd}\sup_{\d \in \Gamma_{0}}
\fronorm{\left[\ProjObs h(\PlTheta/\noisestd) - \ProjObs
    h(\ThetaStar/\noisestd)\right]\pltheta}$, and
\begin{align*}
T_2 \defn \frac{2}{\noisestd}\sup_{\d \in \Gamma_{0}}
\fronorm{\left[\ProjObs h(\ThetaStar/\noisestd)\right]\pltheta}, \quad
\mbox{and} \quad T_3 \defn \frac{2}{\noisestd}\sup_{\d \in \Gamma_{0}}
\fronorm{[\Exs\ProjObs h(\PlTheta/\noisestd)]\pltheta}.
\end{align*}
To prove the lemma, it suffices to show that with probability at least
$1 - \usedim^{-12}$, each of $\{\Term_{1},\Term_{2},\Term_{3}\}$ is
bounded from above by $\frac{1}{12}\alpha\epsnum$. For $\Term_{1}$, we
have
\begin{align*}
\Term_{1} \le \frac{2}{\noisestd}\sup_{\d \in \Gamma_{0}} \usedim
\infnorm{H(\PlTheta/\noisestd) - H(\ThetaStar/\noisestd)} \cdot
\opnorm{\pltheta} & \overset{(i)}{\le}
\frac{16\usedim}{\noisestd}\cdot\obupper{4\obsnr}\cdot\sup_{\d \in
  \Gamma_{0}}\infnorm{\PlTheta/\noisestd - \ThetaStar/\noisestd} \\
& \leq \frac{16d2^{ -
    \usedim^{2}}\obupper{4\obsnr}}{\noisestd^{2}}\overset{(ii)}{\le}
\frac{\sqrt{\rdim\log^{2}\usedim}\obupper{4\obsnr}}{\noisestd}\obsnr\le
\frac{1}{12}\alpha\epsnum,
\end{align*}
where the step (i) follows from the fact that $h$ is element-wise
$4\obupper{4\obsnr}$-Lipschitz over $[-4\obsnr,4\obsnr]$ and that
$\opnorm{\pltheta} \le 2$ for $\d \in \Gamma_{0}$, and in step (ii)
from the definition $\obsnr \defn
\frac{\inco\rdim}{\usedim\noisestd}.$ Since $\fronorm{\pltheta} \le
2\sqrt{\rdim}$ for $\d \in \Gamma_{0}$, we have
\begin{align*}
\Term_{2} \le \frac{2}{\noisestd} \opnorm{\ProjObs
  h(\ThetaStar/\noisestd)}\sup_{\d \in \Gamma_{0}} \fronorm{\pltheta}
\le \frac{4\sqrt{\rdim}}{\noisestd} \opnorm{\ProjObs
  h(\ThetaStar/\noisestd)}
\end{align*}
Note that for each index pair $(i,j)$,
$\Exs_{\Yout}h_{ij}(\ThetaStar_{ij}/\noisestd)= 0$, and that
$\abs{h_{ij}(\ThetaStar_{ij}/\noisestd)} \le4\obupper{4\obsnr}$ since
$\infnorm{\ThetaStar/\sigma} \le4\obsnr$. Therefore, the matrix
$\frac{1}{4\obupper{4\obsnr}}\ProjObs h(\ThetaStar/\noisestd)$ is a
censored sub-Gaussian random matrix satisfying the assumptions in
Lemma~\ref{LemCensorGaussMat}, by which we obtain $\opnorm{\ProjObs
  h(\ThetaStar/\noisestd)} \le C\obupper{4\obsnr}\sqrt{\pobs\usedim}$
with probability at least $1 - \usedim^{-12}$. It follows that with
the same probability, the second term is bounded as $\Term_{2} \le
\frac{4C\obupper{4\obsnr}\sqrt{\pobs\usedim\rdim}}{\noisestd} \le
\frac{1}{12}\alpha\epsnum$.  Finally, the third term $\Term_{3}$ can
be bounded as
\begin{align*}
\Term_{3} & \le \frac{2}{\noisestd}\sup_{\d \in \Gamma_{0}}
\fronorm{[\Exs\ProjObs h(\PlTheta/\noisestd)]} \opnorm{\pltheta} =
\frac{2}{\noisestd}\sup_{\d \in \Gamma_{0}}
\fronorm{\pobs\cdot\frac{2f'(\PlTheta/\noisestd)(f(\ThetaStar/\noisestd)-f(\PlTheta/\noisestd))}{f(\PlTheta/\noisestd)(1-f(\PlTheta/\noisestd))}}
\opnorm{\pltheta}.
\end{align*}
Note that $\opnorm{\pltheta} \le 2$ for $\d \in \Gamma_{0}$. Moreover,
because $f$ satisfies~(\ref{EqnFlatness}), we know that
$\abs{\frac{f'(x)}{f(x)(1-f(x))}} \le \sqrt{\obupper{4\obsnr}}$ and
$\abs{f(x)-f(x')} \le \sqrt{\obupper{4\obsnr}}\abs{x-x'}$ for all
$x,x'\in[-4\obsnr,4\obsnr]$. It follows that
\begin{align*}
\Term_{3} \leq \frac{8 \pobs \obupper{4\obsnr}}{\noisestd^{2}}
\sup_{\d \in \Gamma_{0}} \fronorm{\ThetaStar - \PlTheta} \leq
\frac{8\pobs\obupper{4\obsnr}}{\noisestd^{2}} \cdot
\frac{3}{2^{\usedim^{2}}} \overset{(i)}{\leq} \frac{24 \sqrt{\pobs
    \usedim \rdim} \obupper{4\obsnr}\nu}{\noisestd} \leq \frac{1}{12}
\alpha \epsnum,
\end{align*}
where the step (i) follows from the definition $\obsnr \defn
\frac{\inco\rdim}{\usedim\noisestd}$. This completes the proof of the
lemma.


\section{Proof of inequality~\eqref{EqnLSini}}
\label{AppProofLemLSini}

Recalling that the matrix $\thetastar$ is orthonormal and
$\inco$-incoherent, we have $\infnorm{\OUTER{\thetastar}} \le
\frac{\inco\rdim}{\usedim}$, and hence $\infnorm{\widebar{\Yout} -
  \OUTER{\thetastar}} \leq 2 \frac{\inco\rdim}{\usedim}$.  On the
other hand, we claim that each row and column of the matrix
$\widebar{\Yout} - \OUTER{\thetastar}$ has at most $\kdim$ non-zero
elements.  To see this, let $\StarSupp$ be the set of the non-zero
element of $\SStar$. If $(i,j)\not \in \StarSupp$, then
$\abs{\Yout_{ij}} = \abs{(\OUTER{\thetastar})_{ij}} \le
\frac{\inco\rdim}{\usedim}$, so $\widebar{\Yout}_{ij} = \Yout_{ij} =
(\OUTER{\thetastar})_{ij}$ and thus $(\widebar{\Yout} -
\OUTER{\thetastar})_{ij} = 0$. Therefore, we find that
$\widebar{\Yout} - \OUTER{\thetastar}$ is supported on the elements in
$\StarSupp$, hence the claim. With the above two facts, we apply
Proposition~3 in the paper~\cite{chandrasekaran2011siam} to obtain
that
\begin{align*}
\opnorm{\widebar{\Yout} - \OUTER{\thetastar}} \le
\kdim\infnorm{\widebar{\Yout} - \OUTER{\thetastar}} \le
2\frac{\inco\rdim\kdim}{\usedim}.
\end{align*}
On the other hand, the gap between the $\rdim$-th and $(\rdim + 1)$-th
singular values of the matrix $\OUTER{\thetastar}$ is $1$. Letting $ \widebar{U} $ be the matrix of the top-$\rdim $ singular vectors of $ \widebar{\Yout} $ and using  $ \Theta[\cdot, \cdot] $ to denote the principal angles between two subspaces, we find that 
\begin{align*}
\min_{\thetastar \in \eclass} \opnorm{\widebar{U} - \thetastar} & \le
\sqrt{2} \opnorm{\sin\Theta \big[
    \text{col}(\widebar{U}),\text{col}(\thetastar) \big] } \le
2\opnorm{\widebar{\Yout} - \OUTER{\thetastar}} \le
\frac{4\inco\rdim\kdim}{\usedim},
\end{align*}
where the first step follows from Proposition 2.2 in the
paper~\cite{vu2013minimax} and the second step follows from Wedin's
$\sin\Theta$ theorem~\cite{golub_matrixcomp}.  It follows that
\begin{align*}
\dist{\thetait{0}} \le \dist{\widebar{U}} \le \sqrt{\rdim}
\min_{\thetastar \in \eclass} \opnorm{\widebar{U} - \thetastar} \le
\frac{4\inco\rdim\sqrt{\rdim}\kdim}{\usedim},
\end{align*}
where the first step holds because $ \thetait{0} =
\ProjGenCons(\widebar{U}) $ and projection onto the convex set $
\GenCons $ is non-expansive, which completes the proof of the claim.

\section{Spectral norms of censored sub-Gaussian random matrices}

In this appendix, we state and prove a useful bound on the spectral
norm sub-Gaussian random matrices with censored entries.
\begin{lem}
\label{LemCensorGaussMat}
Suppose $X \in \real^{\usedim\times \usedim}$ is a symmetric random
matrix with $X_{ij} = g_{ij}\Obsind_{ij}$, where $\{g_{ij} \mid i \ge
j\}$ are independent zero-mean sub-Gaussian random variables with
parameter $1$, $\{\Obsind_{ij} \mid i \ge j \}$ are independent
Bernoulli variables with parameter $\pobs$, and they are mutually
independent. Then there exists a universal constant $c > 0$ such that
\begin{subequations}
\begin{align}
\label{EqnCensorGaussMatExp}
\Exs \big[ \opnorm{X} \big] & \le c \big( \sqrt{\pobs\usedim} +
\log\usedim \big), \quad \mbox{and} \\
\label{EqnCensorGaussMatTail}
\mprob \Big[\opnorm X \geq c \big( \sqrt{\pobs\usedim} + \log\usedim
  \big) \Big] & \leq d^{-12}.
\end{align}
\end{subequations}
\end{lem}

Let us now prove this lemma.  By a standard symmetrization argument,
we can assume without loss of generality that each $g_{ij}$ is a
symmetric random variable. To proceed, we need the following result
from Bandeira and van Handel~\cite{bandeira2014randmat}:
\begin{prop}
[Corollaries 3.6 and 3.12 in
  \cite{bandeira2014randmat}]\label{PropBandeira}Let $\tilde{X}$ be
the $\usedim\times \usedim$ symmetric random matrix whose entries
$\tilde{X}_{ij}$ are independent symmetric random variables bounded by
$\tilde{\sigma}_{*}$, and define \mbox{$\tilde{\sigma} \defn \max_{i}
  \sqrt{\sum_{j}\Exs [\tilde{X}_{ij}^{2}]}$.}  Then there exist
universal constants $\tilde{c}_{1}$ and $\tilde{c}_{2}$ such that
\begin{subequations}
\begin{align}
\Exs \opnorm{\tilde{X}} & \leq 3\tilde{\sigma} + \tilde{c}_{1}
\tilde{\sigma}_{*} \sqrt{\log\usedim}, \quad \mbox{and} \\
\mprob \big[ \opnorm{\tilde{X}} \ge3\tilde{\sigma} + t \big] & \le
\usedim \cdot \exp \Big( - \frac{t^{2}}{\tilde{c}_{2}
  \tilde{\sigma}_{*}^{2}} \Big) \qquad \mbox{for each $t \geq 0$.}
\end{align}
\end{subequations}
\end{prop}
To apply the proposition with unbounded entries in $X$, we use a
standard truncation argument. For some constant $b$ to be specified
later, let $\tilde{X}$ be the matrix with $\tilde{X}_{ij} =
X_{ij}1_{\abs{X_{ij}} \le b\sqrt{\log\usedim}}.$ Observe that
$\tilde{X}$ satisfies the assumption in Proposition \ref{PropBandeira}
with $\tilde{\sigma}_{*} \le b\sqrt{\log\usedim}$ and $\tilde{\sigma}
\le \sqrt{\pobs\usedim}$. Applying Proposition \ref{PropBandeira} with
$t= \sqrt{12\tilde{c}_{2}b^{2}}\log\usedim$, we obtain the bounds
\begin{subequations}
\begin{align}
\label{EqnOPExpBounded}
\Exs \opnorm{\tilde{X}} & \le 3\sqrt{\pobs\usedim} + \tilde{c}_{1}
b \log \usedim, \quad \mbox{and} \\
\label{EqnOPConcentrateBounded}
\mprob \big[\opnorm{\tilde{X}} \ge 3\sqrt{\pobs\usedim} +
  \sqrt{12\tilde{c}_{2}b^{2}} \log\usedim \big] & \le \usedim \exp
\Big( - \frac{t^{2}}{\tilde{c}_{2} \tilde{\sigma}_{*}^{2}} \Big) \le
\usedim^{-13},
\end{align}
\end{subequations}
where the last inequality follows from $t\ge
\tilde{\sigma}_{*}\sqrt{12\tilde{c}\log\usedim}$.  On the other hand,
by choosing the constant $b$ sufficiently large and using a standard
bound on the maximum of sub-Gaussian variables, we know that
\begin{align*}
\mprob \big[ \tilde{X} \neq X \big] \le \mprob \Big[ \max_{i,j}
  \abs{g_{ij}} > b\sqrt{\log\usedim} \Big] \le \usedim^{-13}.
\end{align*}
Combining with the tail bound~\eqref{EqnOPConcentrateBounded} yields
\begin{align*}
\mprob \left[\opnorm X\ge3\sqrt{\pobs\usedim} +
  \sqrt{12\tilde{c}b^{2}}\log\usedim\right]\le \mprob
\left[X\neq\tilde{X}\right] + \mprob \left[\opnorm{\tilde{X}}
  \ge3\sqrt{\pobs\usedim} +
  \sqrt{12\tilde{c}b^{2}}\log\usedim\right]\le \usedim^{-12},
\end{align*}
which proves the second inequality in Lemma~\ref{LemCensorGaussMat}.

Turning to the first inequality in the lemma, we let $\breve{X}$ be
the matrix with $\breve{X}_{ij} = X_{ij} - \tilde{X}_{ij} =
X_{ij}1_{\abs{X_{ij}} > b\sqrt{\log\usedim}}$, and observe that by
definition, $\mprob(0 < \max_{i,j}\abs{\breve{X}_{ij}} \le
b\sqrt{\log\usedim})= 0$.  Moreover, by choosing the constant $b$
sufficiently large and using a standard concentration inequality for
convex Lipschitz functions~\cite{Ledoux01}, we find that for each $t
\ge 0$,
\begin{align*}
\mprob \left[\max_{i,j}\abs{\breve{X}_{ij}} > b\sqrt{\log\usedim} +
  t\right] & \le \mprob \left[\max_{i,j} \abs{g_{ij}} \ge
  \Exs\max_{i,j}\abs{g_{ij}} + t + 4 \sqrt{\log\usedim}\right]
\\
& \le 2e^{-(t + 4\sqrt{\log\usedim})^{2}/5} \le
\frac{2}{\usedim^{2}}e^{-t^{2}/5}.
\end{align*}
Integrating these tail bounds gives $\Exs
\big[\max_{i,j}\abs{\breve{X}_{ij}} \big] \le
\frac{\breve{c}\sqrt{\log\usedim}}{\usedim^{2}}$ Combining with
equation~\eqref{EqnOPExpBounded} yields the upper bound
\begin{align*}
\Exs\opnorm X\le \Exs \opnorm{\tilde{X}} + \Exs \opnorm{\breve{X}} \le
\Exs \opnorm{\tilde{X}} + \usedim\Exs\max_{i,j}\abs{\breve{X}_{ij}}
\le3\sqrt{\pobs\usedim} + \tilde{c}_{1}b\log\usedim +
\frac{\breve{c}\sqrt{\log\usedim}}{\usedim},
\end{align*}
which completes the proof of Lemma~\ref{LemCensorGaussMat}.





\printbibliography
\end{document}